\documentclass[11pt,a4paper,final]{article}

\usepackage{amsfonts,amsmath,amssymb,theorem,graphicx,hyperref,color,natbib}
\usepackage{multirow}
\newcommand{\ind}{1\!{\rm l}}

\newcommand{\bu}{\mathbf{u}}
\newcommand{\bc}{\mathbf{c}}
\newcommand{\bw}{\mathbf{\omega}}
\newcommand{\btheta}{\boldsymbol{\theta}}

\newtheorem{Th}{Theorem}

\newtheorem{Prop}{Proposition}
\newtheorem{Lemma}{Lemma}

\newtheorem{Cor}{Corollary}
\newtheorem{remark}{Remark}

\def\1{1\!{\rm l}}

\newcommand{\E}{\ensuremath{\mathbb{E}}}

\renewcommand{\P}{\ensuremath{\mathbb{P}}}

\newcommand{\R}{\ensuremath{\mathbb{R}}}
\newcommand{\Nb}{\ensuremath{N}}
\newcommand{\Xn}{\ensuremath{\mathcal{X}^{(n)}}}

\newcommand{\Nor}{\mathrm{N}}
\renewcommand{\L}[1]{\ensuremath{\mathbb{L}^{#1}}}

\renewcommand{\L}{\ensuremath{\mathbb{L}}}

\renewcommand{\d}{\mathrm{d}}

\newcommand{\pg}[1]{\left\{#1\right\}}
\newcommand{\pq}[1]{\left[#1\right]}
\newcommand{\pt}[1]{\left(#1\right)}
\newcommand{\vir}[1]{``#1''}

\newenvironment{proof}{\noindent{\bf Proof.}}{\hfill
 $\square$\par\noindent
 }

\numberwithin{equation}{section}

\newcommand{\Data}{X^{(n)}}
\newcommand{\data}{x^{(n)}}
\newcommand{\hga}{\hat \gamma_n}

\begin{document}
\title{Posterior concentration rates for empirical Bayes procedures, with applications to Dirichlet Process mixtures}
\author{Sophie Donnet\thanks{Math\'ematiques et Informatique Appliqu\'ees (MIA). 	Institut national de la recherche agronomique (INRA). UMR0518 Ð AgroParisTech}, Vincent Rivoirard\thanks{Universit\'e de Paris Dauphine, France}, Judith Rousseau\thanks{CREST, France} and Catia Scricciolo\thanks{Bocconi University, Italy}}
\maketitle

\begin{abstract}
In this paper we provide general conditions to check on the model and the prior to derive posterior concentration rates for data-dependent priors (or empirical Bayes approaches). We aim at providing conditions that are close to the conditions provided in the seminal paper by  \citet{ghosal:vdv:07}. We then apply the general theorem to two different settings: the estimation of a density using Dirichlet process mixtures of Gaussian random variables with base measure depending on some empirical quantities and the estimation of the intensity of a counting process under the Aalen model. A simulation study for inhomogeneous Poisson processes also illustrates our results. In the former case we also derive some results on the estimation of the mixing density and on the deconvolution problem. In the latter, we provide a general theorem on posterior concentration rates for  counting processes with Aalen multiplicative intensity with priors not depending on the data.
\end{abstract}


{\bf Short title} Empirical Bayes posterior concentration rates.

\section{Introduction}\label{intro}
In a Bayesian approach to inference, the prior distribution should, in principle, be chosen
independently of the data; however, it is not always an easy task to elicit
the values of the prior hyperparameters and a common practice is to estimate
them
by
reasonable empirical quantities. The prior is then data-dependent and the approach falls under the umbrella
of empirical Bayes methods, as opposed to fully Bayesian methods.
More formally, consider a statistical model $( \mathbb P_\theta^{(n)},\, \theta \in \Theta)$ over a sample space $\mathcal X^{(n)}$, together with a family of prior distributions $\pi( \cdot \mid \gamma) $, $\gamma \in \Gamma$, on the parameter space $\Theta$.
A Bayesian statistician would either set $\gamma$ to a specific value or integrate it out using a probability distribution
in a hierarchical specification of the prior for $\theta$. Both approaches would lead to a prior distribution for $\theta$ that does not depend on the data, say $\pi$, resulting in a posterior distribution $\pi( \cdot \mid \Data)$.
However, it is often the case that
knowledge
is not \emph{a priori} available
to either fix a value for $\gamma$ or elicit a prior distribution for it,
so that these hyperparameters are more easily estimated from the data.
Throughout the paper, we will denote by $\hga$ a data-driven selection of the prior hyperparameters.
There are many instances in the literature where empirical Bayes selection of the prior hyperparameters is performed,
sometimes without explicitly mentioning it. Some examples concerning the parametric case can be found in \citet{casella:85}, \citet{berger:1985} and \citet{ahmed:reid:01}. Regarding the nonparametric case,
\citet{richardson:green:1997} propose a default empirical Bayes approach to deal with parametric and nonparametric mixtures of Gaussian random variables; \citet{mcauliffe:blei:jordan} propose another empirical Bayes approach for Dirichlet process mixtures of Gaussian random variables, while in  \citet{knapikSVZ2012} and \citet{szabo:vdv:vzanten:13} an empirical Bayes procedure is proposed in the context of the white noise model to obtain adaptive posterior distributions. There are many other instances of empirical Bayes methods in the literature, especially in applied problems. Quoting \citet{hjort:holmes:mueller:2010}, \vir{\citet{efron:2003} argues that the brightest statistical future may be reserved for \emph{empirical Bayes} methods}.

In this paper, our aim is not to claim that empirical Bayes methods are better than fully Bayesian methods in some way, but rather to provide tools to study frequentist asymptotic properties of empirical Bayes posterior distributions, given their wide use in practice.
Surprisingly, very little is known, in a general framework, on the asymptotic behavior of such empirical Bayes posterior distributions. It is common belief that, if $\hga$ asymptotically converges to some value $\gamma^\ast$, then the empirical Bayes posterior distribution associated with $\hga$ is eventually close to the Bayesian posterior associated with $\gamma^\ast$. Results have been obtained in explicit specific frameworks
in
\citet{knapikSVZ2012} and \citet{szabo:vdv:vzanten:13} for the
white noise model; in \citet{clyde:george:2000} and \citet{cui:george:2008}
for wavelets or variable selection. Recently, \citet{petrone:rousseau:scricciolo:14} have investigated asymptotic properties of empirical Bayes posterior distributions: they have obtained general conditions for frequentist consistency of empirical Bayes posteriors and, in the parametric case, studied when strong merging between Bayesian and maximum marginal likelihood empirical Bayes posterior distributions takes place.

In this work, we are interested in studying the frequentist asymptotic behavior of empirical Bayes posterior distributions $\pi( \cdot \mid \Data,\, \hga) $ in terms of contraction rates. Let $d(\cdot,\,\cdot)$ be a loss function, say a metric or a pseudo-metric, on $\Theta$ and, for $\theta_0 \in \Theta$, let $U_\epsilon$ be a neighborhood of $\theta_0$, \emph{i.e.}, a set of the form
$\{ \theta:\, d(\theta,\, \theta_0) < \epsilon \}$ with $\epsilon>0$. The empirical Bayes posterior distribution is said to concentrate at $\theta_0$, with rate $\epsilon_n$ relative to $d$, where $\epsilon_n $ is a positive sequence converging to $0$, if the empirical Bayes posterior probability of $U_{\epsilon_n}$ tends to $1$ in $\mathbb{P}_{\theta_0}^{(n)}$-probability. In the case of fully Bayesian procedures, there has been so far a vast literature on posterior consistency and contraction rates since the seminal articles of \citet{barron:schervish:wasserman:1999} and \citet{ghosal:ghosh:vdv:00}. Following  ideas of \citet{schwartz}, \citet{ghosal:ghosh:vdv:00} in the case of independent and identically distributed (iid) observations and \citet{ghosal:vdv:07} in the case of non-iid observations have developed an elegant and powerful methodology to assess posterior contraction rates which boils down to lower bound the prior mass of Kullback-Leibler type neighborhoods of $\mathbb P_{\theta_0}^{(n)}$ and to construct exponentially powerful tests for the testing problem $H_0:\, \theta = \theta_0$ versus $H_1:\theta\in\{\theta:\,d(\theta,\, \theta_0)>\epsilon_n\}$. However, this approach cannot be taken to deal with posterior distributions corresponding to data-dependent priors. Therefore, in this paper, we develop a similar methodology for deriving posterior contraction rates in the case where the prior distribution depends on the data through a data-driven choice $\hga$ of the hyperparameters.

In Section \ref{sec:gene}, we provide a general theorem on posterior contraction rates for empirical Bayes
posterior distributions in the spirit of those presented in \citet{ghosal:vdv:07}, which is then applied to nonparametric mixture models. Two main applications are considered: Dirichlet mixtures of Gaussian distributions for the problem of density estimation in Section \ref {EB:DPMG} and Dirichlet mixtures of uniform random variables for estimating the intensity function of counting processes obeying the Aalen model in Section \ref{sec:aalen}.
Dirichlet process mixtures have been introduced by \citet{ferguson:1974} and have proved to be a major tool in Bayesian nonparametrics, see for instance \citet{hjort:holmes:mueller:2010}.

Rates of convergence for fully Bayesian posterior distributions corresponding to
Dirichlet process mixtures of Gaussian distributions have been widely studied and it has been proved that
they lead to minimax-optimal procedures over a wide collection of density functional classes, see \citet{ghosal:vdv:01}, \citet{ghosal:vdv:mixture:07}, \citet{kruijer:rousseau:vdv:10}, \citet{ghosal:shen:tokdar} and \citet{scricciolo:12}. Here, we extend existing results to the case of a Gaussian base measure with data-dependent mean and variance, as advocated for instance by \citet{richardson:green:1997}. Furthermore, due to some inversion inequalities, we get,
as a by-product, empirical Bayes posterior recovery rates for the problem of density deconvolution when the errors are ordinary or super-smooth and the mixing density is modeled as a Dirichlet process mixture of normal densities with a Gaussian base measure having data-driven chosen mean and variance.
The problem of Bayesian density deconvolution when the mixing density is modeled as a Dirichlet process mixture of Gaussian densities and the errors are super-smooth has been recently studied by \citet{sarkar:13}.

In Section \ref{sec:aalen},  we focus on Aalen multiplicative intensity models, which constitutes a major class of counting processes  that  have been extensively used in the analysis of data arising from various fields like medicine, biology, finance, insurance and social sciences.  General statistical and probabilistic literature on such processes is very huge and we refer the reader to \citet{Karr}, \citet{ABGK}, \citet{MR1950431} and \citet{MR2371524} for a good introduction. In the specific setting of nonparametric Bayesian statistics, practical or methodological contributions have been obtained by \citet{MR1256374}, \citet{AMM09} or  \citet{MR3036340}. First quite general theoretical results have been obtained by \citet{belitser:serra:vanzanten} who established the frequentist asymptotic behavior of posterior distributions for intensities of Poisson processes. We extend their results by considering Aalen multiplicative intensity models instead of simple Poisson models (see Theorem \ref{th:gene:aalen}). In Theorem \ref{cor:EB}, we derive rates of convergence for empirical Bayes estimation of monotone non-increasing intensity functions in counting processes satisfying the Aalen multiplicative intensity model using
Dirichlet process mixtures of uniform distributions with a truncated gamma base measure whose scale parameter is estimated from the data. Numerical illustrations are also presented, in this context.

Proofs and technical derivations are postponed to Section \ref{sec:proofs}.
Instrumental and auxiliary results are reported in the Appendix in Section \ref{app:Aalen}.
\paragraph{Notations and context.} Let $\Data\in \Xn$ be the observations with $(\Xn,\,  \mathcal A_n,\, \P_\theta^{(n)},\, \theta \in \Theta)$ a sequence of experiments, where $\Xn$ and $\Theta$ are Polish spaces endowed with their Borel $\sigma$-fields $\mathcal A_n$ and $\mathcal B$ respectively. We assume that there exists a $\sigma$-finite
measure $\mu^{(n)} $ on $\Xn$ (for convenience of notation, we suppress dependence on $n$ in $\mu^{(n)}$ in what follows) dominating all probability measures $\P_\theta^{(n)}$, $\theta \in \Theta$. For each $\theta\in\Theta,$ we denote $\ell_n(\theta)$ the associated log-likelihood. We also denote by $\E_\theta^{(n)} $ the expectation with respect to $\P_\theta^{(n)}$. We consider a family of prior distributions $(\pi(\cdot\mid\gamma))_{\gamma \in \Gamma}$ on $\Theta$, where $\Gamma \subset \R^d$ with $d \geq 1$, and  we denote by $\pi( \cdot\mid \Data,\, \gamma)$ the posterior corresponding to the prior $\pi(\cdot\mid\gamma)$ which is given by
$$\pi( B \mid \Data,\, \gamma) = \frac{ \int_B e^{\ell_n(\theta)} \d\pi(\theta\mid\gamma)}{ \int_\Theta e^{\ell_n(\theta)} \d\pi(\theta\mid\gamma) },\quad B\in \mathcal B. $$
We denote by $\textrm{KL}(\cdot;\, \cdot) $ the Kullback-Leibler divergence. Given $\theta_1,\theta_2\in\Theta_n$, we denote by $V_k(\theta_1;\, \theta_2)$ the recentered $k$-th moment of the log-likelihood difference associated to $\theta_1$ and $\theta_2$. So, we have:
$$\textrm{KL}(\theta_1;\, \theta_2)=\E_{\theta_1}^{(n)}[\ell_n(\theta_1)-\ell_n(\theta_2)],$$ $$V_{k}(\theta_1;\,\theta_2) = \E_{\theta_1}^{(n)}[|\ell_n(\theta_1)-\ell_n(\theta_2)-\E_{\theta_1}^{(n)}[\ell_n(\theta_1)-\ell_n(\theta_2)]|^{k}].$$
We denote by $h(f_1,f_2) $ the Hellinger metric between two densities $f_1$ and $f_2$, \emph{i.e.}, $h^2(f_1,\, f_2) = \int (\sqrt{f_1}(x) - \sqrt{f_2}(x))^2 \d x $. Throughout the text, we denote by $D(\zeta ,\, B,\, d(\cdot,\,\cdot) )$ the $\zeta$-covering number of $B$ by $d$-balls  of radius  $\zeta$, for any set $B$, any positive constant $\zeta$ and
any pseudo-metric $d(\cdot,\,\cdot)$. We denote $\theta_0$ the true parameter.

\section{ General result on posterior concentration rates for empirical Bayes} \label{sec:gene}
Let $\hga$ be a measurable function of the observations. The associated empirical Bayes posterior distribution is then $\pi(\cdot \mid\Data,\, \hga)$. In this section, we present a general theorem  to obtain posterior concentration rates for the empirical Bayes posterior $\pi(\cdot \mid\Data,\, \hga)$. Our aim is to give conditions resembling those usually considered in fully Bayesian approaches. For this purpose, we first define the usual quantities. We assume that, with probability going to 1, $\hga$ belongs to a subset $\mathcal K_n$ of $\Gamma$:
\begin{equation}\label{cond:gamma}
\P_{\theta_0}^{(n)}\left( \hga \in \mathcal K_n^c  \right) =o(1).
\end{equation}
For any positive sequence $(u_n)_n$, let $N_n(u_n)$ stand for the $u_n$-covering number of
$\mathcal K_n$ relative to the Euclidean distance which is denoted by $\|\cdot\|$. For instance, if $\mathcal K_n $ is included in a ball of radius $R_n$ then $N_n(u_n) = O((R_n/u_n)^d)$.

As in \citet{ghosal:vdv:07}, for any  $\epsilon_n>0$, we introduce the $\epsilon_n$-Kullback-Leibler neighborhood of $\theta_0$ defined, for $k\geq 1$, by
$$\bar B_{k,n} = \{ \theta:\, \textrm{KL}(\theta_0;\, \theta)\leq n \epsilon_n^2 , V_k(\theta_0;\, \theta) \leq n^{k/2} \epsilon_n^k\}.$$
Let $d(\cdot,\,\cdot) $ be a loss function on $\Theta$. We consider the posterior concentration rates in terms of $d(\cdot,\,\cdot) $ using the following neighborhoods:
$$U_{M\epsilon_n} = \{ \theta\in\Theta:\, d(\theta_0,\, \theta) \leq M \epsilon_n \}$$
with $M>0$. For any integer $j$, we define $$S_{n,j} =  \{ \theta \in \Theta:\, d(\theta_0,\, \theta) \in (j\epsilon_n,\, (j+1)\epsilon_n]\}.$$
In order to obtain posterior concentration rates with data-dependent priors, we express the impact of $\hga$ on the prior in the following way: for all $\gamma,\, \gamma' \in \Gamma$, we construct a measurable transformation $\psi_{\gamma, \gamma'} : \Theta \rightarrow \Theta$ such that if $\theta \sim \pi(\cdot \mid \gamma) $ then $\psi_{\gamma, \gamma'}(\theta) \sim \pi(\cdot \mid \gamma') $. Let $e_n(\cdot,\,\cdot)$ be another semi-metric on $\Theta\times\Theta$. We consider the following set of assumptions. Let $\Theta_n \subset \Theta $ and $k$ be fixed.

\medskip

\noindent {\bf [A1]}
There exist a sequence $(u_n)_n$ and $\tilde B_n \subset \bar B_{k,n}$  such that
\begin{equation}\label{Nn}
N_n(u_n) = o((n\epsilon_n^2)^{k/2})
\end{equation}
and
\begin{equation}\label{KLcond}
\sup_{\gamma\in \mathcal K_n} \sup_{\theta \in \tilde B_n} \P_{\theta_0}^{(n)}\left\{ \inf_{\gamma': \ \|\gamma' -\gamma\| \leq u_n} \ell_n(\psi_{\gamma, \gamma'}(\theta)) - \ell_n(\theta_0) < -n\epsilon_n^2  \right\} = o(N_n(u_n)^{-1}).
\end{equation}
{\bf [A2]} Defining for all $\gamma \in \mathcal K_n$,
 $\d Q_{\gamma, n}^\theta(\data)= \sup_{\| \gamma'  - \gamma\| \leq u_n} e^{\ell_n(\psi_{\gamma, \gamma'}(\theta))(\data)}\d\mu(\data) $,
we have
\begin{equation} \label{Thetanc}
\sup_{\gamma \in \mathcal K_n} \frac{\int_{\Theta\setminus\Theta_n} Q_{\gamma, n}^\theta(\Xn)\d\pi(\theta\mid\gamma) }{ \pi(\tilde B_n\mid\gamma)}  =o(N_n(u_n)^{-1}e^{- 2n\epsilon_n^2})
\end{equation}
 and there exist constants $\zeta,\,K>0$ such that
 \begin{itemize}
 \item for all $j$ large enough,
    \begin{equation}\label{ratio:mass}
    \sup_{\gamma \in \mathcal K_n} \frac{ \pi(S_{n,j} \cap \Theta_n\mid\gamma ) }{ \pi(\tilde B_n\mid\gamma ) }  \leq e^{ K n j^2 \epsilon_n^2/2},  
  \end{equation}
  \item  for all $\epsilon>0$,  for all $\theta \in \Theta_n$  with $d(\theta_0,\, \theta)>\epsilon$ there exist tests $\phi_n(\theta)$ satisfying
 \begin{equation} \label{test}
  \E_{\theta_0}^{(n)}\left[ \phi_{n}(\theta) \right] \leq e^{-K n \epsilon^2 } ,  \quad \sup_{\gamma\in{\mathcal K}_n}
\sup_{ e_n(\theta',\theta) \leq \zeta \epsilon } \int_{\Xn} [1-\phi_n(\theta)] \d Q_{\gamma, n}^{\theta'}\leq e^{-K n \epsilon^2 },
\end{equation}
 \item for all $j$ large enough, \begin{equation}\label{eq:entropy}
 \log D(\zeta j\epsilon_n,\, S_{n,j},\, e_n(\cdot,\,\cdot)) \leq K(j+1)^2 n\epsilon_n^2 /2.
 \end{equation}
\end{itemize}
We then have the following theorem.
\begin{Th}\label{th:gene:EB}
Let $\theta_0\in \Theta$. Assume that the estimator $\hga$ satisfies \eqref{cond:gamma} and the prior satisfies assumptions [A1] and [A2], with $n\epsilon_n^2\rightarrow\infty $ and $\epsilon_n\rightarrow0 $. Then, for $J_1$ large enough,
$$\E_{\theta_0}^{(n)} [\pi( U_{J_1\epsilon_n}^c \mid\Data ,\, \hga)] = o(1),$$
where $U_{J_1\epsilon_n}^c$ is the complementary of $U_{J_1\epsilon_n}$ in $\Theta$.
\end{Th}
The proof of Theorem \ref{th:gene:EB} shows that the posterior concentration rate of the empirical Bayes posterior based on $\hga \in \mathcal K_n$ is bounded from above by the worst concentration rate over the classes of posterior distributions corresponding to the priors $(\pi(\cdot \mid \gamma),\, \gamma \in \mathcal K_n)$. In particular, assume that each posterior $\pi(\cdot \mid\Data ,\,  \gamma)$ converges at rate $\epsilon_n(\gamma) = (n/\log n)^{-\alpha(\gamma)}$, where $\gamma\to\alpha(\gamma)$ is Lipschitzian,  and $\hga$ converges to $\gamma^*$ at rate $v_n,$ with $v_n=o((\log n)^{-1})$ namely $\mathcal K_n=[\gamma^*-v_n,\gamma^*+v_n]$, then the posterior concentration rate of the empirical Bayes posterior is of the order $O(\epsilon_n(\gamma^*) )$. Indeed, in this case, $\sup_{\gamma\in \mathcal K_n}\{\epsilon_n(\gamma) /\epsilon_n(\gamma^*)\}=\sup_{\gamma\in \mathcal K_n} (n/\log n)^{|\alpha(\gamma)-\alpha(\gamma^*|)}=O(1)$.
This is of special interest when $\epsilon_n(\gamma^*) $ is optimal. Proving that the posterior distribution has optimal posterior concentration rate then boils down to proving that $\hga$ converges to  the oracle value $\gamma^*$.

\begin{remark} \label{rem:loss}
As in \citet{ghosal:vdv:07}, we can replace conditions \eqref{test}  and \eqref{eq:entropy} by the existence of a global test $\phi_n$ over $S_{n,j}$ similar to equation (2.7) of \citet{ghosal:vdv:07} satisfying
$$  \E_{\theta_0}^{(n)}\left[ \phi_{n}\right] = o(N_n(u_n)^{-1}) , \quad \sup_{\gamma\in{\mathcal K}_n}
 \sup_{\theta\in S_{n,j}} \int_{\Xn} (1-\phi_n) \mathrm{d} Q_{\gamma, n}^\theta \leq e^{-K n j^2 \epsilon_n^2 }
$$
without modifying the conclusion. Note also that when the loss function $d(\theta,\, \theta_0)$ is not bounded it is often the case that obtaining exponential control on the error rates in the form $e^{-K n \epsilon^2 }$ or $e^{-K n j^2 \epsilon_n^2 }$ is not possible for large values of $j$. It is enough in that case to consider a
 modification $\tilde d(\theta,\, \theta_0)$ of the loss  which affects only the values of $\theta$ for which  $d(\theta,\, \theta_0) $ is large and to prove  \eqref{test} and \eqref{eq:entropy} for $\tilde d( \theta,\, \theta_0)$, defining  $S_{n,j}$ and the covering number $D(\cdot ) $ with respect to  $\tilde d(\cdot,\, \cdot ) $. As an illustration of this remark,  see the proof of Theorem \ref{cor:EB}.
\end{remark}

The assumptions of [A2] are very similar to those for establishing posterior concentration rates proposed for instance by \citet{ghosal:vdv:07} (see their Theorem 1 and the associated proof). We need to strengthen some conditions to take into account that we only know that $\hga$ lies in the compact set $\mathcal{K}_n$ with high probability.

The key idea here is to construct a transformation
$\psi_{\gamma, \gamma'}$ which allows to transfer the dependence on the data from the prior to the likelihood, similarly to what is considered in \citet{petrone:rousseau:scricciolo:14}. The only difference with the general  theorems of \citet{ghosal:vdv:07} lies in the control of
the log-likelihood difference $\ell_n(\psi_{\gamma, \gamma'}(\theta)) - \ell_n(\theta_0)$ when $\|\gamma -\gamma'\|\leq u_n$. In nonparametric cases where $n\epsilon_n^2$ is a power of $n$, $u_n$ can be chosen very small as soon as $k$ can be chosen large enough so that controlling this difference uniformly is not such a drastic condition. In parametric models, where at best $n\epsilon_n^2$ is a power of $\log n$, this become more involved and $u_n$ needs to be large or $\mathcal K_n$ needs to be small enough so that $N_n(u_n)$ can be chosen of order $O(1)$. Note that in parametric models it is typically easier to use a more direct control of $\pi(\theta \mid \gamma)/\pi( \theta \mid \gamma')$,
the ratio of the prior densities with respect to a common dominating measure.  In nonparametric prior models this is usually not possible since no such dominating measure exists in most cases. In Sections  \ref{EB:DPMG} and  \ref{sec:aalen}, we apply Theorem \ref{th:gene:EB} to two different types of Dirichlet process mixture models: Dirichlet process mixtures of Gaussian distributions used to model smooth densities and Dirichlet process mixture model of uniforms to model monotone non-increasing intensities in the context of Aalen point processes. It is interesting to note that in the case of general nonparametric mixture models there exists a general construction of 
$\psi_{\gamma, \gamma'}$. More precisely, consider a mixture model in the form
 \begin{equation} \label{NPmixt}
 f(\cdot)  = \sum_{j=1}^K p_j h_{\theta_j}(\cdot) , \quad K \sim \pi_K,
 \end{equation}
and, conditionally on $K$, $p= (p_j)_{j=1}^K \sim \pi_p$ and $\theta_1,\, \ldots,\, \theta_K$ are iid with cdf $G_\gamma$. The Dirichlet process mixture corresponds to $\pi_K  = \delta_{(+\infty)}$ and $\pi_p$ is the GEM distribution obtained from the stick-breaking construction, see for instance \citet{ghosh:ramamo:2003}. Models in the form \eqref{NPmixt} also cover priors on curves if the $(p_j)_{j=1}^K$ are not restricted to the simplex. Denote by $\pi( \cdot \mid \gamma)  $ the prior probability on $f$ induced by \eqref{NPmixt}. Then, for all $\gamma, \gamma' \in \Gamma$, if $ f$ is represented as \eqref{NPmixt} and is thus distributed according to $\pi( \cdot \mid \gamma)$, then
 $$ f^{'}(\cdot) =  \sum_{j=1}^K p_j h_{\theta_j^{'}}(\cdot) , \quad \mbox{ with }\,\,\,\theta_j^{'} = G_{\gamma^{'}}^{-1}( G_\gamma( \theta_j) ), $$
is distributed according to $\pi( \cdot \mid \gamma')$, where $G_{\gamma^{'}}^{-1}(\cdot)$ denotes the generalized inverse of the cumulative distribution function.

\medskip

We now give the proof of Theorem \ref{th:gene:EB}.

\medskip

\begin{proof}  We consider a chaining argument and we split ${\mathcal K}_n$ into $N_n(u_n)$ balls of radius $u_n$ and we denote by $(\gamma_i)_{i=1,\,\ldots,\,N_n(u_n)}$ the centers of these balls. We have
\begin{equation}\label{start}
\E_{\theta_0}^{(n)} [ \pi( U_{J_1\epsilon_n}^c \mid\Data,\, \hga)] \leq\E_{\theta_0}^{(n)}[ \max_i \rho_n(\gamma_i)]
\leq N_n(u_n)\max_{i=1,\,\ldots,\,N_n(u_n)}\E_{\theta_0}^{(n)} \left[ \rho_n(\gamma_i)\right],
\end{equation}
with
$$\rho_n(\gamma_i) :=   \sup_{\|\gamma - \gamma_i\|\leq u_n} \pi( U_{J_1\epsilon_n}^c \mid\Data,\, \gamma)=  \sup_{\|\gamma - \gamma_i\|\leq u_n}
 \frac{ \int_{ U_{J_1\epsilon_n}^c} e^{\ell_n(\psi(\gamma_i, \gamma)(\theta)) - \ell_n(\theta_0) } \d\pi(\theta\mid \gamma_i ) }{ \int_{ \Theta }  e^{\ell_n(\psi(\gamma_i, \gamma)(\theta)) - \ell_n(\theta_0) } \d\pi(\theta\mid \gamma_i) }.$$
We now study, for any $i$, $\E_{\theta_0}^{(n)} \left[ \rho_n(\gamma_i)\right] $.
We mimic the proof of Lemma~9 of  \citet{ghosal:vdv:07}. For every $j$ large enough, by (\ref{eq:entropy}), there exist  $L_{j,n}$ balls of centers $\theta_{j,1},\,\dots,\,\theta_{j,L_{j,n}}$, with radius  $\zeta j\epsilon_n$ relative to the $e_n$-distance, covering $S_{n,j}\cap\Theta_n$, with $L_{j,n}\leq \exp(K(j+1)^2n\epsilon_n^2 /2)$. We then consider tests $\phi_n(\theta_{j,\ell})$ satisfying (\ref{test}) with $\epsilon=j\epsilon_n$. By setting
$$\phi_n= \max_{j\geq J_1} \max_{\ell\in\{1,\,\ldots,\, L_{j,n}\}} \phi_n(\theta_{j,\ell}),$$
we obtain
\begin{equation}\label{test0}
\E_{\theta_0}^{(n)}[\phi_n] \leq \sum_{j \geq J_1} L_{j,n} e^{-Knj^2 \epsilon_n^2 }  =O(e^{-KJ_1^2n\epsilon_n^2/2}).
\end{equation}
Moreover, for any $j\geq J_1$, any $\theta\in S_{n,j}\cap\Theta_n$ and any $i$,
\begin{equation}\label{test1}
\int_{\Xn} (1-\phi_n) \d Q_{\gamma_i,n}^\theta\leq e^{-K n j^2\epsilon_n^2}.
\end{equation}
Since for all $i$ we have $\rho_n(\gamma_i)\leq 1$, using the usual decomposition,
\begin{equation}\label{decomp}
\E_{\theta_0}^{(n)}\left[ \rho_n(\gamma_i) \right]\leq \E_{\theta_0}^{(n)}[\phi_n] + \P_{\theta_0}^{(n)}(A_{n,i}^c)+\frac{ e^{2n \epsilon_n^2 } }{ \pi(\tilde B_n\mid\gamma_i) } C_{n,i},
\end{equation}
with
$$A_{n,i}=\left\{\inf_{\|\gamma - \gamma_i\|\leq u_n }\int_{ \Theta }  e^{\ell_n(\psi(\gamma_i, \gamma)(\theta)) - \ell_n(\theta_0) } \d\pi(\theta\mid \gamma_i) \geq e^{-2n\epsilon_n^2}\pi(\tilde B_n\mid\gamma_i)\right\}$$
and
$$C_{n,i}=\E_{\theta_0}^{(n)}\left[(1-\phi_n)\sup_{\|\gamma - \gamma_i\|\leq u_n} \int_{ U_{J_1\epsilon_n}^c} e^{\ell_n(\psi(\gamma_i, \gamma)(\theta)) - \ell_n(\theta_0) } \d\pi(\theta\mid \gamma_i )\right].$$
We study each term of (\ref{decomp}). First, by using (\ref{Nn}) and (\ref{test0}),
$$\E_{\theta_0}^{(n)}[\phi_n]=o(N_n(u_n)^{-1}).$$
Secondly, since
$$e^{n\epsilon_n^2}\inf_{\|\gamma - \gamma_i\|\leq u_n }e^{\ell_n(\psi(\gamma_i, \gamma)(\theta)) - \ell_n(\theta_0) }\geq 1_{\{\inf_{\|\gamma - \gamma_i\|\leq u_n }e^{\ell_n(\psi(\gamma_i, \gamma)(\theta)) - \ell_n(\theta_0) }\geq e^{-n\epsilon_n^2}\}},$$
we have
\begin{eqnarray*}
\P_{\theta_0}^{(n)}(A_{n,i}^c)&\leq&\P_{\theta_0}^{(n)}\left\{\int_{\tilde B_n }  \inf_{\|\gamma - \gamma_i\|\leq u_n }e^{\ell_n(\psi(\gamma_i, \gamma)(\theta)) - \ell_n(\theta_0) } \frac{\d\pi(\theta\mid \gamma_i)}{\pi(\tilde B_n\mid\gamma_i)} < e^{-2n\epsilon_n^2}\right\}\\
&\leq&\P_{\theta_0}^{(n)}\left\{\int_{\tilde B_n }1_{\{\inf_{\|\gamma - \gamma_i\|\leq u_n }e^{\ell_n(\psi(\gamma_i, \gamma)(\theta)) - \ell_n(\theta_0) }\geq e^{-n\epsilon_n^2}\}} \frac{\d\pi(\theta\mid \gamma_i)}{\pi(\tilde B_n\mid\gamma_i)} < e^{-n\epsilon_n^2}\right\}\\
&\leq&\P_{\theta_0}^{(n)}\left\{\int_{\tilde B_n }1_{\{\inf_{\|\gamma - \gamma_i\|\leq u_n }e^{\ell_n(\psi(\gamma_i, \gamma)(\theta)) - \ell_n(\theta_0) }< e^{-n\epsilon_n^2}\}} \frac{\d\pi(\theta\mid \gamma_i)}{\pi(\tilde B_n\mid\gamma_i)} >1- e^{-n\epsilon_n^2}\right\}\\
&\leq&(1-e^{-n\epsilon_n^2})^{-1}\int_{\tilde B_n }\P_{\theta_0}^{(n)}\left\{ \inf_{\|\gamma-\gamma_i\| \leq u_n} \ell_n(\psi_{\gamma_i, \gamma}(\theta)) - \ell_n(\theta_0) \leq -n\epsilon_n^2  \right\}\frac{\d\pi(\theta\mid \gamma_i)}{\pi(\tilde B_n\mid\gamma_i)}\\
&=&o(N_n(u_n)^{-1}),
\end{eqnarray*}
by (\ref{KLcond}).
Also, we have
\begin{eqnarray*}
C_{n,i}&\leq&\E_{\theta_0}^{(n)}\left[(1-\phi_n)\int_{ U_{J_1\epsilon_n}^c} \sup_{\|\gamma - \gamma_i\|\leq u_n} e^{\ell_n(\psi(\gamma_i, \gamma)(\theta)) - \ell_n(\theta_0) } \d\pi(\theta\mid \gamma_i )\right]\\
&\leq&\int_{U_{J_1\epsilon_n}^c} \int_{\Xn} (1-\phi_n) \d Q_{\gamma_i,n}^\theta \d\pi(\theta\mid\gamma_i)\\
&\leq & \int_{\Theta\setminus\Theta_n}Q_{\gamma_i,n}^\theta (\Xn)\d\pi(\theta\mid\gamma_i) +\sum_{j\geq J_1}\int_{S_{n,j}\cap \Theta_n} \int_{\Xn} (1-\phi_n) \d Q_{\gamma_i,n}^\theta \d\pi(\theta\mid\gamma_i).
\end{eqnarray*}
By using (\ref{Thetanc}), (\ref{ratio:mass}) and (\ref{test1}),
\begin{eqnarray*}
C_{n,i}&\leq& \sum_{j\geq J_1}e^{-Knj^2\epsilon_n^2}\pi(S_{n,j}\cap\Theta_n\mid\gamma_i)+o(N_n(u_n)^{-1}e^{-2n\epsilon_n^2}\pi(\tilde B_n\mid\gamma_i))\\
&\leq&\sum_{j\geq J_1}e^{-Knj^2\epsilon_n^2/2}\pi(\tilde B_n\mid\gamma_i)+o(N_n(u_n)^{-1}e^{-2n\epsilon_n^2}\pi(\tilde B_n\mid\gamma_i))\\
&=&o(N_n(u_n)^{-1}e^{-2n\epsilon_n^2}\pi(\tilde B_n\mid\gamma_i))
\end{eqnarray*}
and
$$\max_{i=1,\,\ldots,\,N_n(u_n)}\frac{ e^{2n \epsilon_n^2 } }{ \pi(\tilde B_n\mid\gamma_i) } C_{n,i}=o(N_n(u_n)^{-1}).$$
Having controlled each term in (\ref{decomp}), by (\ref{start}),
the proof of Theorem \ref{th:gene:EB} is achieved.
\end{proof}
\section{Adaptive posterior contraction rates for empirical Bayes Dirichlet process mixtures of Gaussian densities} \label{EB:DPMG}
Let $\Data = (X_1,\,\ldots,\,X_n)$ be $n$ iid observations from an unknown Lebesgue density $p$ on $\mathbb{R}$.
Consider the following prior distribution on the class of densities $\mathcal P = \{p:  \R \rightarrow \R_+\mid\, \|p\|_1=1\}$:
\begin{equation}\label{DPMGaussian}
\begin{split}
p(\cdot) = p_{F,\sigma}(\cdot) &: =\int_{-\infty}^\infty\phi_\sigma(\cdot-\theta)\,\mathrm{d}F(\theta),\\
F\sim \mathrm{DP}(\alpha_{\mathbb{R}}\Nor_{(m, s^2)}), &\quad
\sigma \sim\mathrm{IG}(\nu_1,\,\nu_2), \quad \nu_1 ,\, \nu_2 >0,
\end{split}
\end{equation}
where $\alpha_{\mathbb{R}}$ is a positive constant,
$\phi_\sigma(\cdot)=\sigma^{-1}\phi(\cdot/\sigma)$, with $\phi(\cdot)$ the density of a standard Gaussian distribution, and
$\Nor_{(m,s^2)}$ denotes the Gaussian distribution with mean
$m$ and variance $s^2$. Set $\gamma = (m,\,s^2)\in \Gamma\subseteq\mathbb{R}\times \mathbb{R}_+^*$, let $\hga:\,\mathbb{R}^n\rightarrow\Gamma$ be some measurable function of the observations.
Typical choices are $\hga= (\bar X_n,\, S_n^2) $, with the sample mean $\bar X_n = \sum_{i=1}^n X_i/n$ and the sample variance $S_n^2=\sum_{i=1}^n(X_i-\bar{X}_n)^2/n$, and
$\hga= (\bar X_n,\, R_n)$, with the range $R_n = \max_i X_i - \min_i X_i$
as in \citet{richardson:green:1997}. Let $\mathcal K_n \subset \R\times \R_+^{*}$ be compact, independent of the data $\Data$ and such that
\begin{equation}\label{cond:DPM:Kn}
\P^{(n)}_{p_0}\left( \hga\in \mathcal K_n \right) = 1 + o(1).
\end{equation}
Throughout this section, we assume that the true density
$p_0$ satisfies the following tail condition:
\begin{equation}\label{cond:tail1}
p_0(x)\lesssim e^{-c_0|x|^\tau} \quad \mbox{ for }\,  |x|\mbox{ large enough},
\end{equation}
with finite constants $c_0,\,\tau>0$. Let $\mathbb{E}_{p_0}[X_1]= m_0\in\mathbb{R}$ and $\mathbb{V}\textrm{ar}_{p_0}[X_1]=\tau_0^2\in\mathbb{R}^*_+$. If $\hga = (\bar X_n,\, S_n^2)$, then condition \eqref{cond:DPM:Kn} is satisfied with $\mathcal K_n = [m_0- (\log n)/\sqrt{n},\, m_0+(\log n)/\sqrt{n}] \times  [\tau_0^2-(\log n)/\sqrt{n},\, \tau_0^2 +(\log n)/\sqrt{n}]$,  while, if $\hga= (\bar X_n,\, R_n) $, then $\mathcal K_n = [m_0- (\log n)/\sqrt{n},\, m_0+(\log n)/\sqrt{n}] \times [a,\, 2(2c_0^{-1} \log n)^{1/\tau }]$.

Mixtures of Gaussian densities have been extensively used and studied in Bayesian nonparametric literature. Posterior contraction rates have been first investigated by \citet{ghosal:vdv:01} and \citet{ghosal:vdv:mixture:07}. Following an idea of \citet{rousseau:09}, \citet{kruijer:rousseau:vdv:10} have proved that nonparametric location mixtures of Gaussian densities lead to adaptive posterior contraction rates over the full scale of locally H\"older log-densities on $\R$. This result has been extended to the multivariate set-up by \citet{ghosal:shen:tokdar} and to super-smooth densities by \citet{scricciolo:12}.
The key idea behind these results is that,
for an ordinary $\beta$-smooth density $p_0$,
given $\sigma>0$ small enough, there exists a finite mixing distribution $F^*$,
with $N_\sigma=O(\sigma^{-1}|\log \sigma|^{\rho_2})$
support points in $[-a_\sigma,\, a_\sigma]$, for $a_\sigma = O(|\log \sigma|^{1/\tau})$,
such that the corresponding Gaussian mixture density $p_{F^*, \sigma}$ satisfies
\begin{equation}\label{Holder}
n^{-1}\textrm{KL}(p_0;\, p_{F^*, \sigma})
\lesssim \sigma^{2\beta} , \qquad n^{-k/2}V_k (p_0;\, p_{F^*, \sigma}) \lesssim \sigma^{k\beta}, \quad k\geq 2,
\end{equation}
see, for instance, Lemma 4 in \citet{kruijer:rousseau:vdv:10}. In all these articles, only the case where $k=2$
has been treated for the inequality on the right-hand side (RHS) of \eqref{Holder}, but the extension to
any $k>2$ is straightforward.
The regularity assumptions considered in \citet{kruijer:rousseau:vdv:10}, \citet{ghosal:shen:tokdar} and \citet{scricciolo:12} to verify \eqref{Holder} are slightly different. For instance, \citet{kruijer:rousseau:vdv:10} assume that $\log p_0$ satisfies some locally H\"older conditions, while \citet{ghosal:shen:tokdar} consider H\"older-type conditions on $p_0$ and \citet{scricciolo:12} Sobolev-type assumptions. To avoid taking into account all these special cases, we state \eqref{Holder} as a condition. We then have the following theorem, where the distance $d$ defining the ball $U_{J_1\epsilon_n}$ with center at $p_0$ and radius $J_1\epsilon_n$ can equivalently be the Hellinger or the $\mathbb{L}_1$-distance. Note that the constant $J_1$ may be different for each one of the results stated below.

\begin{Th}\label{Th:DPMGaussian1}
Suppose that $p_0$ satisfies the tail condition \eqref{cond:tail1} and that the inequality on the RHS of
\eqref{Holder} holds with $k>8(2\beta+1)$. Consider a prior distribution of the form
\eqref{DPMGaussian}, with an empirical Bayes selection $\hga$ for $\gamma$, and assume that $\hga$ and $\mathcal K_n$ satisfy condition \eqref{cond:DPM:Kn}, where
$\mathcal K_n \subseteq [m_1,\, m_2] \times [a_1,\, a_2(\log n)^{b_1}]$ for some constants $m_1,\,m_2 \in \R$,  $a_1,\, a_2>0$ and $b_1\geq 0$. Then, for a sufficiently large constant $J_1>0$,
$$\E_{p_0}^{(n)} [\pi( U_{J_1\epsilon_n}^c \mid\Data ,\,\hga) ] = o(1), \quad \mbox{ with } \epsilon_n = n^{-\beta / (2\beta+1)} (\log n)^{a_3},$$
for some constant $a_3>0$.
\end{Th}

In Theorem \ref{Th:DPMGaussian1}, the constant $a_3$ is the same as that appearing in the rate of convergence for the
posterior distribution corresponding to a non data-dependent prior with fixed $\gamma$.

\medskip

The crucial step for assessing posterior contraction rates in the case where $p_0$ is super-smooth, which, as in the ordinary smooth case,
consists in proving the existence of a finitely supported Gaussian mixture density that approximates $p_0$ with an error of the appropriate order,
requires some refinements. We suppose that $p_0$ has
Fourier transform $\widehat{p_0}(t)=\int_{-\infty}^\infty e^{itx}p_0(x)\d x$, $t\in\mathbb{R}$, that satisfies for some finite
constants $\rho,\,L>0$ the integrability condition
\begin{equation}\label{Ftransf}
\int_{-\infty}^\infty|\widehat{p_0}(t)|^2e^{-2(\rho|t|)^{r}}\d t\leq 2\pi L^2,\quad \mbox{ with $r\in[1,\,2]$,}
\end{equation}
where the regularity of $p_0$, which is measured through a scale of integrated tail 
bounds on the Fourier transform $\widehat{p_0}$, is related to
the characteristic exponent $r$. Densities satisfying condition \eqref{Ftransf} are analytic on $\mathbb{R}$ and
increasingly \vir{smooth} as $r$ increases. They
form a larger class than that of analytic densities, including relevant statistical examples like
the Gaussian distribution which corresponds to $r=2$, the Cauchy distribution which corresponds to $r=1$, all symmetric stable laws,
Student's-$t$ distribution, distributions with characteristic
functions vanishing outside a compact set as well as their mixtures and convolutions.
In order to state a counter-part of requirement \eqref{Holder} for the super-smooth case,
we consider, for $\alpha\in(0,\,1]$, the $\rho_\alpha$-divergence of a density $p$ from $p_0$ which is defined as
$\rho_\alpha(p_0;\,p)=\alpha^{-1}\{\mathbb{E}_{p_0}[(p_0/p)^\alpha(X_1)]-1\}$. Following the trail of Lemma 8 in \citet{scricciolo:12}, it can be proved that, for any density $p_0$ satisfying condition \eqref{Ftransf}, together with
the monotonicity and tail conditions $(b)$ and $(c)$, respectively, of Section 4.2 in \citet{scricciolo:12},
for $\sigma>0$ small enough,
there exists a finite mixing distribution $F^*$, with $N_\sigma$ support points in $[-a_\sigma,\, a_\sigma]$, so that
\begin{equation}\label{analytic}
\rho_\alpha(p_0;\,p_{F^*, \sigma})\lesssim e^{-c(1/\sigma)^{r}} \quad \mbox{ for every $\alpha\in(0,\,1]$,}
\end{equation}
where $a_\sigma =O(\sigma^{-r/(\tau\wedge2)})$ and  $N_\sigma=O((a_\sigma/\sigma)^2)$.
Since
$$n^{-1}\textrm{KL}(p_0;\,p_{F^*, \sigma})=\lim_{\beta\rightarrow0^+}\rho_{\beta}(p_0;\,p_{F^*, \sigma})\leq \rho_\alpha(p_0;\,p_{F^*, \sigma}) \quad\mbox{ for every $\alpha\in(0,\,1]$,}$$
inequality \eqref{analytic}
is stronger than that on the LHS of \eqref{Holder} and allows to obtain an almost sure lower bound on the denominator of the ratio
defining the empirical Bayes posterior probability of the set $U_{J_1\epsilon_n}^c$, see Lemma 2 of \citet{shen:wasserman:01}.


\begin{Th}\label{Th:DPMGaussian2}
Suppose that $p_0$ satisfies condition \eqref{Ftransf} and that the tail condition
\eqref{cond:tail1} holds with $\tau>1$ such that $(\tau-1)r\leq\tau$.
Suppose also that the monotonicity condition $(b)$
of Section 4.2 in \citet{scricciolo:12} is satisfied.
Consider a prior distribution of the form
\eqref{DPMGaussian}, with an empirical Bayes selection $\hga$ for $\gamma$, and assume that $\hga$ and $\mathcal K_n$ satisfy condition \eqref{cond:DPM:Kn}, where $\mathcal K_n \subseteq [m_1,\, m_2] \times [a_1,\, a_2(\log n)^{b_1}]$ for some
constants $m_1,\, m_2 \in \R$,  $a_1,\, a_2>0$ and $b_1\geq 0$. Then, for a sufficiently large constant $J_1>0$,
$$\mathbb{E}^{(n)}_{p_0}[\pi( U_{J_1\epsilon_n}^c \mid \Data ,\, \hga)]=o(1),\quad \mbox{ with $\epsilon_n = n^{-1/2}(\log n)^{5/2+3/r}$.}$$
\end{Th}

\medskip

We now present some results on empirical Bayes posterior recovery rates for mixing distributions. These results are
derived from Theorem~\ref{Th:DPMGaussian1} and Theorem~\ref{Th:DPMGaussian2} via some inversion inequalities.
We first consider the case where the sampling density $p_0$ is itself a mixture of Gaussian densities and derive rates for recovering the true mixing distribution relative to any Wasserstein metric of order $1\leq q<\infty$. We then assess empirical Bayes
posterior recovery rates relative to the $\mathbb{L}_2$-distance for the
density deconvolution problem in the ordinary and super-smooth cases.

We begin by recalling the definition of \emph{Wasserstein distance of order $q$}.
For any $1\leq q<\infty$, define the \emph{Wasserstein distance of order $q$}
between any two Borel probability measures $\nu$ and $\nu'$ on $\Theta$ with finite $q$th-moment,
\emph{i.e.}, $\int_\Theta d^q(\theta,\,\theta_0)\nu(\mathrm{d}\theta)<\infty$ for some (hence any $\theta_0$) in $\Theta$,
as
$W_q(\nu,\,\nu'):=(\inf_{\mu\in\Pi(\nu,\,\nu')}\int_{\Theta \times \Theta} d^q(\theta,\,\theta')\mu(\mathrm{d}\theta,\,\mathrm{d}\theta'))^{1/q}$, where $\mu$ runs over the set $\Pi(\nu,\,\nu')$ of all joint probability measures on $\Theta \times \Theta$
with marginal distributions $\nu$ and $\nu'$. When $q=2$, we take $d$ to be the Euclidean distance $\|\cdot\|$.
Posterior rates, relative to Wasserstein metrics, for recovering mixing distributions have
been recently investigated by \citet{nguyen:2013}. In the following result, the prior probability measure corresponds to the product of a Dirichlet process, with a data-dependent base measure $\alpha_{\mathbb{R}}\Nor_{\hga}$, and a point mass at a given $\sigma_0$, in symbols, $\mathrm{DP}(\alpha_{\mathbb{R}}\Nor_{\hga})\times \delta_{\sigma_0}$.

\begin{Cor}\label{Mixing}
Suppose that $p_0= F_0\ast \phi_{\sigma_0}$, where the true mixing distribution $F_0$ satisfies the tail condition
$
F_0(\theta:\,|\theta|>t)\lesssim
e^{-c_0t^2}$ for $t$ large enough
and $\sigma_0$ denotes the true value for the scale.
Consider a prior distribution of the form
$\mathrm{DP}(\alpha_{\mathbb{R}}\Nor_{\hga})\times \delta_{\sigma_0}$ and assume that $\hga$ and $\mathcal K_n$ satisfy condition \eqref{cond:DPM:Kn}, where $\mathcal K_n \subseteq [m_1,\, m_2] \times [a_1,\, a_2(\log n)^{b_1}]$ for some
constants $m_1,\, m_2 \in \R$,  $a_1,\, a_2>0$ and $b_1\geq 0$. Then,
for every $1\leq q<\infty$, there exists a sufficiently large
constant $J_1>0$ so that $$\mathbb{E}^{(n)}_{p_0}[\pi(F:\,W_q(F,\,F_0)\geq J_1(\log n)^{-1/2}\mid \Data,\,\hga)]=o(1).$$
\end{Cor}
\medskip

The result implies that optimal recovery of mixing distributions is possible using Dirichlet process mixtures of Gaussian densities with an empirical Bayes selection for the prior hyper-parameters of the base measure. \citet{dedecker:michel:2013} have shown that, for the deconvolution problem with super-smooth errors, the rate $(\log n)^{-1/2}$
is minimax-optimal over a slightly larger class of probability measures than the one herein considered.

\medskip

We now assess adaptive recovery rates for empirical Bayes density deconvolution when the errors are ordinary or super-smooth and the mixing density is modeled as a Dirichlet process mixture of Gaussian densities with a data-driven choice for the prior hyper-parameters of the base measure. The problem of deconvoluting a density when the mixing density is modeled as a Dirichlet process mixture of Gaussian densities and the errors are super-smooth has been recently investigated by \citet{sarkar:13}. In a frequentist set-up, rates for density deconvolution have been studied by \citet{carroll:hall:1988}, \citet{fan:1992},
\citet{fan:1991a}, \citet{fan:1991b}. Consider the following model
\[X=Y+\varepsilon,\]
where $Y$ and $\varepsilon$ are independent random variables.
Let $p_Y$ denote the density of $Y$ and $K$ the density
of the error measurement $\varepsilon$. 
The density of $X$ is the convolution $p_X(\cdot)=(K\ast p_Y)(\cdot)=\int
K(\cdot-y)p_Y(y)\d y$. The density
$K$ is assumed to be completely known and its Fourier transform $\widehat{K}$ to satisfy either
\begin{equation}\label{ordinarysmoothkern}
|\widehat{K}(t)|\gtrsim(1+t^2)^{-\eta/2},
\quad t\in\mathbb{R},\qquad \mbox{(ordinary smooth case)}
\end{equation}
for some $\eta>0$, or
\begin{equation}\label{supersmoothkern}
|\widehat{K}(t)|\gtrsim e^{-\varrho|t|^r}, \quad t\in\mathbb{R},\qquad \mbox{(super-smooth case)}
\end{equation}
for some $\varrho,\,r>0$.
The density $p_Y$ is modeled as a Dirichlet process mixture of Gaussian densities as in \eqref{DPMGaussian} with an empirical Bayes choice
$\hga$ of $\gamma$.
Assuming data $\Data=(X_1,\,\ldots,\,X_n)$ are iid  observations from a density
$p_{0X}=K\ast p_{0Y}$ such that the
Fourier transform $\widehat{p_{0Y}}$ of the mixing density $p_{0Y}$ satisfies
\begin{equation}\label{Fourierpoly}
\int_{-\infty}^\infty (1+t^2)^\beta|\widehat{p_{0Y}}(t)|^2\,\d t<\infty \quad \mbox{for some $\beta>1/2$,}
\end{equation}
we derive adaptive empirical Bayes posterior convergence rates for recovering $p_{0Y}$.
\begin{Cor}\label{Mixingdensity}
Suppose that $\hat{K}$ satisfies either condition \eqref{ordinarysmoothkern} (ordinary smooth case) or condition \eqref{supersmoothkern} (super-smooth case) and
that $\widehat{p_{0Y}}$ satisfies the integrability condition \eqref{Fourierpoly}.
Consider a prior for $p_Y$ of the form \eqref{DPMGaussian}, with an empirical Bayes selection $\hga$ for $\gamma$.
Suppose that $p_{0X}=K\ast p_{0Y}$ satisfies the conditions of Theorem \ref{Th:DPMGaussian1} in the ordinary smooth case, with $d$ being the $\mathbb{L}_2$-distance, or those of Theorem \ref{Th:DPMGaussian2} in the super-smooth case.
Then, there exists a sufficiently large
constant $J_1>0$ so that
\[\mathbb{E}^{(n)}_{p_{0X}}[\pi(\|p_Y-p_{0Y}\|_2\geq J_1\epsilon_n\mid \Data ,\,\hga)]=o(1),\]
where, for some constant $\kappa_1>0$,
\[
\epsilon_n=\left\{
  \begin{array}{ll}
        n^{-\beta/{[2(\beta+\eta)+1]}}(\log n)^{\kappa_1}, & \hbox{if $\hat{K}$ satisfies \eqref{ordinarysmoothkern},}\\[5pt]
        (\log n)^{-\beta/r}, & \hbox{if $\hat{K}$ satisfies \eqref{supersmoothkern}.}
  \end{array}
\right.
\]
\end{Cor}

\medskip

The obtained rates are minimax-optimal, up to a logarithmic factor, in the ordinary smooth case and minimax-optimal
in the super-smooth case.
Inspection of the proof of Corollary \ref{Mixingdensity}
shows that, since the result is based on inversion inequalities that relate the $\mathbb{L}_2$-distance between
the true mixing density and the random approximating mixing density to the $\mathbb{L}_2$/$\mathbb{L}_1$-distance between the corresponding mixed densities,
once adaptive rates are known for the direct problem
of Bayes or empirical Bayes estimation of the sampling density $p_{0X}$, the proof can be applied to
yield adaptive recovery rates for either the Bayes or the empirical Bayes density deconvolution problem.
If compared to the approach followed by \citet{sarkar:13},
the present strategy simplifies the derivation
of adaptive recovery rates in a Bayesian density deconvolution problem.
The ordinary smooth case seems to be first treated here even for the fully Bayesian adaptive
density deconvolution problem.

\section{Applications to counting processes with Aalen multiplicative intensity} \label{sec:aalen}
In this section, we illustrate our results on counting processes with Aalen multiplicative intensity. Bayesian nonparametric methods have so far been mainly adopted
to explore possible prior distributions on intensity functions
with the aim of showing that Bayesian nonparametric inference for inhomogeneous Poisson
processes can give satisfactory results in applications, see, \emph{e.g.},
\citet{kottas:sansò:2007}. First frequentist asymptotic behaviors of posteriors like consistency or computations of rates of convergence have been obtained by \citet{belitser:serra:vanzanten} still for inhomogeneous Poisson
processes. As explained in introduction, Theorems~\ref{cor:EB} and \ref{th:gene:aalen} in Section \ref{EB:aalen} extend these results. Section \ref{sec:numerical aalen} illustrates our procedures on artificial data.
\subsection{Notations and setup} \label{subsec:notations}
Let $N$  be a counting process adapted to a filtration $(\mathcal G_t)_t$ with compensator $\Lambda$ so that  $(N_t-\Lambda_t)_t$ is a zero-mean $(\mathcal G_t)_t$-martingale. A counting process satisfies the {\it Aalen multiplicative intensity model} if
$\d\Lambda_t=Y_t\lambda(t)\d t$,
where $\lambda$ is a non-negative deterministic function called, with a slight abuse, the intensity function in the sequel and $(Y_t)_t$ is a non-negative predictable process. Informally,
\begin{equation}\label{informally}
\E[N[t,\,t+\d t]\mid \mathcal G_{t^-}]=Y_t\lambda(t)\d t,
\end{equation}
see \cite{ABGK}, Chapter III. We assume that $\Lambda_t<\infty$ almost surely for every $t$. We also assume that the processes $N$ and $Y$ both depend on an integer $n$ and we consider estimation of $\lambda$ (not depending on $n$) in the asymptotic perspective $n\to\infty$, while $T$ is kept fixed.
The following cases illustrate the interest for this model.
\begin{itemize}
\item[-] {\bf Inhomogeneous Poisson process.} We observe $n$ independent Poisson processes with common intensity $\lambda$. This model is equivalent to the model where we observe a Poisson process with intensity $n\times\lambda$, so it corresponds to the case $Y_t\equiv n$.
\item[-] {\bf Survival analysis with  right-censoring.} This model is popular in biomedical problems. We have $n$ patients and, for each patient $i$, we observe  $(Z_i,\,\delta_i)$ with $Z_i=\min\{X_i,\,C_i\}$, where $X_i$ represents the lifetime of the patient, $C_i$ is the independent censoring time and $\delta_i=1_{X_i\leq C_i}$. In this case, we set $N_t^i=\delta_i\times 1_{Z_i\leq t}$, $Y_t^i=1_{Z_i\geq t}$ and $\lambda$ is the hazard rate of the $X_i$'s: if $f$ is the density of $X_1$,
$\lambda(t)=f(t)/\P(X_1\geq t).$
Then, $N$ (respectively $Y$) is obtained by aggregating the $n$ independent processes $N^i$'s (respectively the $Y^i$'s): for any $t\in [0,\,T]$,
$N_t=\sum_{i=1}^n N^i_t$ and $  Y_t=\sum_{i=1}^n Y^i_t$.
\item[-] {\bf Finite state Markov process.} Let  $X=(X(t))_t$ be a Markov process with finite state space $\mathbb{S}$ and right-continuous sample paths. We assume the existence of  integrable {\it transition intensities} $\lambda_{hj}$ from state $h$ to state $j$ for $h\not= j$. We assume we are given $n$ independent copies of the process $X$, denoted by $X^1,\,\ldots,\,X^n$. For any $i\in\{1,\,\ldots,\,n\}$,
let $N_t^{ihj}$ be the number of direct transitions for $X^i$ from $h$ to $j$  in $[0,\,t]$, for $h\not=j$. Then, the intensity of the multivariate counting process $\frak{N}^i=(N^{ihj})_{h\not=j}$ is $(\lambda_{hj}Y^{ih})_{h\not=j},$ with $Y^{ih}_t=1_{\{X^i(t^-)=h\}}$. As previously, we can consider $\frak{N}$ (respectively $Y^h$) by aggregating the processes $\frak{N}^i$ (respectively the $Y^{ih}$'s): $\frak{N}_t =\sum_{i=1}^n \frak{N}^{i}_t$, $Y_t^h=\sum_{i=1}^n Y^{ih}_t$ and $t\in [0,\,T]$.
The intensity of each component $(N_t^{hj})_t$ of $(\frak{N}_t)_t$  is then $(\lambda_{hj}(t)Y_t^h)_t$.
We refer the reader to \cite{ABGK}, p.~126, for more details. In this case, $N$ is either one of the $N^{hj}$'s or the aggregation of some processes for which the $\lambda_{hj}$'s are equal.
\end{itemize}

\medskip

We denote $\lambda_0$ the true intensity and we define
$\mu_n(t):=\E_{\lambda_0}^{(n)}\left[ Y_t\right]$ and $\tilde\mu_n(t)  :=n^{-1}\mu_n(t)$.
We assume the existence of a non-random
set $\Omega \subset [0,\, T]$ such that there are two constants $m_1,\,m_2$ satisfying for any $n$,
\begin{equation}\label{ass:Y1}
m_1\leq \inf_{t\in\Omega}\tilde \mu_n(t)\leq \sup_{t\in\Omega}\tilde \mu_n(t)\leq m_2
\end{equation}
and $\alpha\in (0,\,1)$ such that if
$\Gamma_n:=\{ \sup_{t\in\Omega}|n^{-1}Y_t-\tilde\mu_n(t)|\leq \alpha m_1\}\cap \{\sup_{t\in\Omega^c}  Y_t = 0\}$,
where $\Omega^c$ is the complementary set of $\Omega$ in $[0,\,T]$, then
\begin{equation}\label{ass:Y2}
\lim_n \P_{\lambda_0}^{(n)}\left(\Gamma_n\right)= 1.
\end{equation}
Assumption \eqref{ass:Y1} implies that  on~$\Gamma_n$,
\begin{equation}\label{ass:Y3}
\forall\, t\in \Omega,\quad (1-\alpha)\tilde\mu_n(t) \leq  \frac{Y_t}{n} \leq (1+\alpha) \tilde\mu_n(t).
\end{equation}
\begin{remark}
Since our results are asymptotic in nature, we can assume, without loss of generality, that \eqref{ass:Y1} is true only for $n$ large enough.
\end{remark}
\begin{remark}
Our assumptions are satisfied, for instance, for the first two illustrative models introduced above. For inhomogeneous Poisson processes, \eqref{ass:Y1} and \eqref{ass:Y2} are obviously satisfied with $\Omega=[0,\,T]$ since $Y_t\equiv\mu_n(t)\equiv n.$ For right-censoring models, with $Y_t^i=1_{Z_i\geq t}$, $i=1,\,\ldots,\,n$, we denote by $\Omega$ the support of the $Z_i$'s and by $M_{\Omega}=\max\Omega\in \overline\R_+.$ Then, \eqref{ass:Y1} and \eqref{ass:Y2} are satisfied if $M_{\Omega}>T$ or $M_{\Omega}\leq T$ and $\P(Z_1=M_{\Omega})>0$ (the concentration inequality is implied by using the DKW inequality).
\end{remark}
Recall that the log-likelihood for Aalen processes is given by,  see \cite{ABGK},
\begin{equation}\label{loglik:aalen}
\ell_n(\lambda)=\int_0^T\log(\lambda(t))\d N_t-\int_0^T\lambda(t)Y_t\d t.
\end{equation}
Since $N$ is empty on $\Omega^c$ almost surely, we only consider estimation over $\Omega$.
So, we set
  $$\mathcal F = \left\{ \lambda : \Omega \rightarrow \R_+\, : \ \int_{\Omega} \lambda(t)\d t < \infty \right\}$$
  endowed with the classical $\L_1$-norm: for all $\lambda,\,\lambda'\in \mathcal F$, let
$\| \lambda - \lambda' \|_1  = \int_{\Omega} | \lambda(t)  -  \lambda'(t)  | \d t$.
We assume that the true intensity $\lambda_0$ satisfies $\lambda_0\in\mathcal F$ and, for any $\lambda \in \mathcal F$,
we write $\lambda  = M_\lambda\bar \lambda$, where $M_\lambda= \int_{\Omega} \lambda(t)\d t$ and $\bar \lambda  \in \mathcal F_1$, with
$\mathcal F_1 = \{ f \in \mathcal F: \ \int_\Omega f(t) \d t = 1\}.$
Note that a prior probability measure $\pi$ on $\mathcal F$ can be written as $\pi_1 \otimes \pi_M$, where $\pi_1$ is a probability distribution on $\mathcal F_1$ and $\pi_M $ is a probability distribution on $\R_+$. This representation will be used in next section.
\subsection{Empirical Bayes and Bayes  posterior concentration rates for monotone non-increasing intensities} \label{EB:aalen}
In this section, we concentrate on estimation of monotone non-increasing intensities, which is equivalent to considering $\bar \lambda $ monotone non-increasing in the above described parameterization. To construct a prior on the set of monotone non-increasing densities on $[0,\,T]$, we use their representation as mixtures of uniform densities as provided by \citet{williamson:56} and we consider a Dirichlet process prior on the mixing distribution:
\begin{equation}\label{DPM:unif}
\bar \lambda(x) = \int_0^\infty \frac{\1_{(0,\,\theta)}(x)}{ \theta} \d P(\theta) , \quad P | G_\gamma,\, A\sim \textrm{DP}(A G_\gamma),
\end{equation}
where $G_\gamma $ is a distribution on $[0,\,T]$. This prior has been studied by \citet{salomond:13} for estimating monotone non-increasing densities. Here, we extend his results to the case of a monotone non-increasing intensity of an Aalen process with a data-dependent $\gamma$. We denote by $\pi(\cdot\mid N)$ the posterior distribution
given the observations of the process $N$.

We study the family of distributions  $G_\gamma$, with density denoted $g_\gamma$, belonging to one of  the following families of densities with respect to the Lebesgue measure: for $ \gamma >0$, $ a> 1$,
 \begin{equation}\label{base:mono}
 \begin{split}
 g_\gamma (\theta) &= \frac{1}{ G(T\gamma ) } \gamma^a  \theta^{a-1} e^{-\gamma \theta} \1_{\{0\leq\theta \leq T\}}  \quad
 \mbox{ or }  \quad \left(\frac{ 1 }{ \theta } -  \frac{ 1 }{T}\right)^{-1} \sim\Gamma( a,\, \gamma),
 \end{split}
\end{equation}
where $G$ is the cumulative distribution function of a $\Gamma(a,\,1)$ random variable.
Assume that, with probability going to 1,  $\hga$ belongs to a  fixed compact  subset of $(1,\, \infty)$ denoted by $\mathcal K$.
We then have the following theorem, which is an application of Theorem \ref{th:gene:EB}.
\begin{Th}\label{cor:EB}
Let $\bar \epsilon_n = (n/\log n)^{-1/3}$.  Assume that the prior $\pi_M$ on the mass $M$  is absolutely continuous with respect to the Lebesgue measure with positive and continuous density on $\R_+$ and that it has a finite Laplace transform in a neighbourhood of 0. Assume that the prior $\pi_1(\cdot \mid \gamma) $ on $\bar \lambda$ is  a Dirichlet process mixture of uniform distributions defined by \eqref{DPM:unif}, with $ A>0$ and the base measure  $G_\gamma$ defined by \eqref{base:mono}.  Let $\hga$ be a measurable function of the observations satisfying $\P_{\lambda_0}^{(n)}( \hga\in \mathcal K ) = 1+ o(1)$ for some fixed compact subset  $\mathcal K \subset (1,\,\infty)$. Assume also that (\ref{ass:Y1}) and (\ref{ass:Y2}) are satisfied and for any $k\geq 1$  there exists $C_{1k}>0$ such that
\begin{equation}\label{moment}
\E_{\lambda_0}^{(n)}\left[ \left(  \int_\Omega ( Y_t-\mu_n(t))^2\d t \right)^k  \right] \leq C_{1k} n^k
\end{equation}
Then, there exists $J_1>0$ such that
\begin{equation*}\label{eq:EB:aalen}
\E_{\lambda_0}^{(n)}[\pi(\lambda: \ \| \lambda - \lambda_0\|_{1} >  J_1 \bar \epsilon_n \mid N,\, \hga)] = o(1)
\end{equation*}
and
\begin{equation*}\label{eq:FB:Aalen}
\sup_{\gamma \in \mathcal K} \E_{\lambda_0}^{(n)}[\pi(\lambda: \ \| \lambda - \lambda_0\|_{1} >  J_1 \bar \epsilon_n \mid N,\,\gamma)] = o(1).
\end{equation*}
\end{Th}
The proof of Theorem \ref{cor:EB} is given in Section \ref{sec:pr:corEB}. It consists in verifying conditions [A1] and [A2] of Theorem \ref{th:gene:EB} and is based on  a general theorem on posterior concentration rates for Aalen processes which is presented below since it is of interest on its own.
Let $v_n$ be a positive sequence going to 0 such that $n v_n^2\rightarrow\infty$.  For all $j \geq 1$, we define
$$\bar S_{n,j} = \left\{ \bar \lambda \in \mathcal F_1:\,\| \bar \lambda - \bar \lambda_0 \|_1 \leq \frac{2(j+1)v_n}{M_{\lambda_0} }\right\},$$
where
$M_{\lambda_0}=\int_{\Omega}\lambda_0(x)\d x$, $\bar\lambda_0=M_{\lambda_0}^{-1}\lambda_0$.
For $H>0$ and $k\geq 2$,  if $k_{[2]} = \min \{ 2^\ell:\, \ell\in \mathbb N,\, 2^\ell \geq k\}$, we define
\begin{small}
\begin{equation*}
\begin{split}
\bar B_{k,n}(H,\, \bar \lambda_0,\,v_n)=&  \left\{ \bar \lambda \in \mathcal F_1:\, \ h^2(\bar \lambda_0,\, \bar \lambda)
\leq \frac{v_n^2}{1 + \log  \| \frac{ \bar \lambda_0 }{ \bar \lambda } \|_\infty }, \, \max_{2\leq j \leq k_{[2]} }E_j(\bar \lambda_0 ,\, \bar \lambda) \leq v_n^2, \
\left\| \frac{ \bar \lambda_0 }{ \bar \lambda } \right\|_\infty \leq n^H,\, \left\|\bar \lambda\right\|_\infty \leq H \right\},
\end{split}
\end{equation*}
\end{small}
where, for every integer $j$,
$E_j(\bar \lambda_0,\, \bar \lambda) = \int_{\Omega} \bar \lambda_0(x) \left| \log \bar \lambda_0  (x) - \log \bar \lambda (x)  \right|^j\d x$ and
$\|\cdot\|_\infty$ stands for the sup-norm.
\begin{Th}\label{th:gene:aalen}
Let $v_n$ be a positive sequence satisfying $v_n^2 \geq \log n/n$ and $v_n = o(1)$. Assume that \eqref{ass:Y1} and  \eqref{ass:Y2}  are satisfied and that the prior $\pi_M$ on the mass $M$  is absolutely continuous with respect to the Lebesgue measure with positive and continuous density on $\R_+$. Assume also that \eqref{moment} is satisfied for some $k\geq 2$. Finally, assume that the prior $\pi_1$ on $\bar \lambda$ satisfies the following assumptions for some $H>0$.
\begin{itemize}
\item[$(i)$] There exists $\mathcal F_n \subset \mathcal F_1$ such that
$$\pi_1\left( \mathcal F_n^c \right) \leq e^{ -(\kappa_0+2) n v_n^2}\pi_1(\bar B_{k,n}(H,\, \bar \lambda_0,\,v_n)),  $$
with $\kappa_0$ as defined in \eqref{kappa0}, and for all $\xi,\,\delta>0$,
$$\log D( \xi,\,\mathcal F_n,\,\| \cdot \|_1 )\leq n \delta \quad\mbox{ for all $n$ large enough.}$$

\item[$(ii)$] For all $\zeta,\,\delta,\,\beta>0$, there exists $J_0>0$ such that, for all $j$ such that $J_0\leq j\leq \beta/v_n$,
$$\frac{ \pi_1(\bar S_{n,j}) }{ \pi_1(\bar B_{k,n}(H,\, \bar \lambda_0,\,v_n) )} \leq e^{\delta (j+1)^2 n v_n^2 } $$
and
$$ \log D( \zeta jv_n,\,\bar S_{n,j} \cap \mathcal F_n,\,\| \cdot \|_1 ) \leq \delta (j+1)^2 n v_n^2 .$$
\end{itemize}
Then, there exists $J_1>0$ such that, when $n\rightarrow\infty$,
$$ \E_{\lambda_0}^{(n)}[\pi(\lambda: \ \| \lambda - \lambda_0\|_{1} >  J_1 v_n \mid N)] = O((n v_n^2)^{-k/2}).$$
\end{Th}
To the best of our knowledge, the only other paper dealing with posterior concentration rates in related models is \citet{belitser:serra:vanzanten} which considers inhomogeneous Poisson processes. Theorem \ref{th:gene:aalen} differs from their general Theorem 1 in two aspects. First, we do not restrict ourselves to inhomogeneous Poisson processes. Second, and more importantly, our set of conditions is quite different. More specifically, we do not need to assume that $\lambda_0$ is bounded from below and we do not need to bound from below the prior mass of neighborhoods of $\lambda_0$ for the sup-norm, but merely the  prior mass of neighborhoods for the Hellinger distance, as in Theorem 2.2 of \citet{ghosal:ghosh:vdv:00}. Our aim, in Theorem  \ref{th:gene:aalen}, is to propose a set of conditions to derive posterior concentration rates on the intensity which are as close as possible to the conditions used in the density model by parameterizing $\lambda $ as $\lambda=M_\lambda \bar \lambda$, where $\bar \lambda$ is a probability density on $\Omega$.
\begin{remark}
Note that if $\bar \lambda \in \bar B_{2,n}(H,\, \bar \lambda_0,\,v_n)$, then, for any $j > 2$,
$ E_j(\bar \lambda_0 ,\, \bar \lambda) \leq H^{j-2} v_n^2 (\log n)^{j-2}$
so that, using Proposition \ref{prop:KL:Aalen}, if we replace $\bar B_{k,n}(H,\, \bar \lambda_0,\,v_n)$ with $\bar \lambda \in \bar B_{2,n}(H,\, \bar \lambda_0,\,v_n)$ in the assumptions of Theorem \ref{th:gene:aalen},
we obtain the same type of conclusion: for any $k\geq 2$, such that condition \eqref{moment} is satisfied,
$$ \E_{\lambda_0}^{(n)} [\pi(\lambda: \ \| \lambda - \lambda_0\|_{1} \geq  J_1 v_n \mid N)] = O((n v_n^2)^{-k/2} (\log n)^{k(k_{[2]}-2)/2} )$$
with an extra $(\log n)$-term on the right-hand side above.
\end{remark}

\begin{remark}
We now prove that Assumption~\eqref{moment} is reasonable by considering the previous examples. Obviously, \eqref{moment} is satisfied in the case of  inhomogeneous Poisson processes since $Y_t=n$ for every $t$. For the censoring model, $Y_t=\sum_{i=1}^n1_{Z_i\geq t}.$ We set, for any $i$,
$V_i=1_{Z_i\geq t}-\P(Z_1\geq t).$
Then, for $k\geq 1$,
\begin{eqnarray*}
\E_{\lambda_0}^{(n)}\left[ \left(  \int_\Omega ( Y_t-\mu_n(t))^2 \mathrm{d} t \right)^k  \right] &\leq&\E_{\lambda_0}^{(n)}\left[ \left(  \int_0^T \left(\sum_{i=1}^nV_i\right)^2\mathrm{d} t \right)^k  \right]\\ &\leq& T^{k-1} \int_0^T\E_{\lambda_0}^{(n)}\left[  \left(\sum_{i=1}^nV_i\right)^{2k}\right] \mathrm{d} t \\
&\leq&C(k,\,T)\int_0^T\left(\sum_{i=1}^n \E_{\lambda_0}^{(n)}[V_i^{2k}]+\left(\sum_{i=1}^n\E[V_i^2]\right)^k\right)\mathrm{d} t\\
&\leq& C_{1k}n^k,
\end{eqnarray*}
where $C(k,\,T)$ and $C_{1k}$ depend on $k$ and $T$ by using the H\"{o}lder and Rosenthal inequalities (see, for instance, Theorem~C.2 of \cite{HKPT}). Under some mild conditions, similar computations can be derived for finite state Markov processes.
\end{remark}


\subsection{Numerical illustration}\label{sec:numerical aalen}

We now present a numerical experiment to  highlight the impact of  an empirical prior distribution for finite sample sizes $n$, in the case of an inhomogeneous Poisson process.   Let $(W_{i})_{i= 1 \dots N(T)}$ be  the  observed points of the process $N$ over $[0,T]$, $N(T)$ being the observed number of jumps over $[0,T]$. We recall that the intensity function has the form $n \lambda_0(t) = n M_{\lambda_0} \overline{\lambda_0}(t)  $ ($n$ being known) where $\int_{0}^T \overline{\lambda}(t)dt=1$. \\

\noindent The estimations of $M_{\lambda_0}$ and $\overline{\lambda_0}$ can be done separately, given the factorisation in \eqref{loglik:aalen}. Provided the use of a Gamma prior distribution on $M_{\lambda_0}$ ($M_{\lambda_0} \sim \Gamma(a_M,b_M)$),   we have  
$$M_{\lambda_0} |N \sim \Gamma(a_M+N(T),  b_M + n). $$
The non-parametric Bayesian estimation of $\overline{\lambda}$ is more involved. However in the case of Dirichlet process mixtures of uniforms as a prior model on $\bar \lambda$ we can use the same algorithms as those considered for density estimation. In this Section we restrict ourselves to the case where the base measure in the Dirichlet process is the second possibility in \eqref{base:mono}, namely $G_{\gamma} =_{\mathcal{D}} \left[\frac{1}{T} + \frac{1}{ \Gamma(a,\gamma)}\right]^{-1}$. It satisfies the assumptions of Theorem \ref{cor:EB} and also presents computational advantages. 
Hence, three hyperparameters are involved in this prior, namely $A$ the mass of the Dirichlet process, $a$ and $\gamma$.   The hyperparameter $A$ strongly conditions the number of classes in the posterior distribution of $\overline{\lambda}$. In order  to mitigate its  influence on  the posterior,  we propose to consider a hierarchical approach and set a Gamma prior distribution on $A$: $A \sim \Gamma(a_A, b_A)$. In the absence of additional information, we set $a_A = b_A = 1/10$, corresponding to a weakly informative prior. Theorem \ref{cor:EB} applies for any $a>1$. We arbitrarily set $a= 2$, the influence of $a$ is not studied in this paper. 
We compare three strategies for the determination of $\gamma$ in our simulation study. 
\begin{enumerate}
 \item[] \emph{Strategy 1: Empirical Bayes -} 
We propose the following empirical estimator: 
\begin{equation}\label{eq:formula gamma emp}
\widehat{\gamma}_n = \Psi^{-1} \left[\overline{W}_{N(T)} \right]
\end{equation}
which comes from the following equality:
  \begin{eqnarray*}
 \mathbb{E}\left[\overline{W}_{N(T)} | N(T)\right]  = \mathbb{E}\left[  \left. \frac{1}{N(T)} \sum_{i=1}^{N(T)} W_i \;  \right\rvert N(T) \right]
 &=&  \int_{0}^T t\int_{\theta}  \frac{1}{\theta} \ind_{[0, \theta]}(t)dG_{\gamma}(\theta) \; dt
 \\
& =& \frac{\gamma^a}{2\Gamma(a)}  \int_{\frac{1}{T}}^{\infty}   \frac{ e^{-\gamma/(\nu-\frac{1}{T})}}{(\nu-\frac{1}{T})^{(a+1)}}      \frac{1}{\nu}  d\nu  := \Psi(\gamma)
 \end{eqnarray*} 
\item[] \emph{Strategy 2:   Fixed $\gamma$ -} In order to avoid the empirical prior, one can fix $\gamma=\gamma_0$. In order to study the impact of a bad choice of $\gamma$ on the behaviour of the posterior distribution, we propose to choose $\gamma_0$  far from a well calibrated  value for $\gamma$, which would be 
$\gamma^* =   \Psi^{-1}(E_{theo})$ with   $E_{theo} = \mathbb{E}\left[W | N(T)\right]$. Since,  in the simulation study, we know the true $\lambda_0$, $\gamma^*$ can be computed  and  we consider
\begin{equation}\label{eq:formula gamma 0}
\gamma_0 = \rho  \cdot  \Psi^{-1}(E_{theo}), \quad \mbox{ where } \quad E_{theo} = \int_0^T t \overline{\lambda}(t) dt, \quad \rho \in \{ 0.01, 30, 100\}
 \end{equation}

  \item[] \emph{Strategy 3: Hierarchical Bayes -} We consider a hierarchical prior on $\gamma$ given by:
 $$  \gamma \sim \Gamma(a_\gamma, b_\gamma)$$
In order to make a fair comparison with the empirical posterior, we center the prior distribution on $\widehat{\gamma}$.  Besides, in the simulation study, we set two different hierarchical hyperparameters $(a_{\gamma},b_{\gamma})$, corresponding to two  prior variances. More precisely, 
$(a_{\gamma},b_{\gamma})$ are such that   the prior expectation is equal to $ \widehat{\gamma}$ and the prior variance is $ \sigma_{\gamma}^2$, the values of $\sigma_{\gamma}$ being specified in Table \ref{tab:num_comput_times}.   
 \end{enumerate}


\noindent Samples from the posterior distribution of $(\overline{\lambda}, A, \gamma)| \Nb$  are generated  using a Gibbs algorithm, decomposed into two or three steps: 
$$[1] \quad  \overline{\lambda} | A, \gamma, \Nb \quad \quad [2] \quad    A | \overline{\lambda}, \gamma, \Nb   \quad\quad [3] ^\dagger  \quad \gamma | A, \overline{\lambda}, \Nb,  $$
step $ [3]^\dagger$ only existing  if  a fully Bayesian strategy is adopted on the hyperparameter $\gamma$ (strategy 3).   
 We use the algorithm developed by \cite{fall2012} whose details are given in  Appendix  \ref{sec:details algo}. 

\begin{figure}[htpb]
 \begin{center}
  \includegraphics[width = \linewidth]{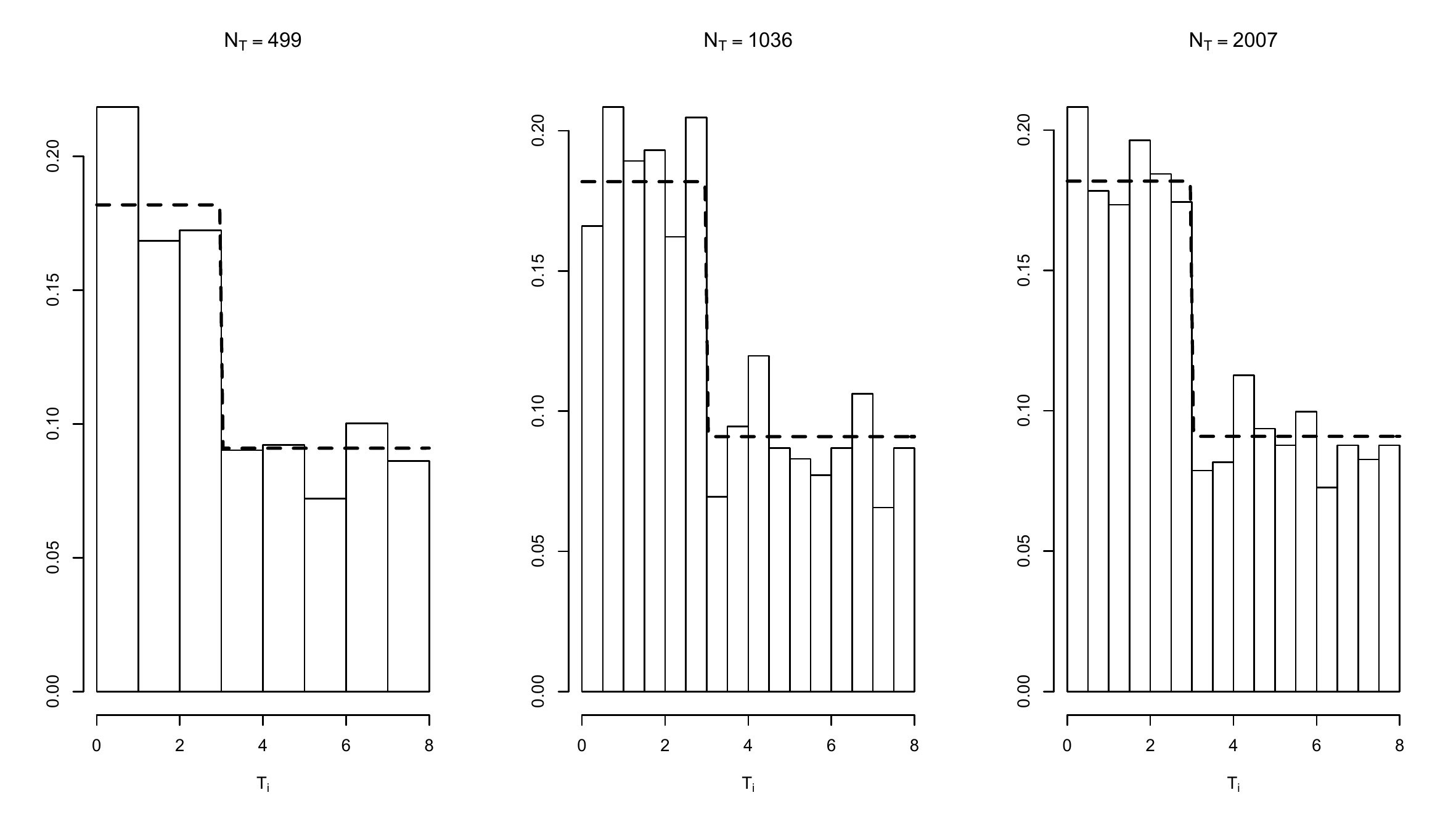}\\
 \includegraphics[width = \linewidth]{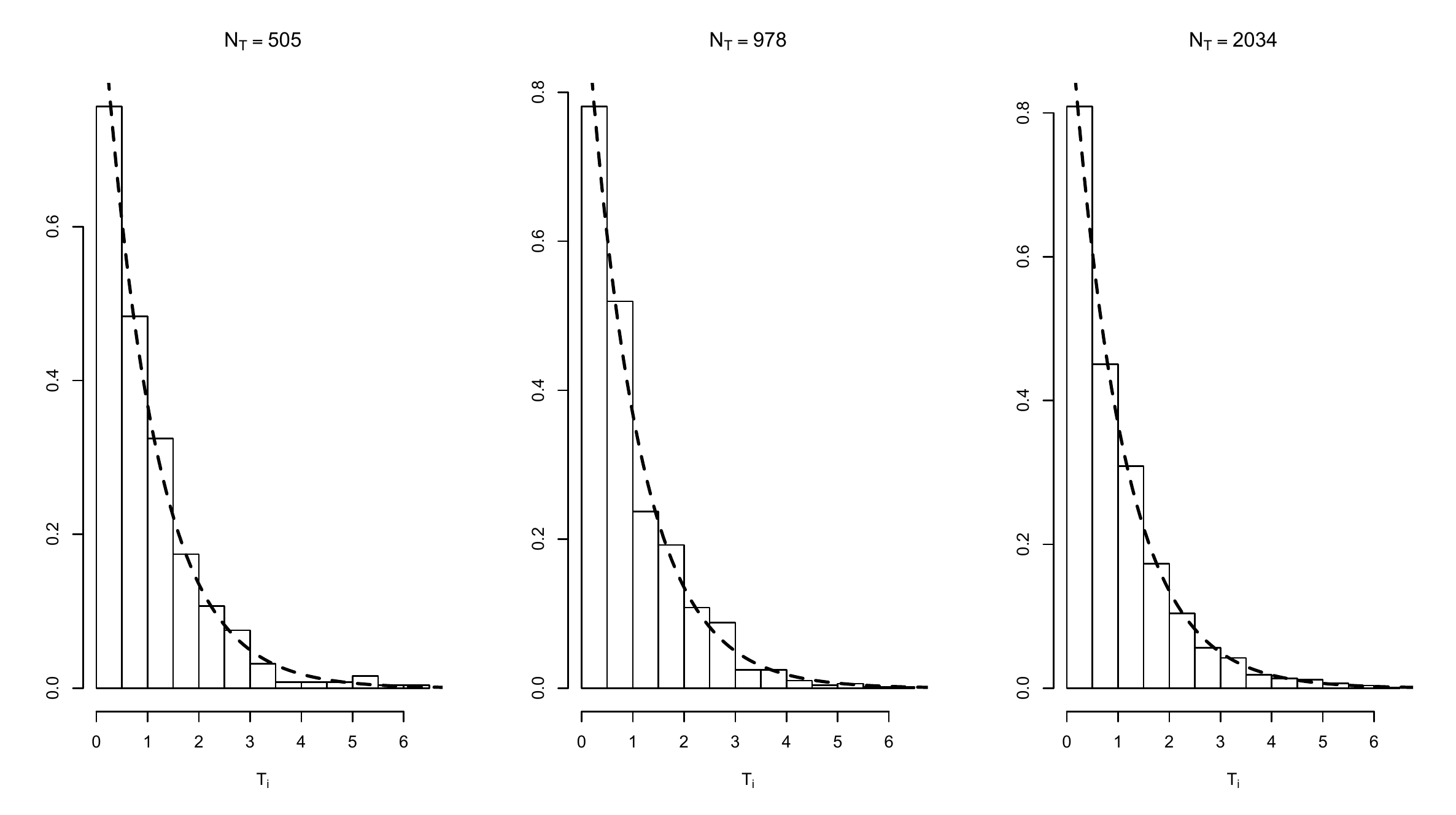}\\
  \includegraphics[width = \linewidth]{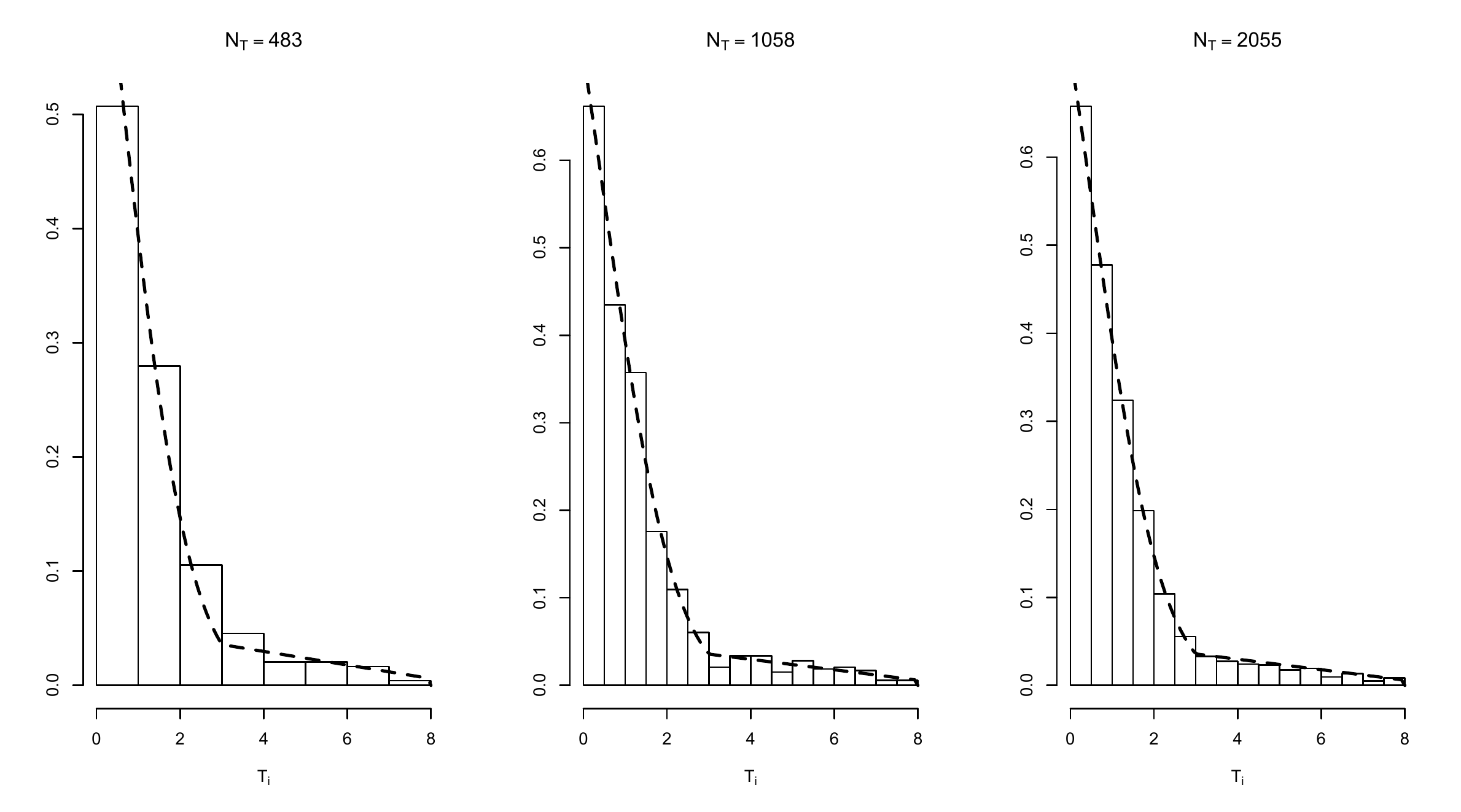}\\
 \end{center}
   \caption{\emph{Simulation study on the Poisson process}: simulated data from intensity  $n\lambda_{0,1}$  (first line), $n \lambda_{0,2}$ (second line), $n \lambda_{0,3}$  (third line), with $n=500$ (first column), $n=1000$ (second column), $n=2000$ (third column). }
\label{fig:data aalen}
\end{figure}

\noindent The various strategies for calibrating the hyperparameter $\gamma$ are tested on $3$ different intensity functions (non null over $[0,T]$, with $T=8$): 
\begin{eqnarray*}
 \lambda_{0,1}(t) &=&  \left[ 4 \; \ind_{[0,3]}(t) + 2\; \ind_{[3,8]}(t) \right] \\
\lambda_{0,2}(t) &=& e^{-0.4 t}\\
\lambda_{0,3}(t) &=&\left[ \cos^{-1}\Phi(t)\ind_{[0,3]}(t)  - ( \frac{1}{6} \cos^{-1}\Phi(3) t - \frac{3}{2} \cos^{-1}\Phi(3)) \ind_{[3,8]}(t) \right] 
 \end{eqnarray*} 
 where $\Phi$ is the probability function of the standard normal distribution. 
 For each intensity $\lambda_{0,1}$, $\lambda_{0,2}$ and $\lambda_{0,3}$,  we simulate $3$ datasets corresponding to $n=500, 1000,2000$ respectively.  The histograms of the 9 datasets  with the true corresponding normalized intensities  $\overline{\lambda_{0,i}}$ are plotted on Figure \ref{fig:data aalen}.   In the following we denote by $\mathbf{D}^i_{n}$ the dataset associated with $n$ and intensity $ \lambda_{0,i}$. 
 For each dataset, we adopt the 3 exposed strategies to calibrate $\gamma$.  The posterior distributions are sampled 
 using $30000$ iterations of which are removed a burn-in period of $15000$ iterations.
 
\begin{figure}
\centering
\begin{tabular}{ccc}
$D^1_{500}$ & $D^1_{1000}$& $D^1_{2000}$ \\
\includegraphics[width=0.33\textwidth]{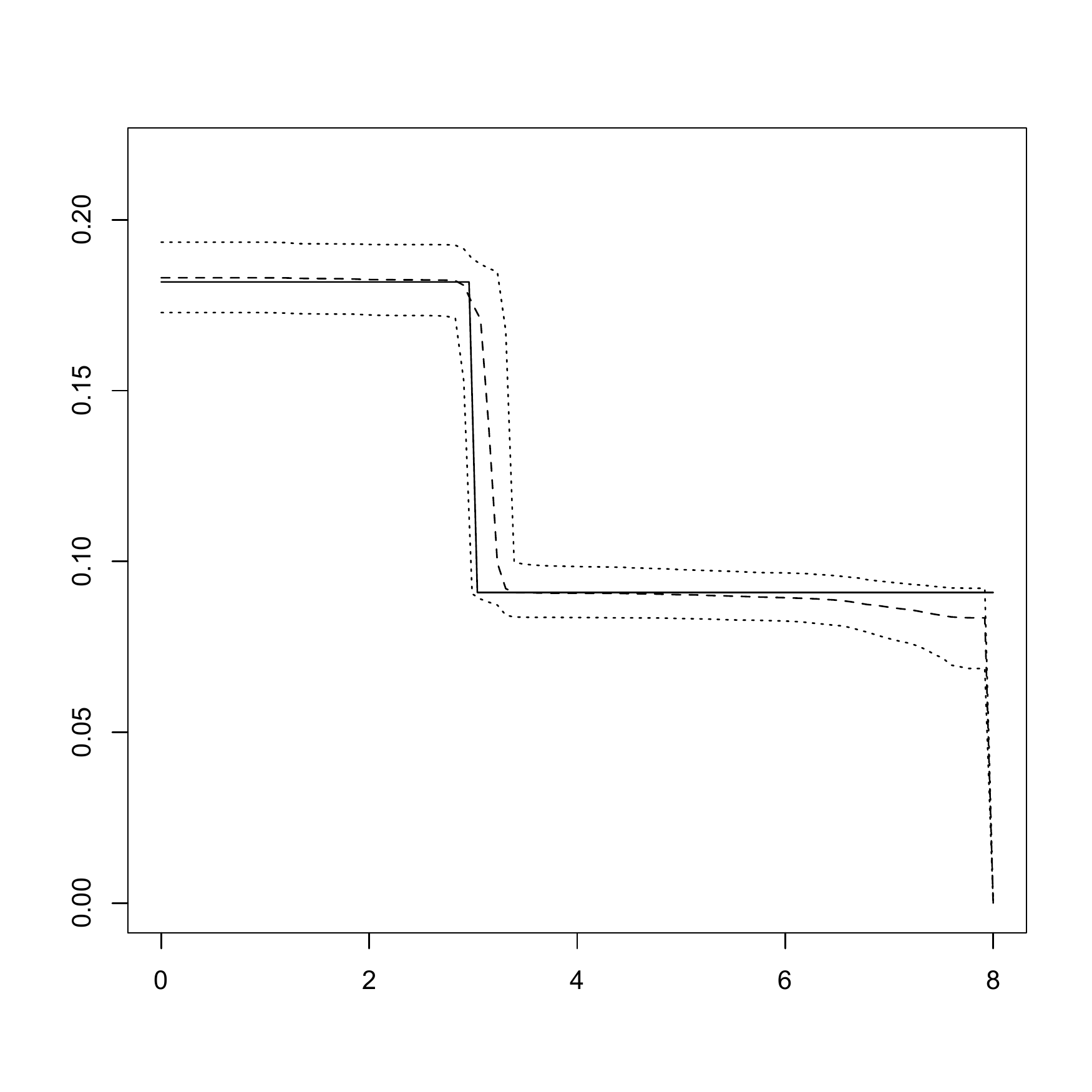} &
\includegraphics[width=0.33\textwidth]{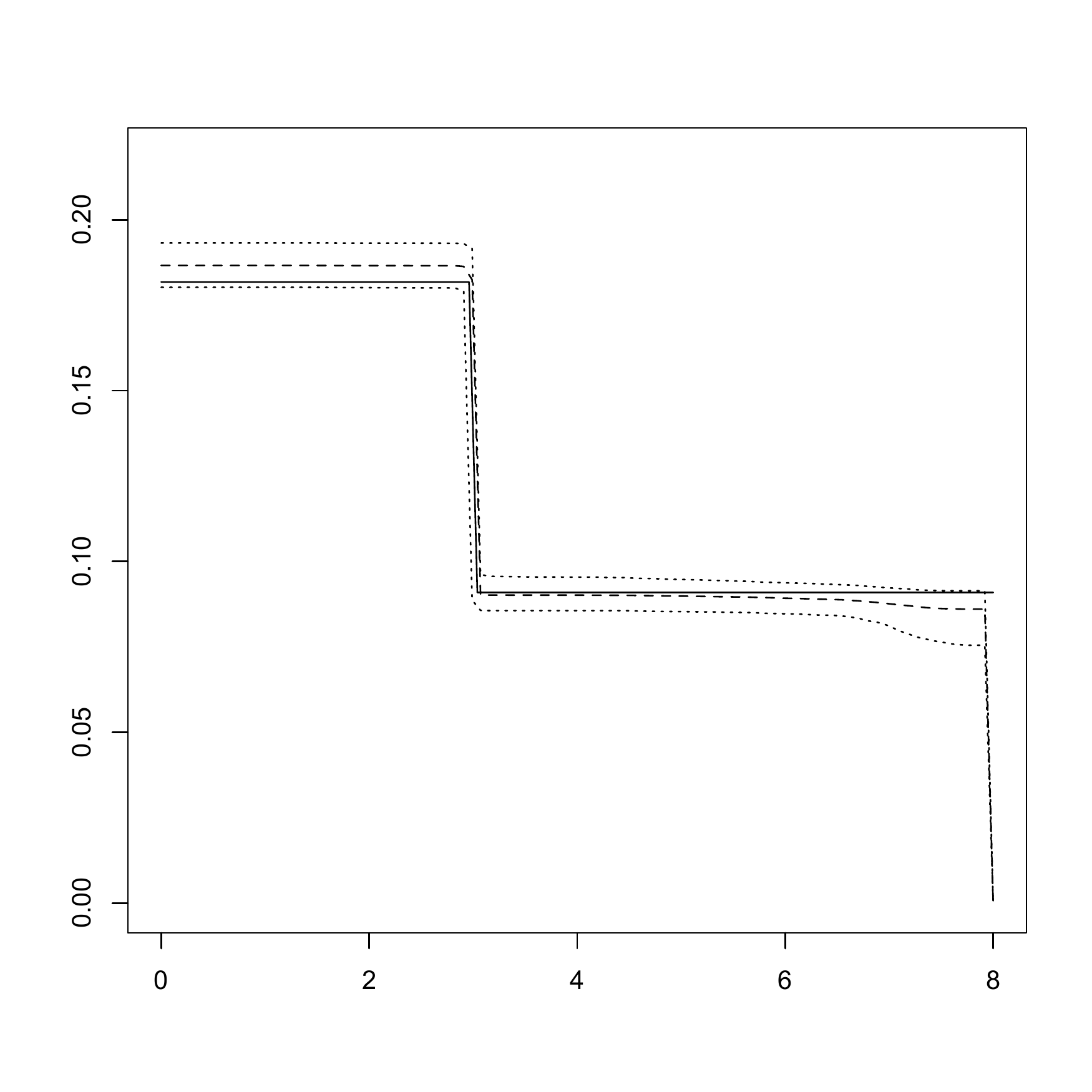} &
\includegraphics[width=0.33\textwidth]{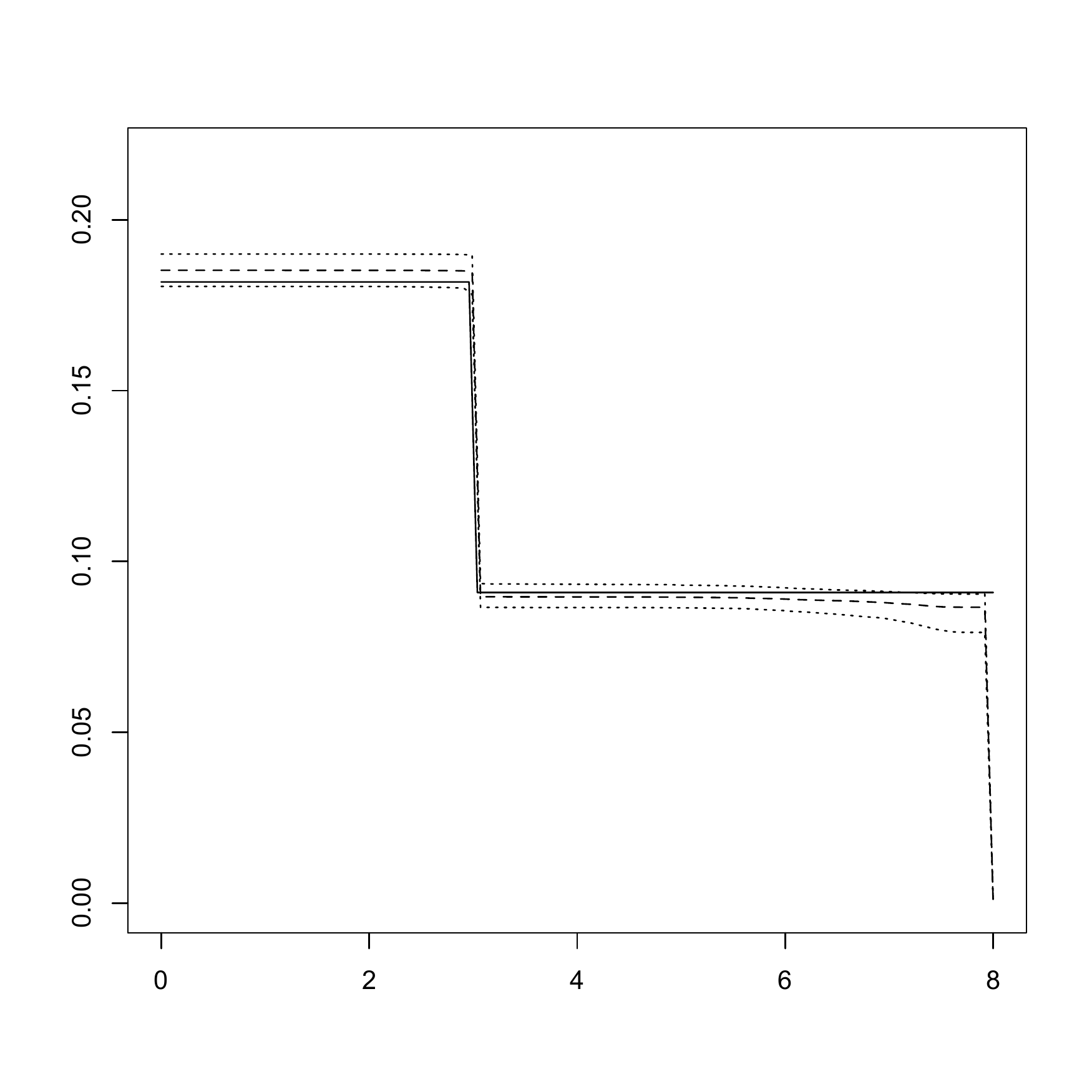} \\\includegraphics[width=0.33\textwidth]{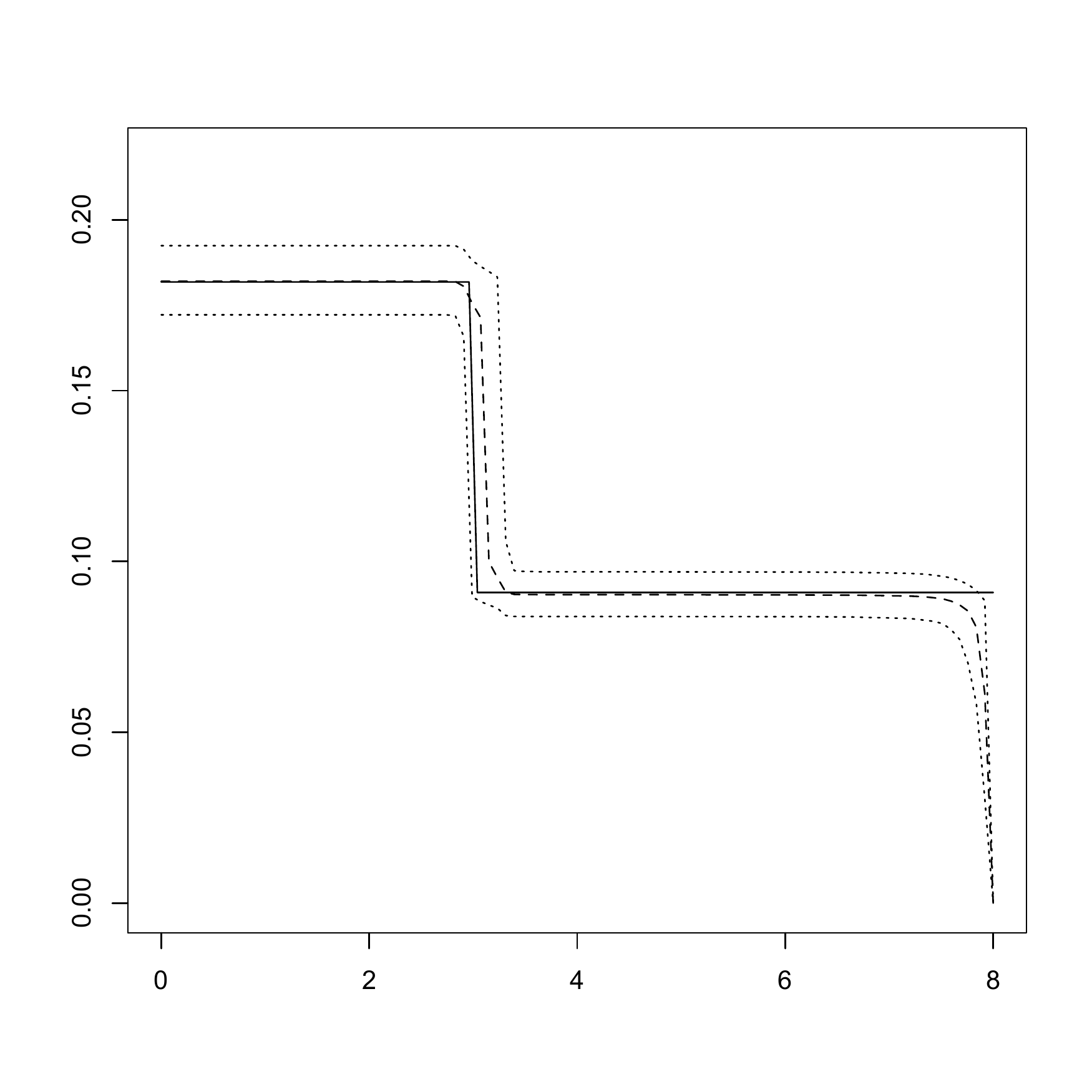} &
\includegraphics[width=0.33\textwidth]{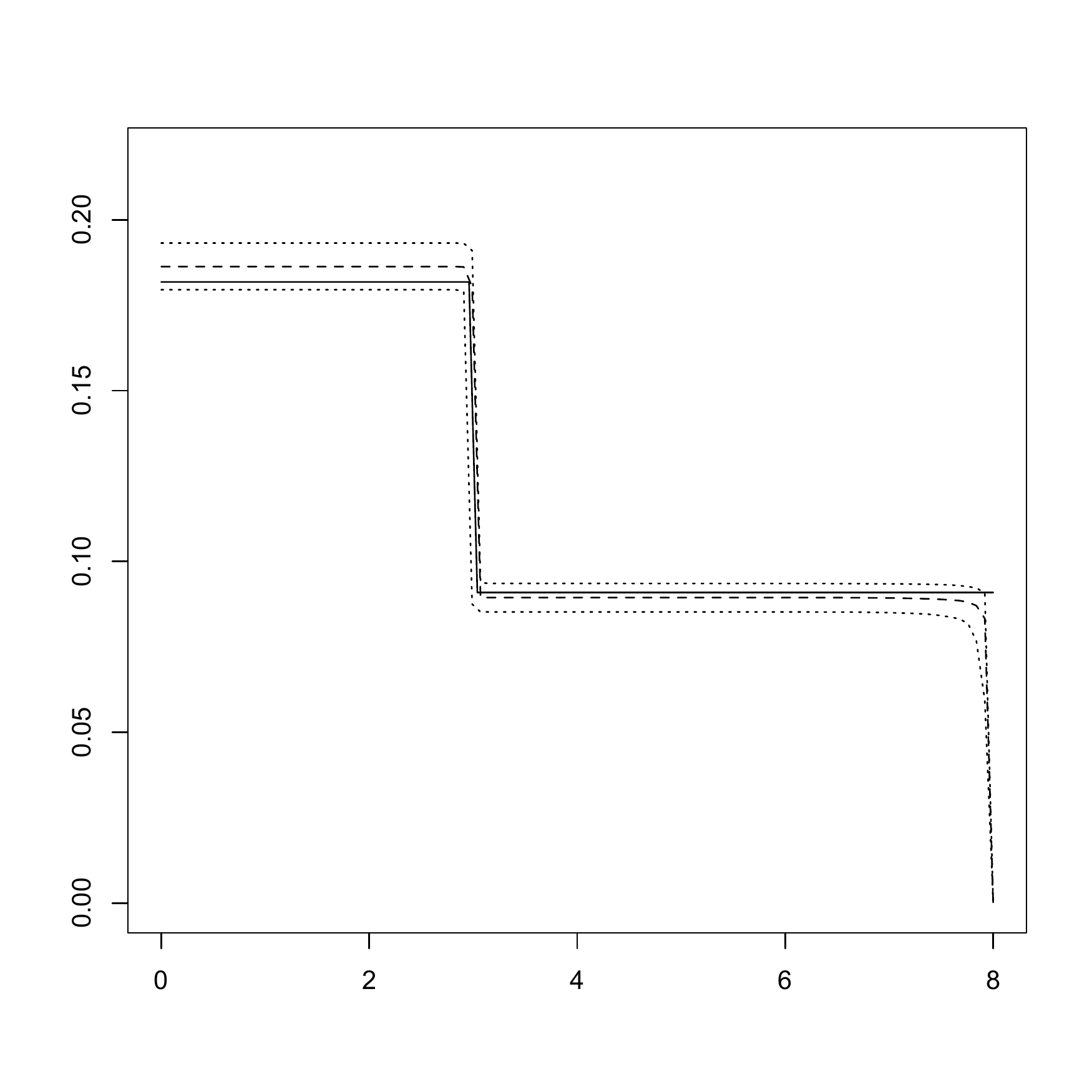} &
\includegraphics[width=0.33\textwidth]{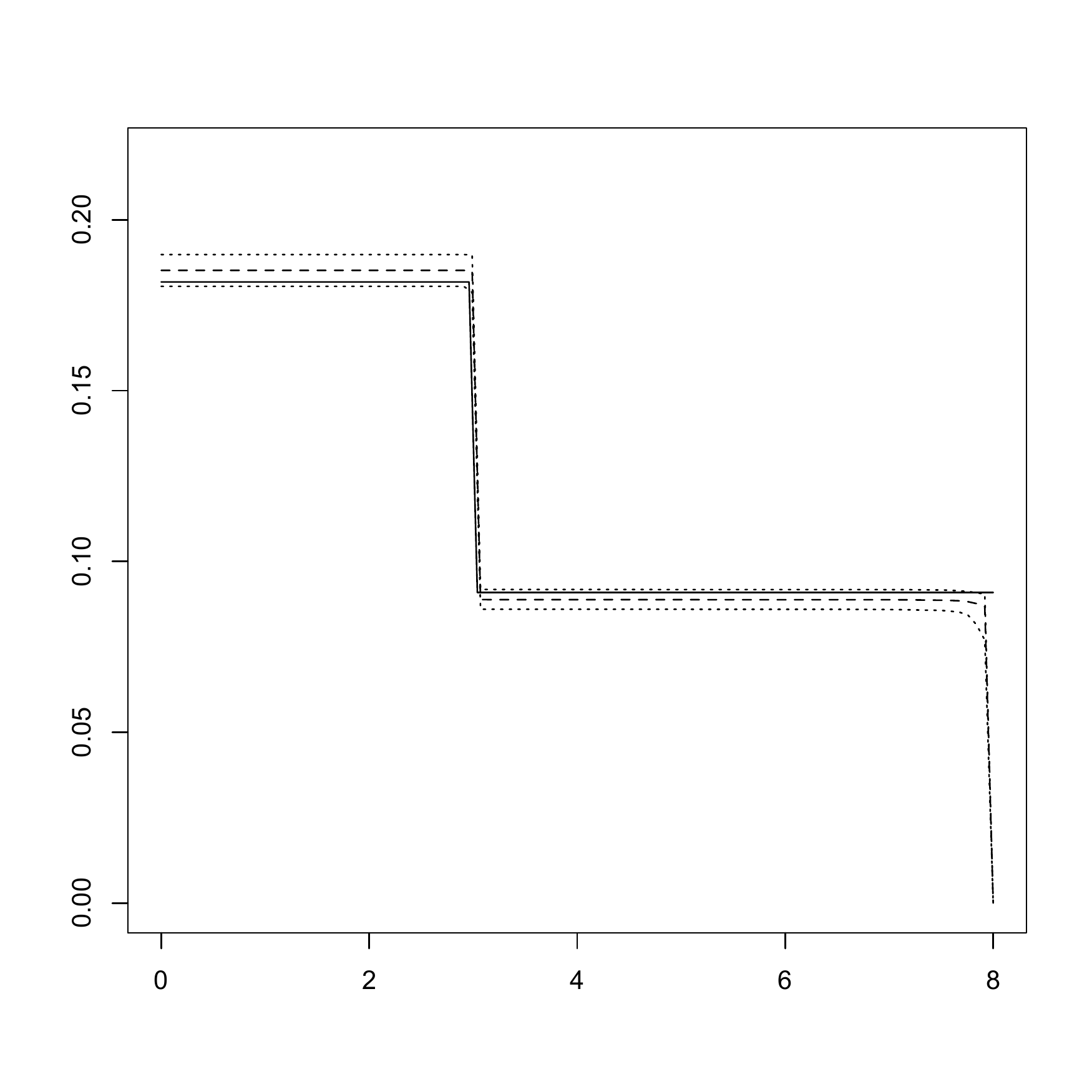} \\
\includegraphics[width=0.33\textwidth]{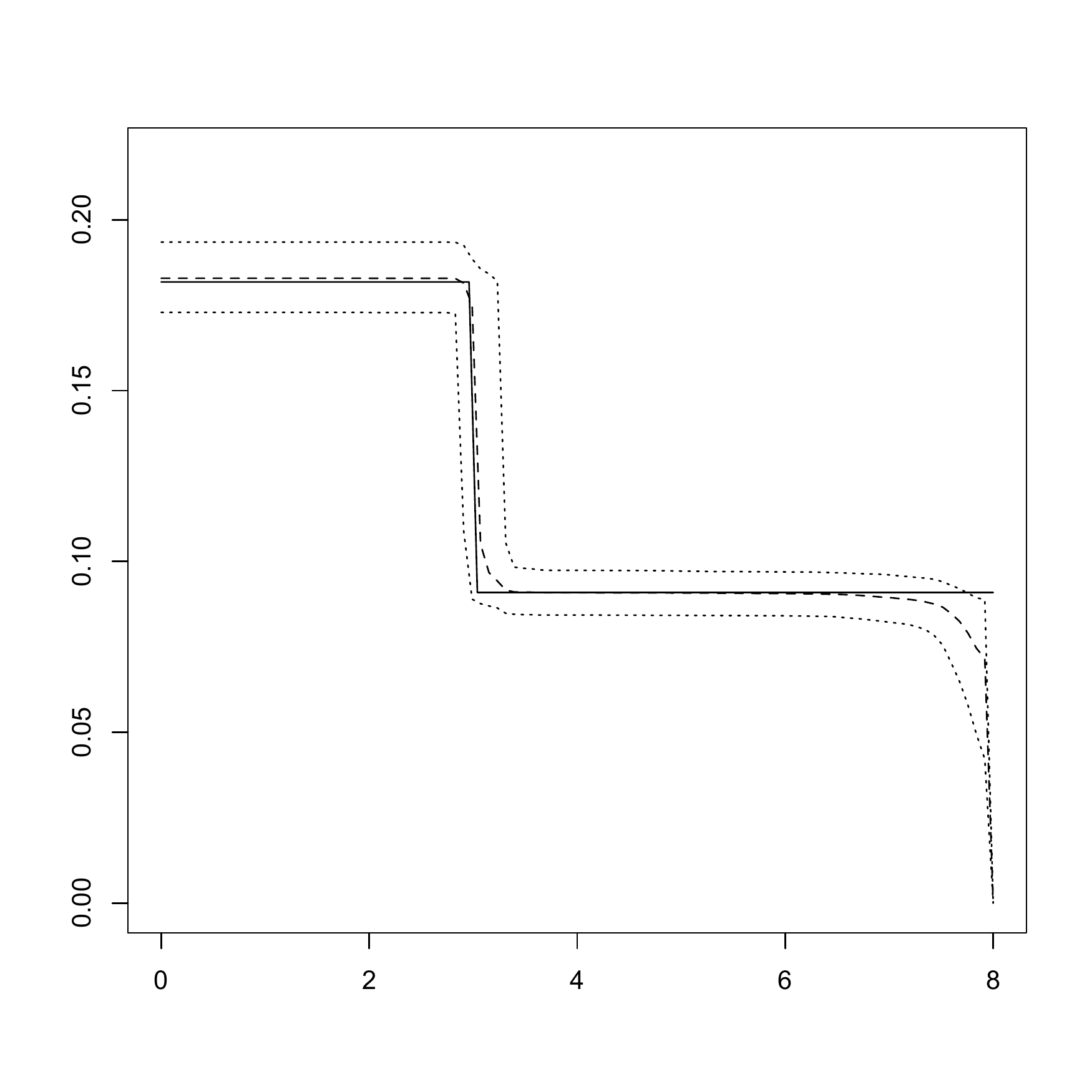} &
\includegraphics[width=0.33\textwidth]{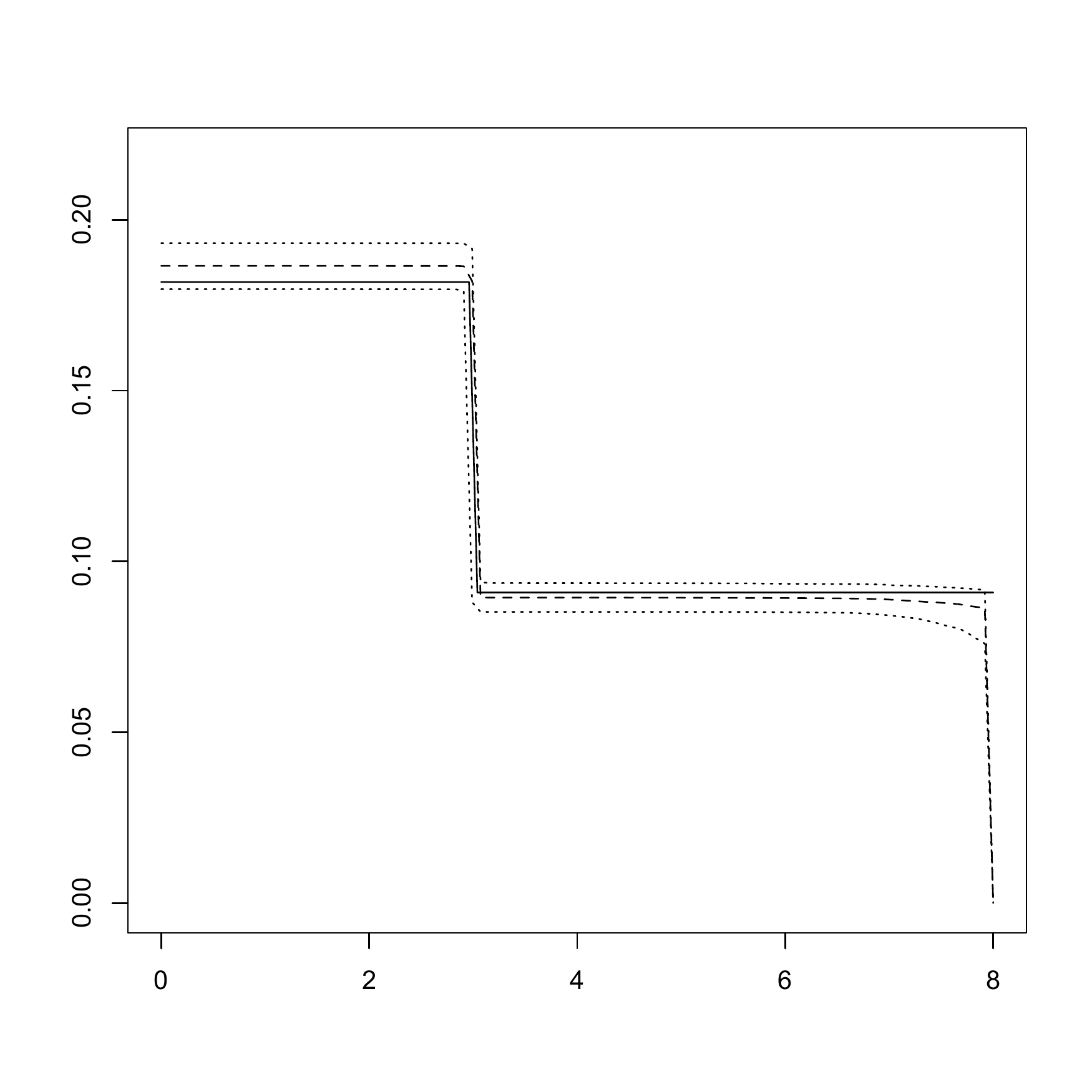} &
\includegraphics[width=0.33\textwidth]{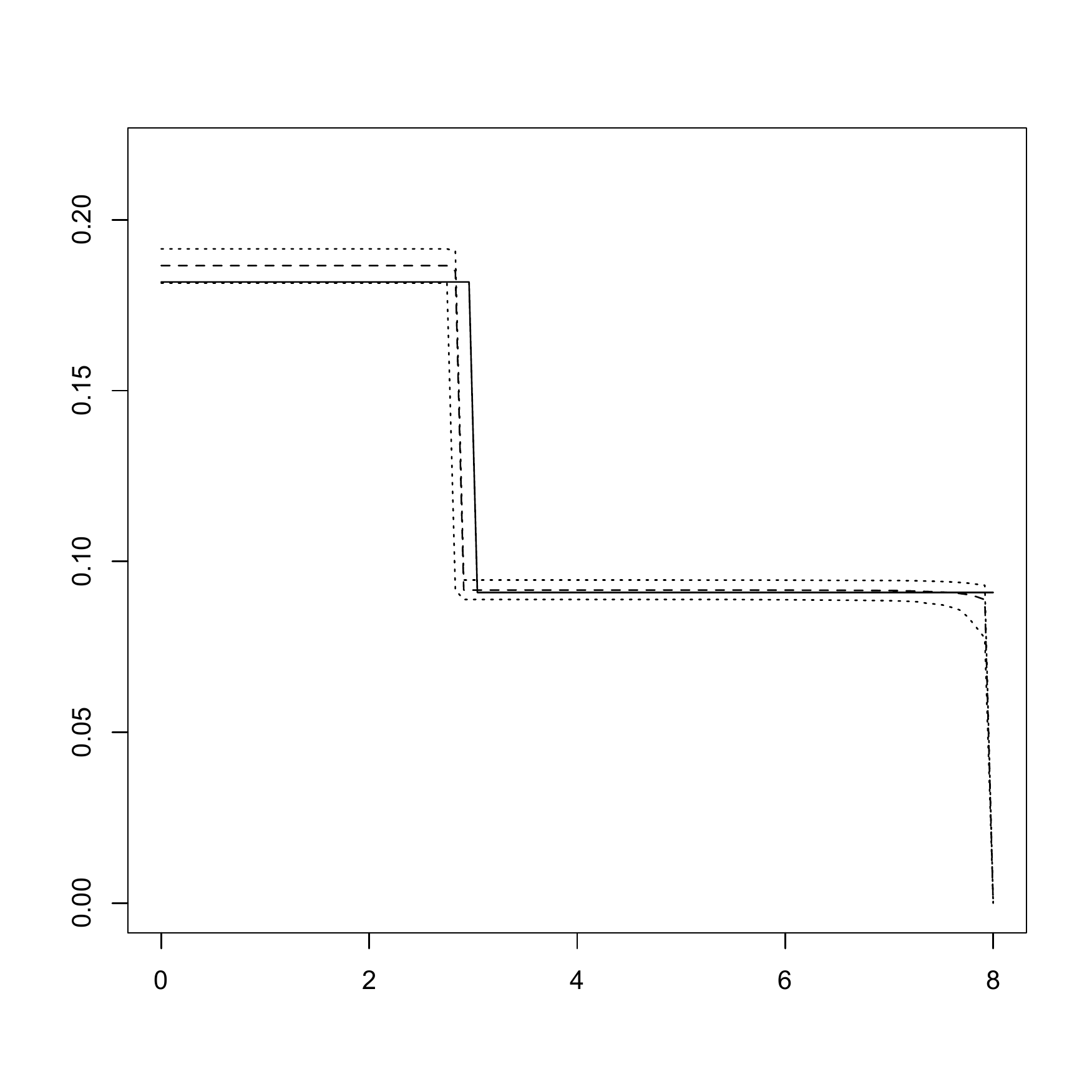} \\
\includegraphics[width=0.33\textwidth]{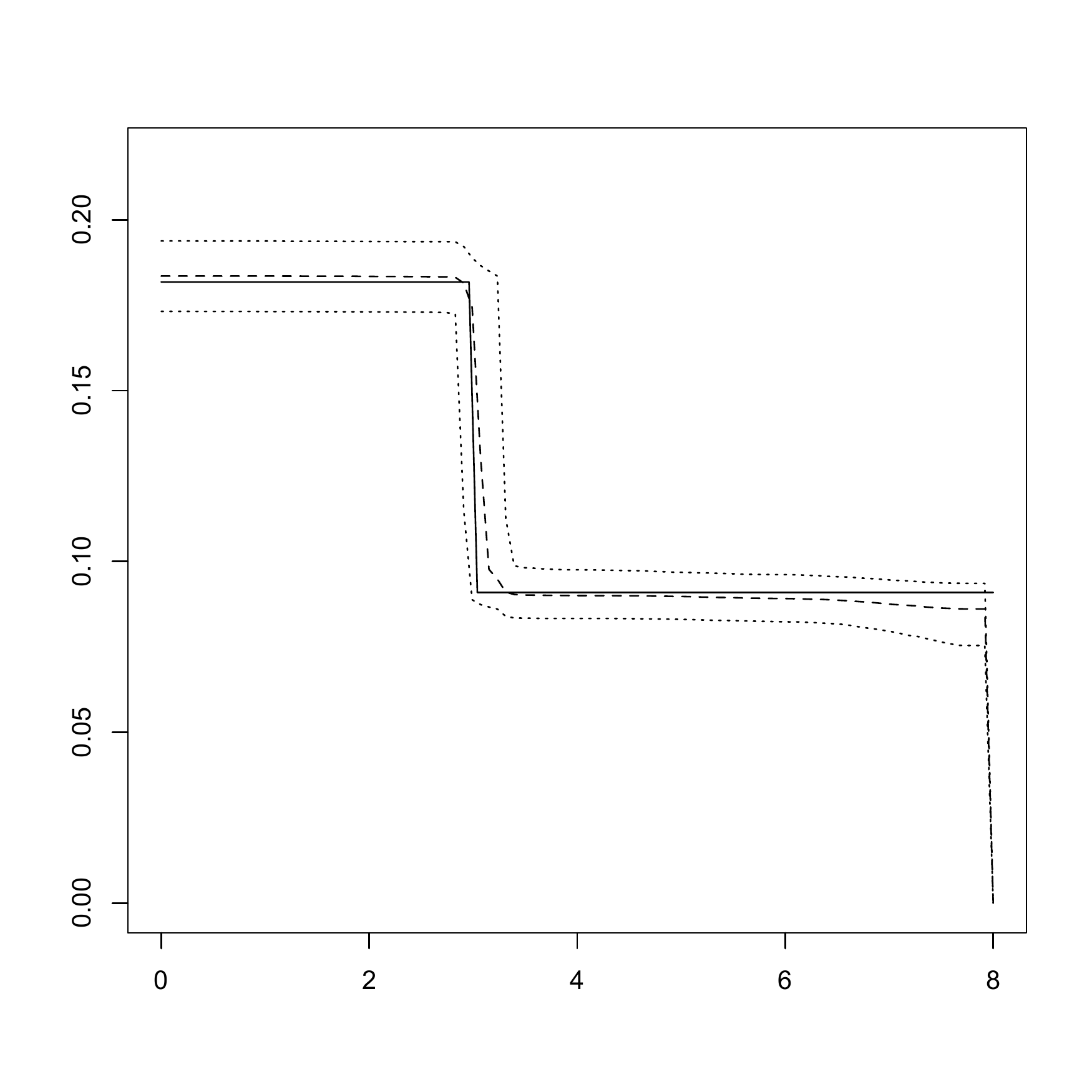} &
\includegraphics[width=0.33\textwidth]{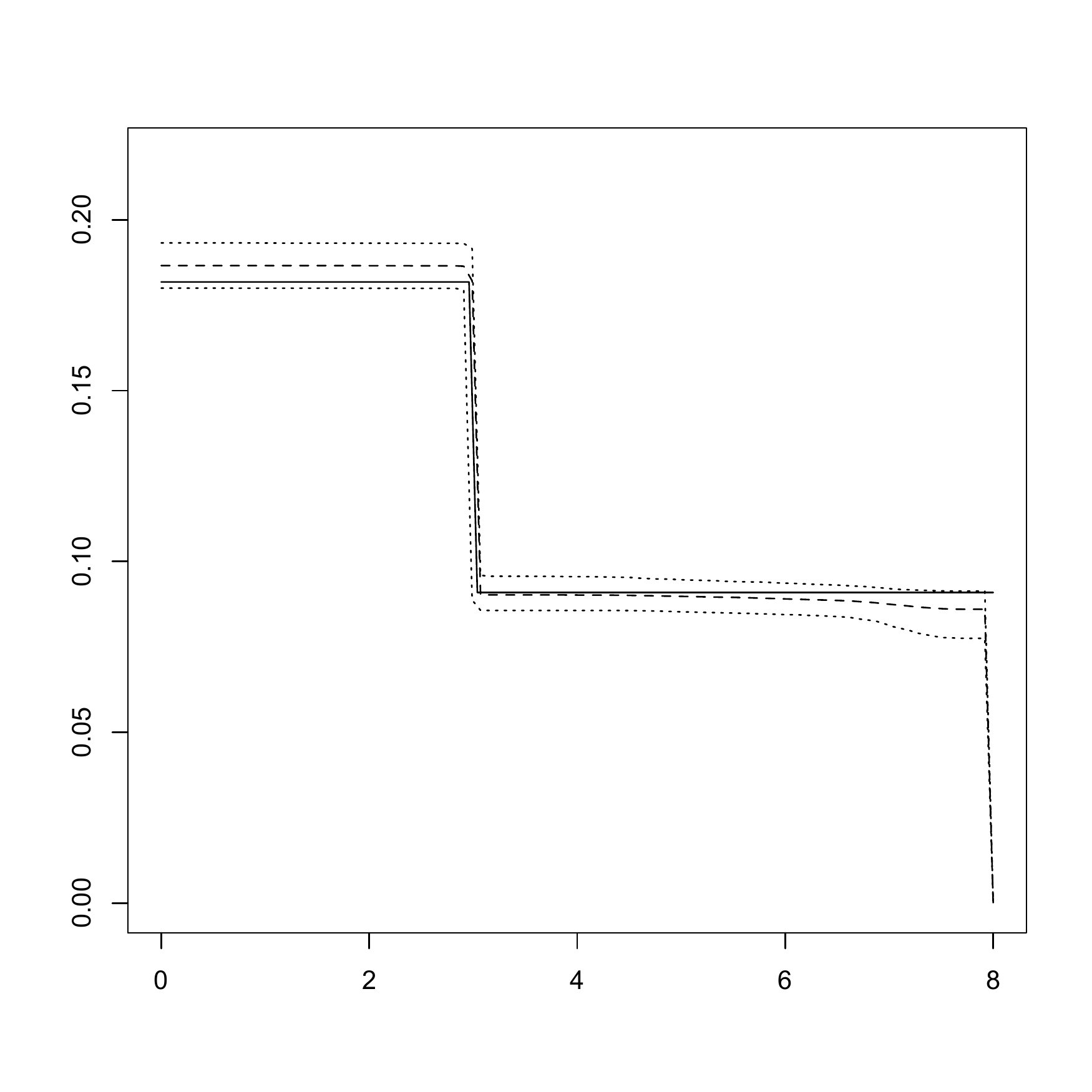} &
\includegraphics[width=0.33\textwidth]{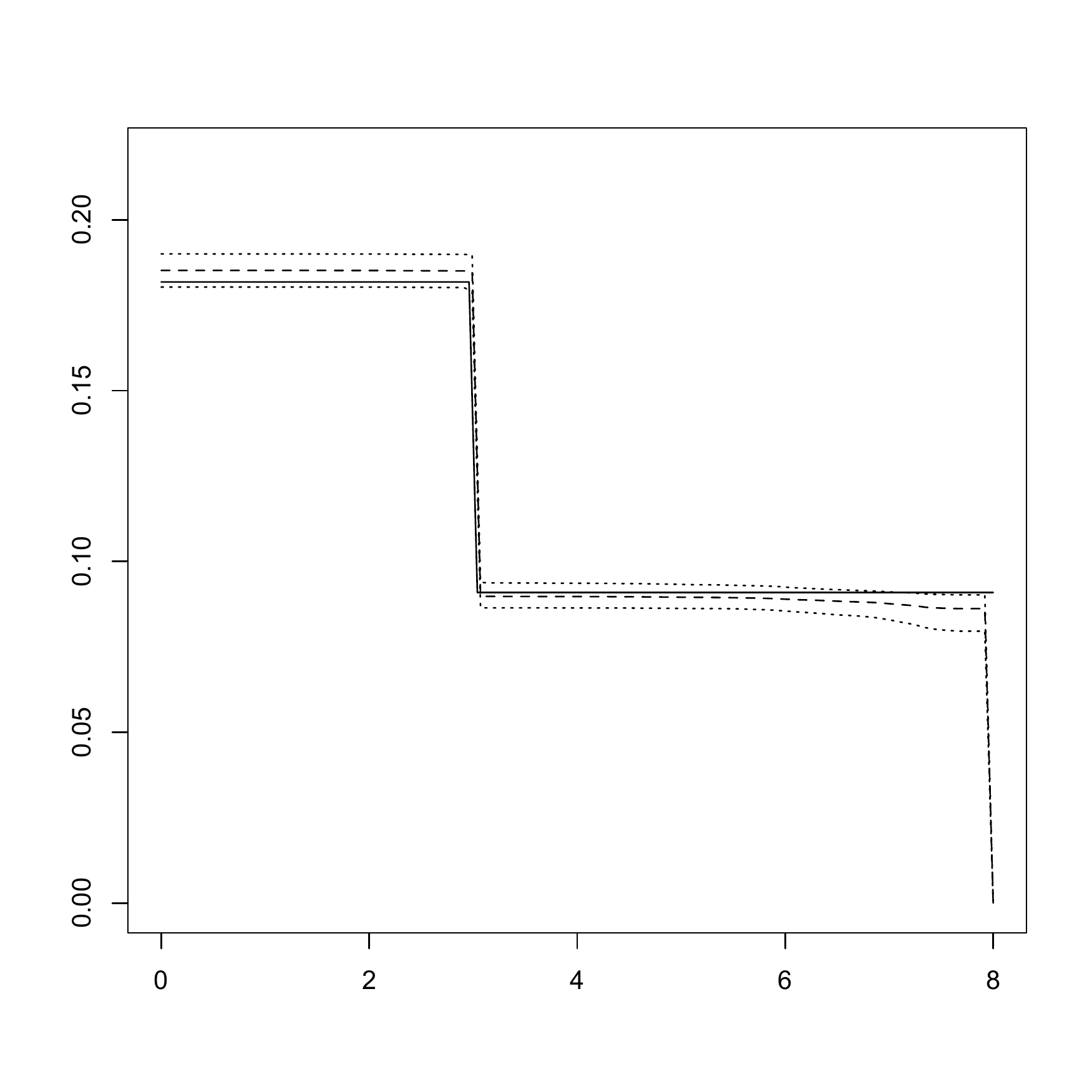}  
\end{tabular}
\caption{Estimation of $\overline{\lambda}_1$ from $D_{500}^1$ (first column), $D_{1000}^1$ (second column) and $D_{2000}^1$ (third column) using the four strategies: empirical prior (line 1), fixed $\gamma$ (line 2), hierarchical empirical prior (line 3), concentrated hierarchical empirical prior (line 4) . True density (plain line), estimation (dashed line) and confidence band (dotted lines)     }
\label{fig:estim1}
\end{figure}

\begin{figure}
\centering
\begin{tabular}{ccc}
$D^2_{500}$ & $D^2_{1000}$& $D^2_{2000}$ \\
\includegraphics[width=0.33\textwidth]{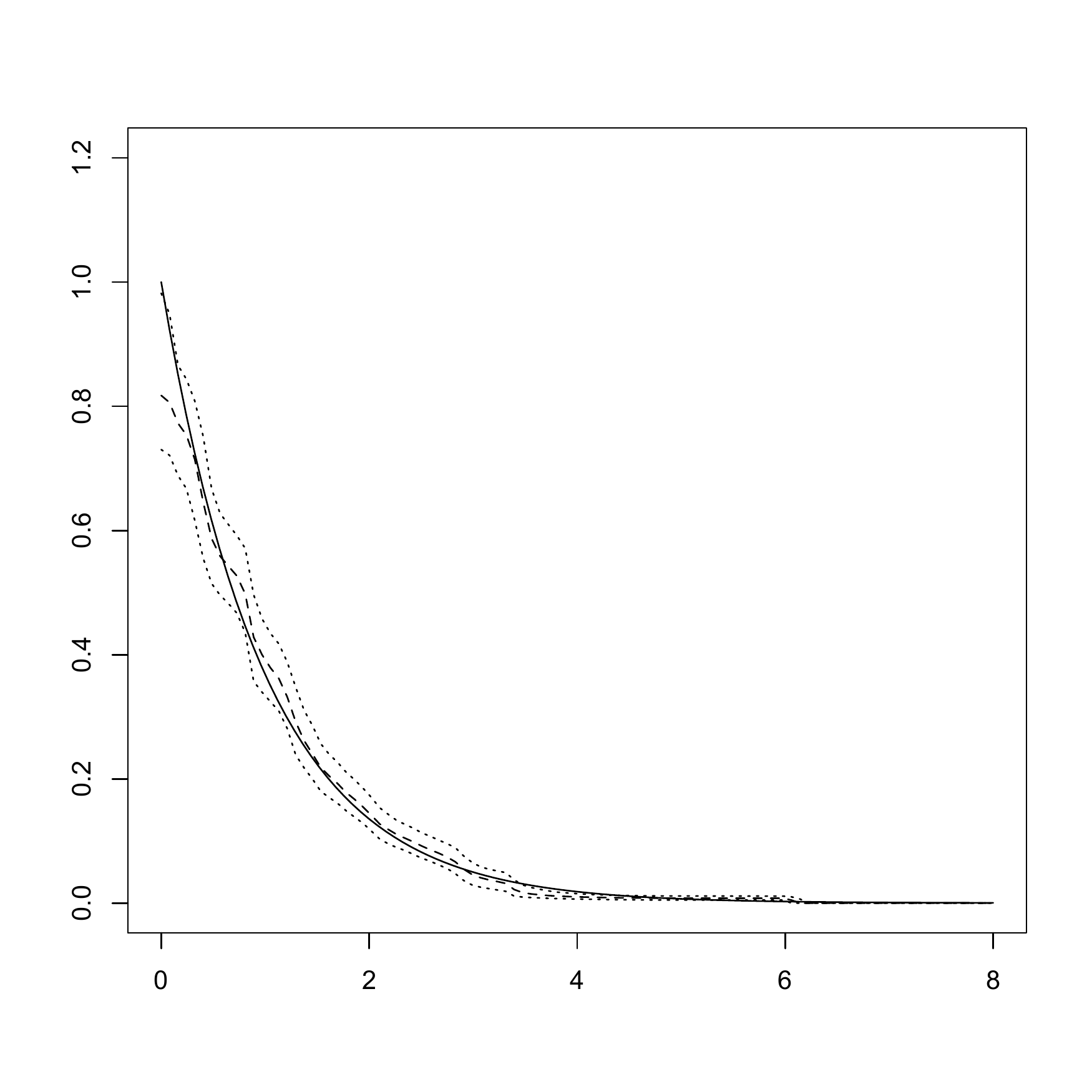} &
\includegraphics[width=0.33\textwidth]{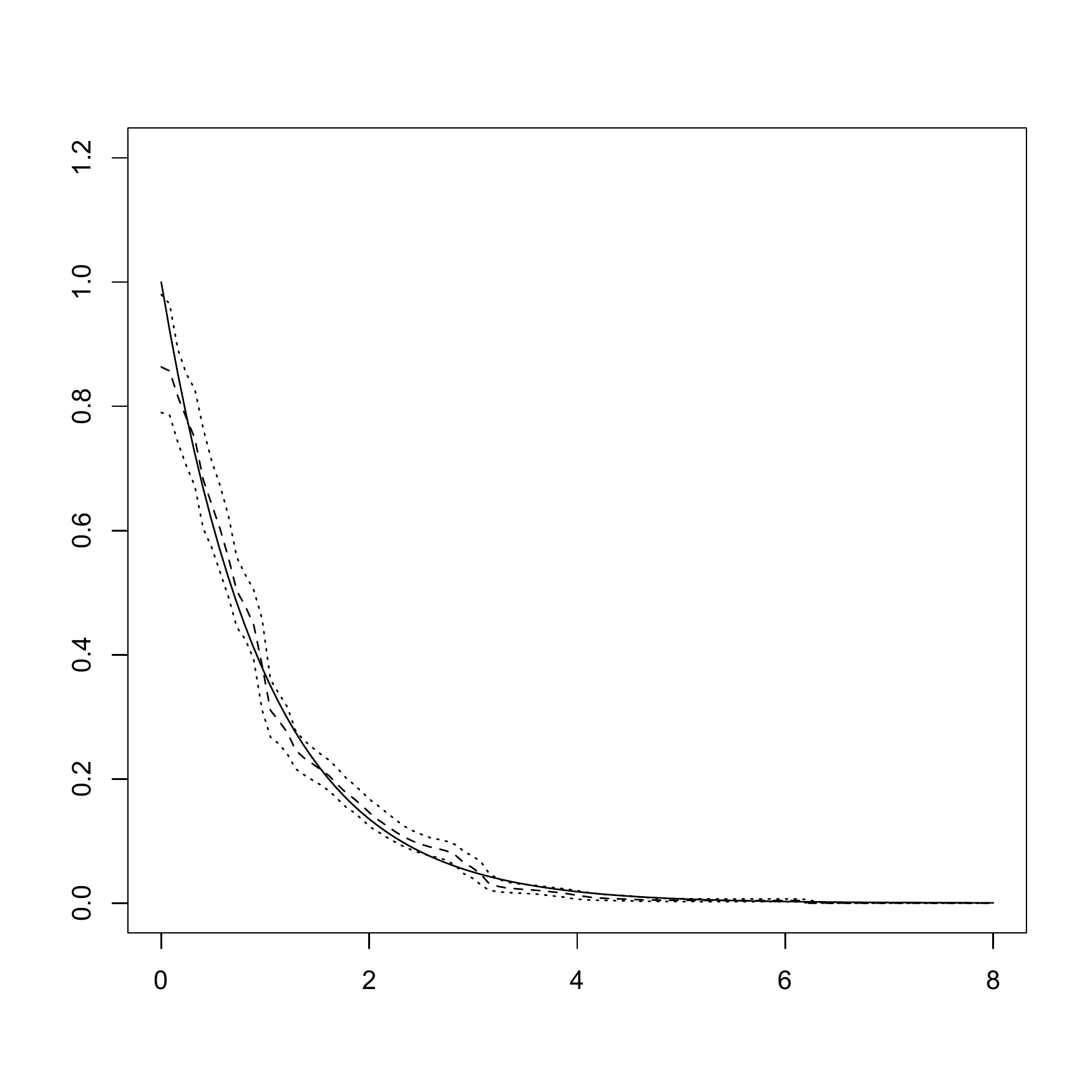} &
\includegraphics[width=0.33\textwidth]{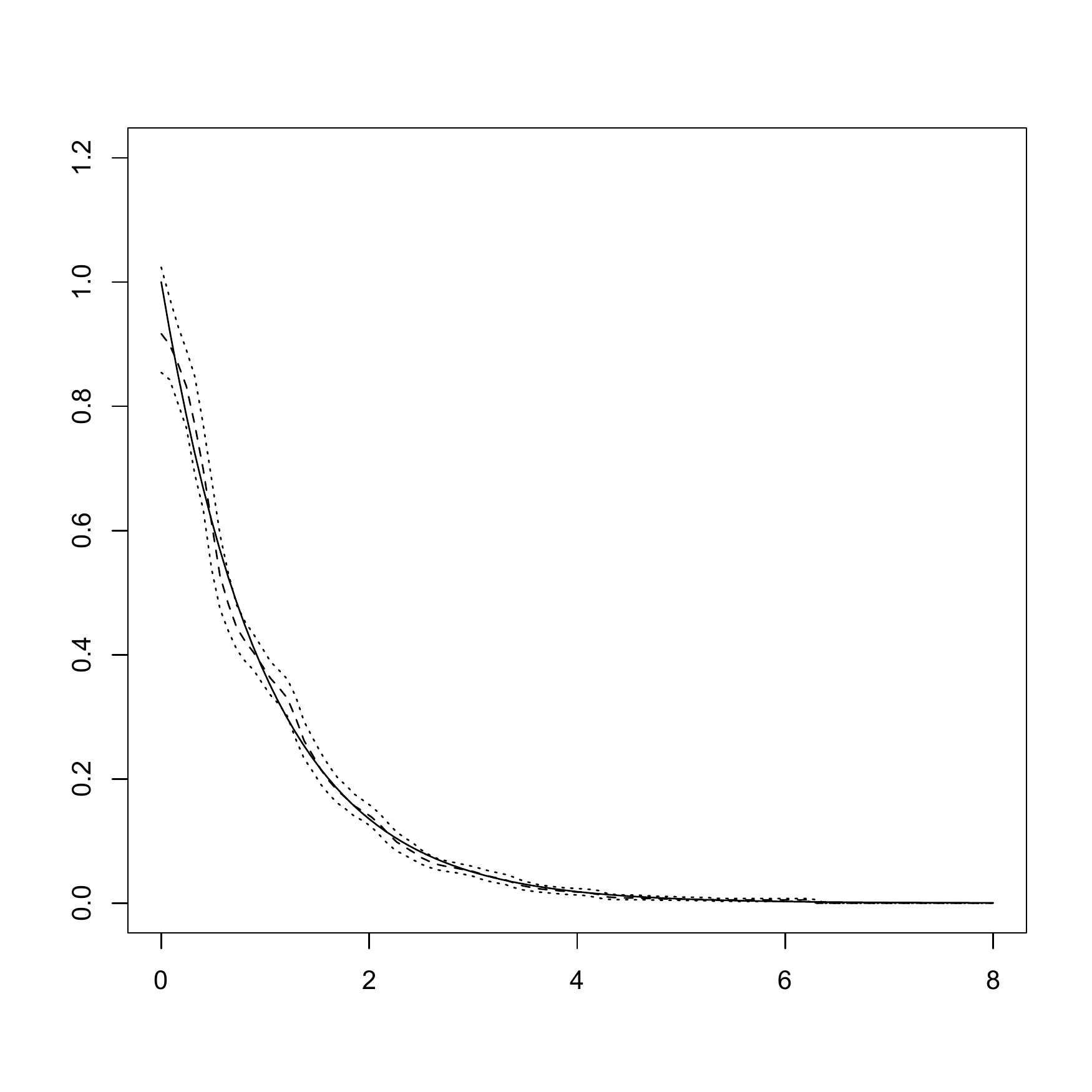} \\\includegraphics[width=0.33\textwidth]{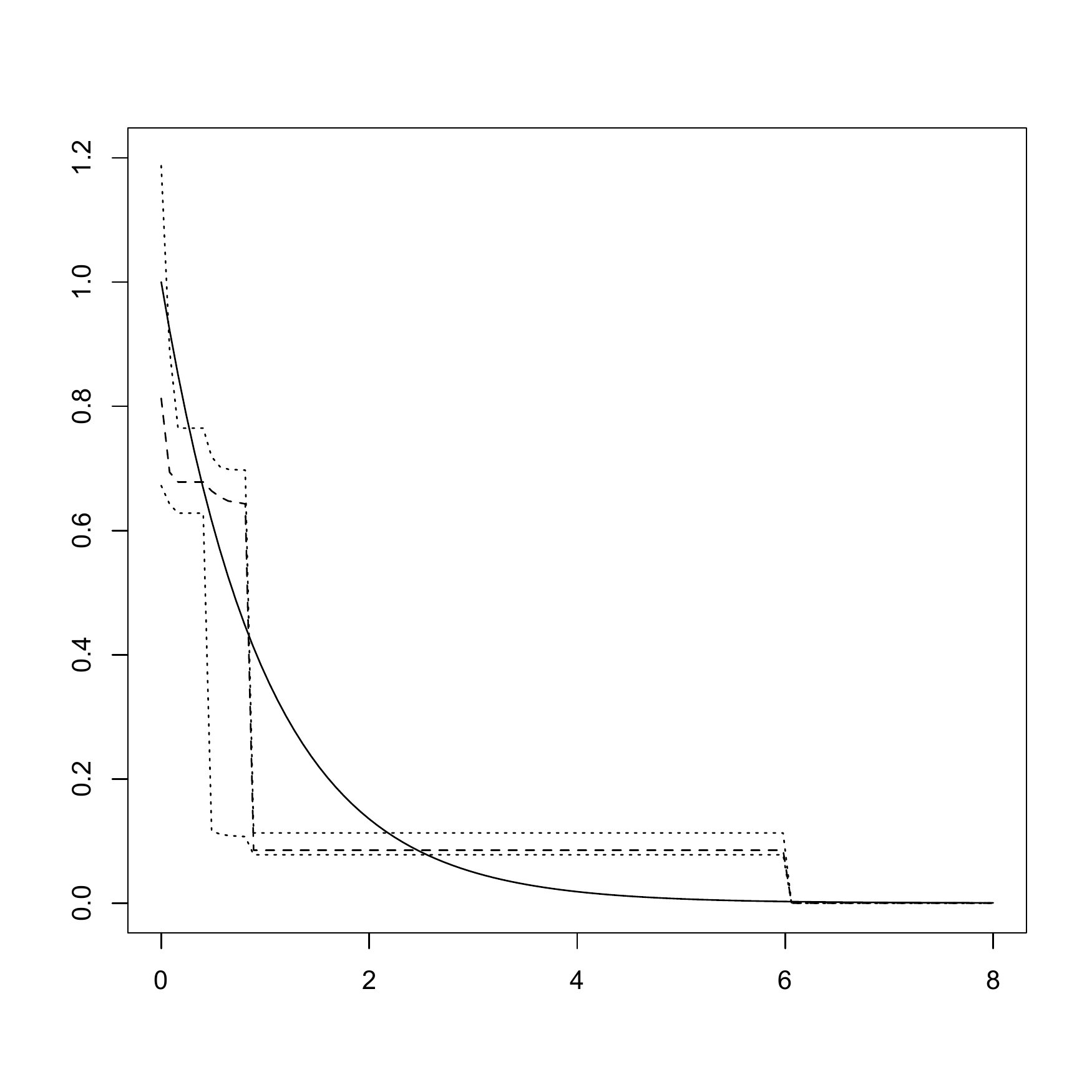} &
\includegraphics[width=0.33\textwidth]{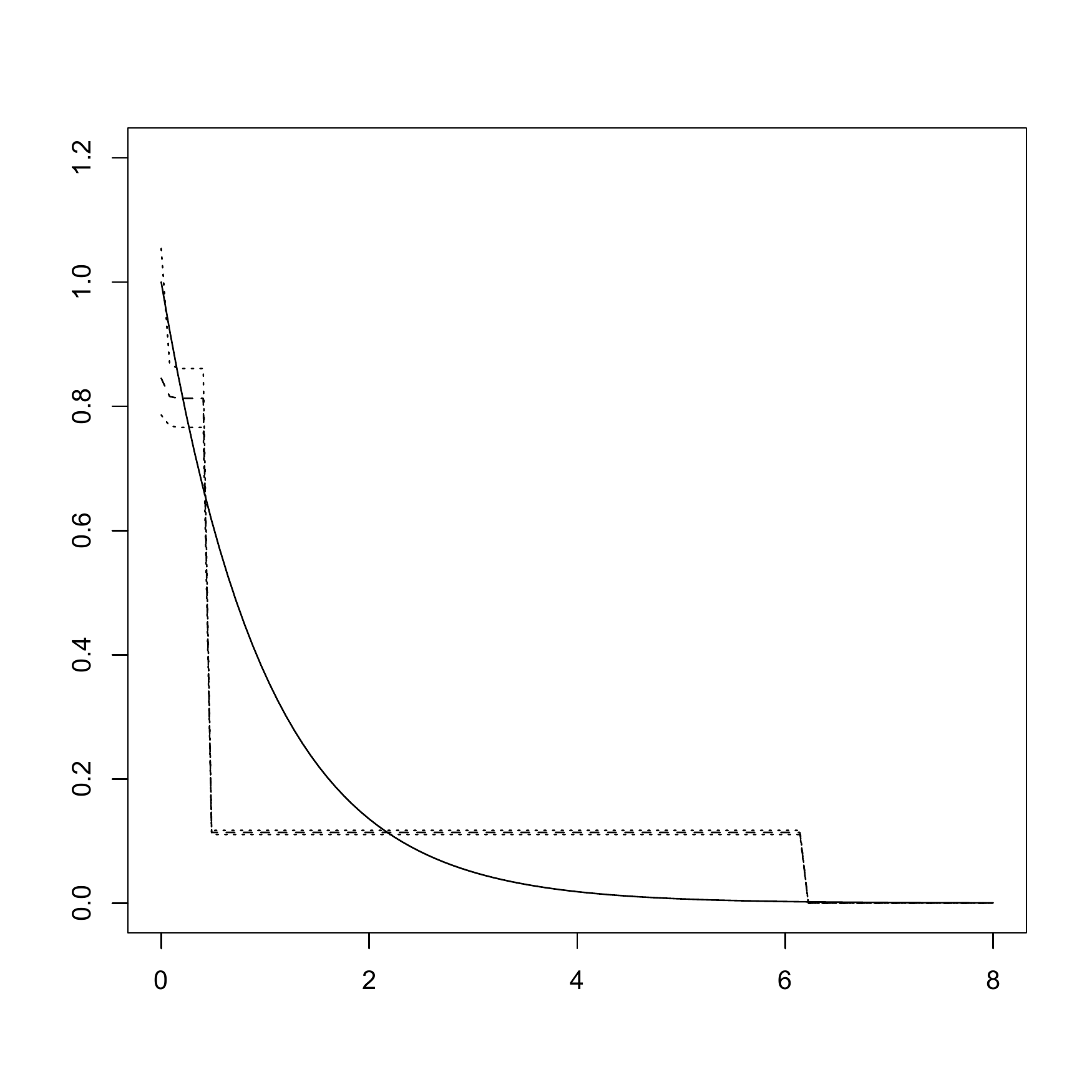} &
\includegraphics[width=0.33\textwidth]{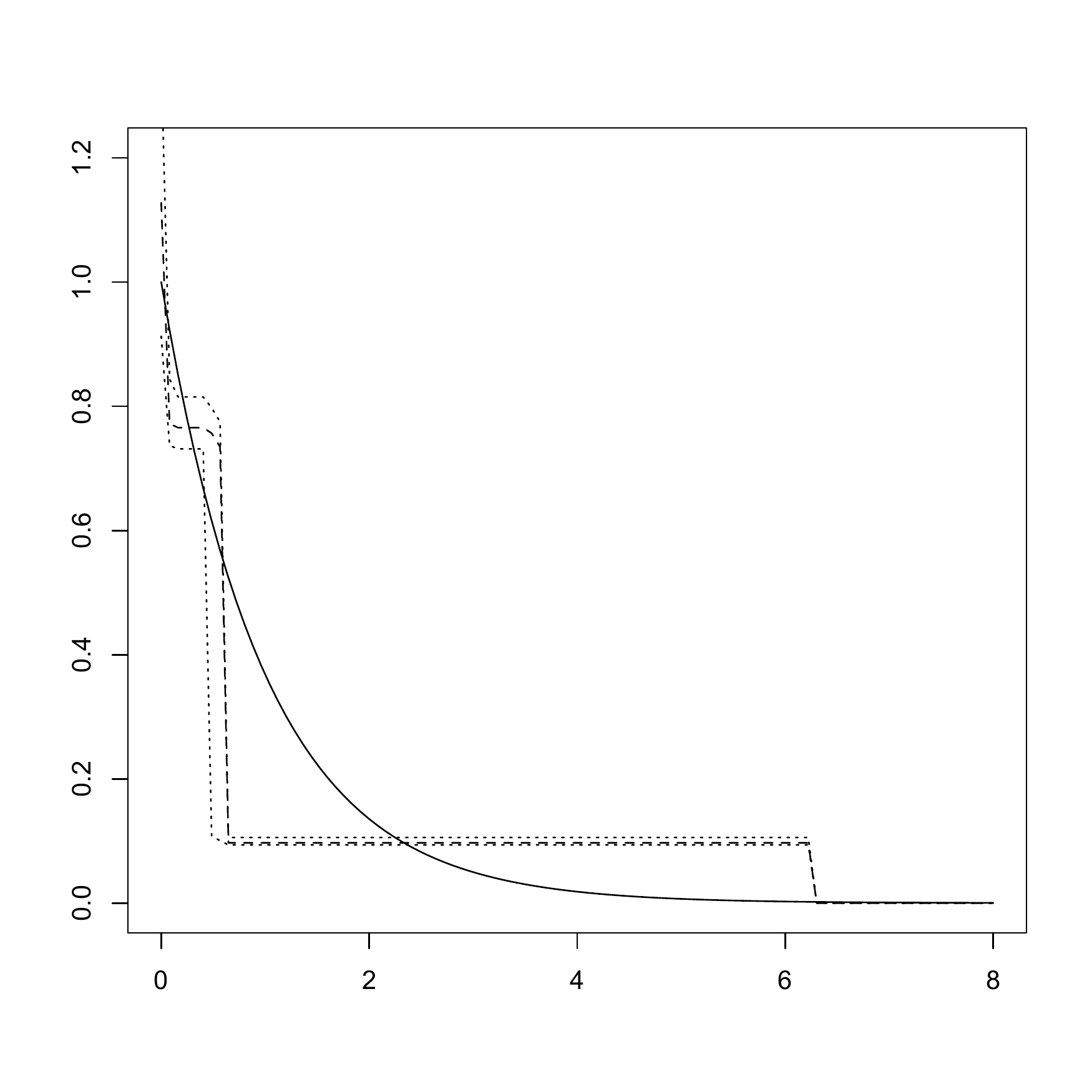} \\
\includegraphics[width=0.33\textwidth]{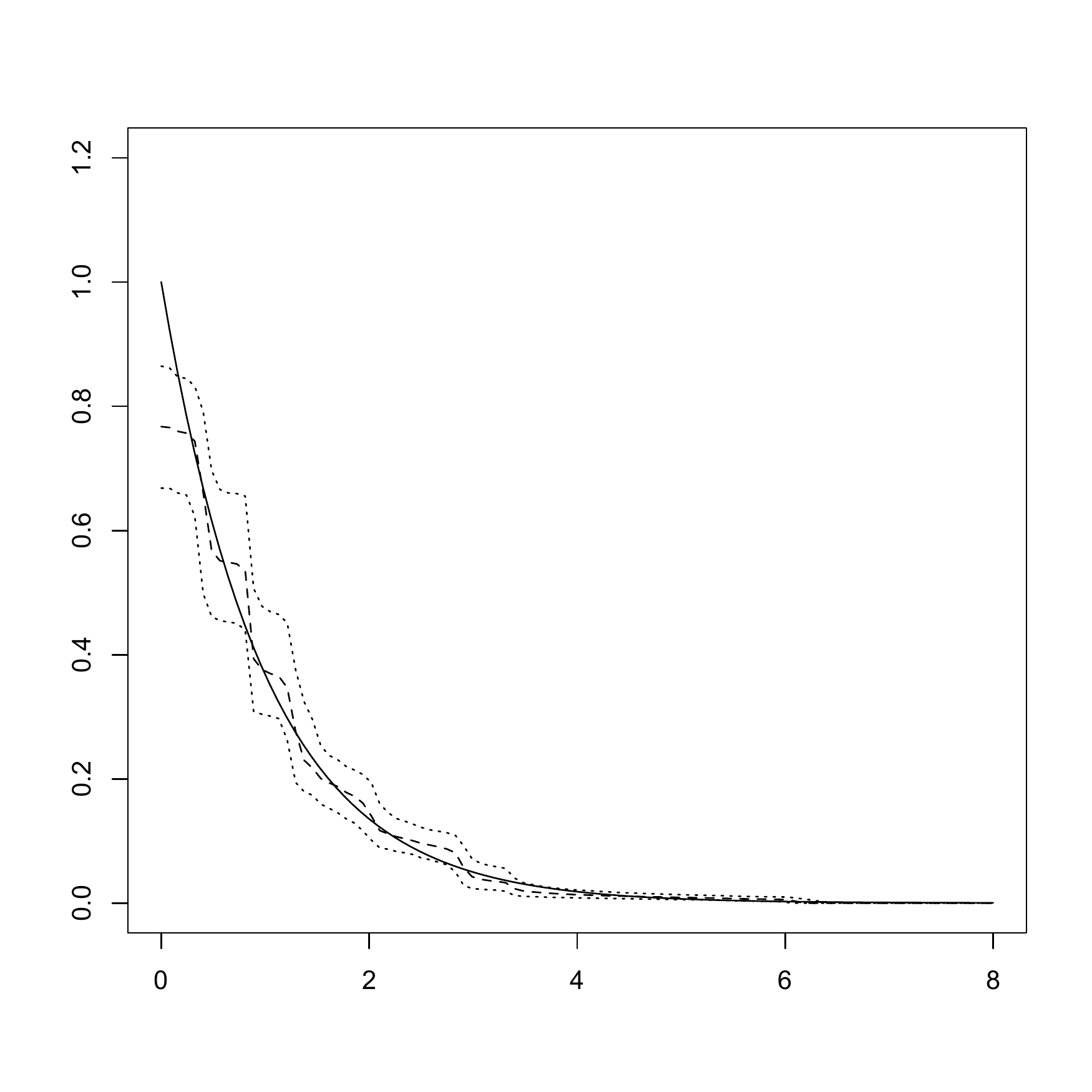} &
\includegraphics[width=0.33\textwidth]{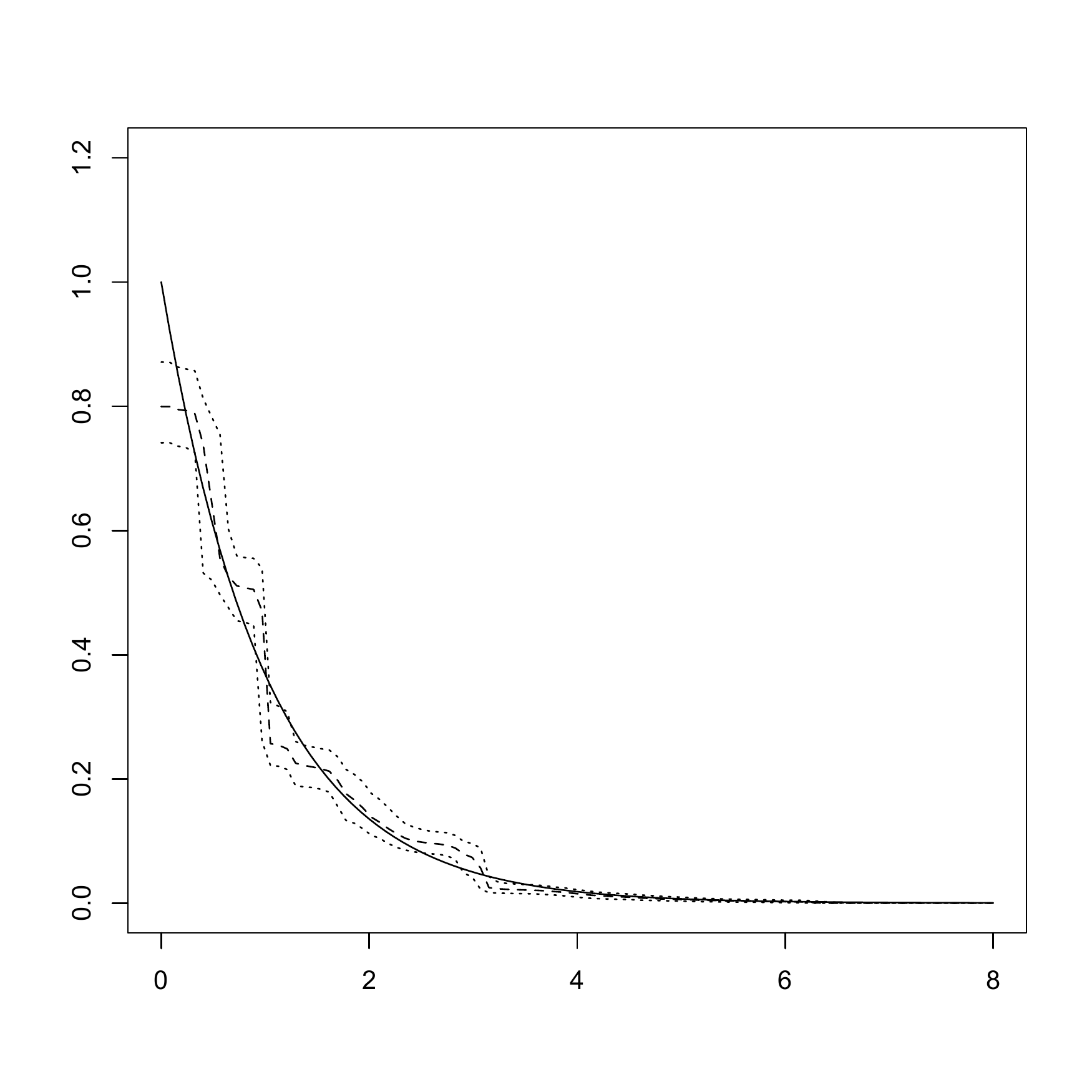} &
\includegraphics[width=0.33\textwidth]{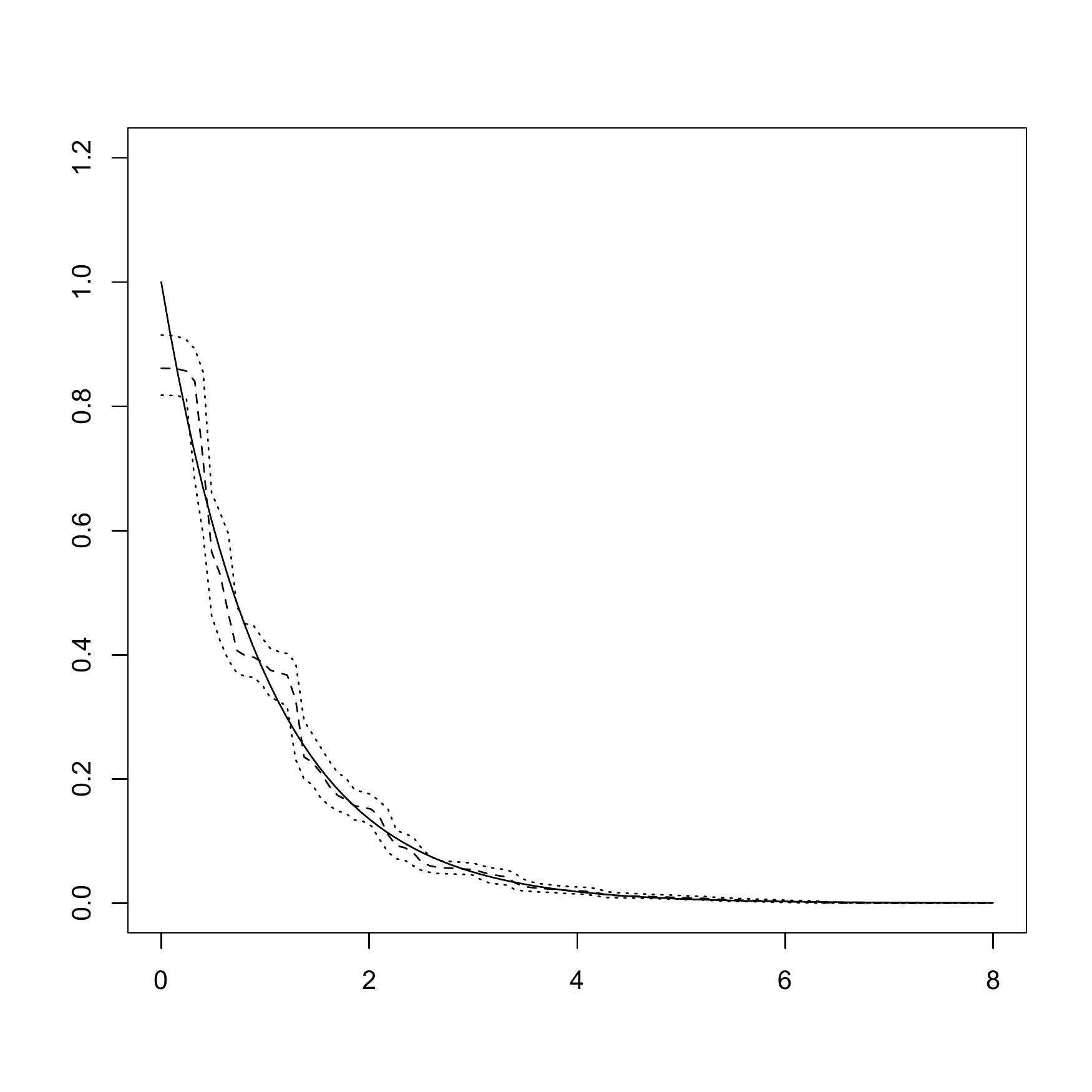} \\
\includegraphics[width=0.33\textwidth]{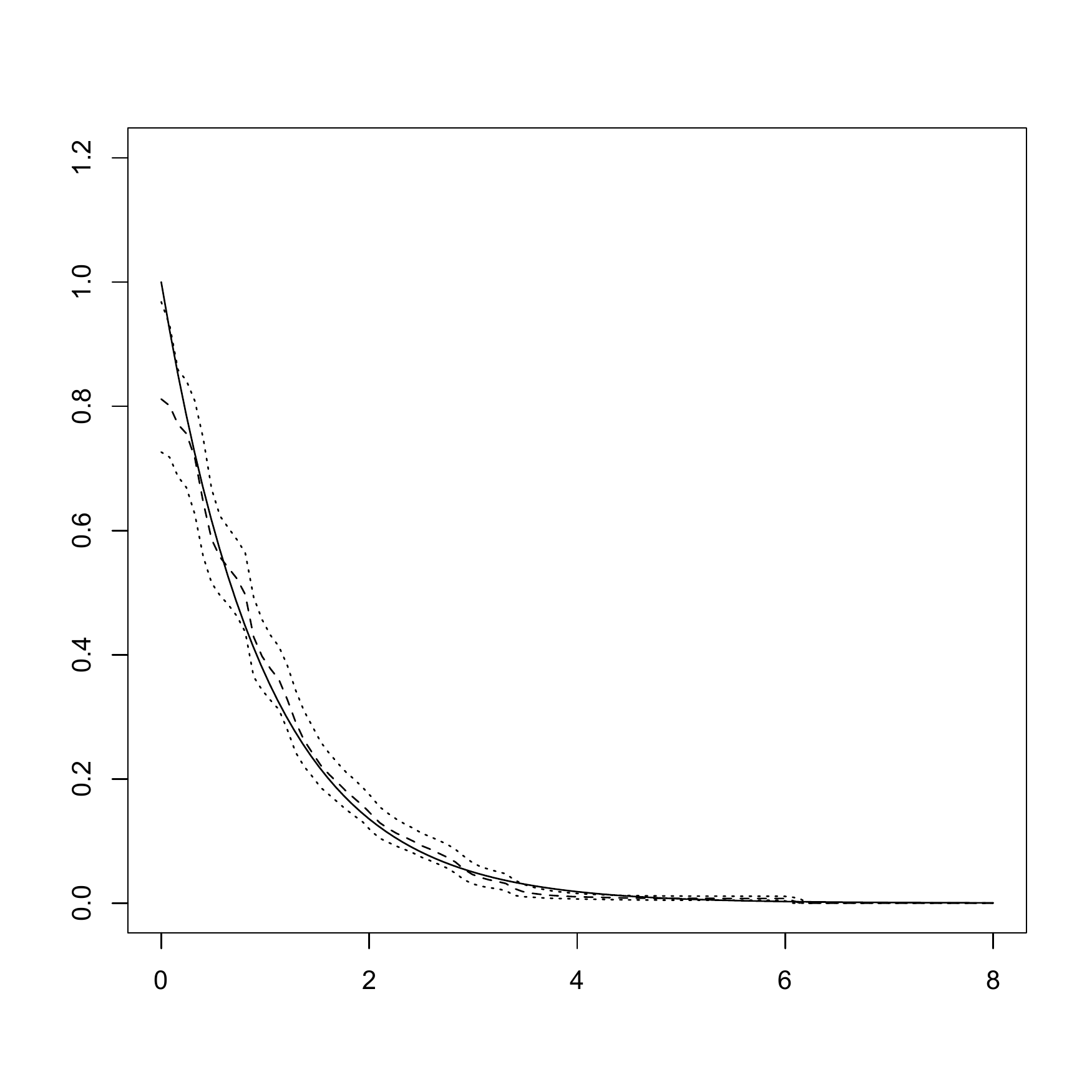} &
\includegraphics[width=0.33\textwidth]{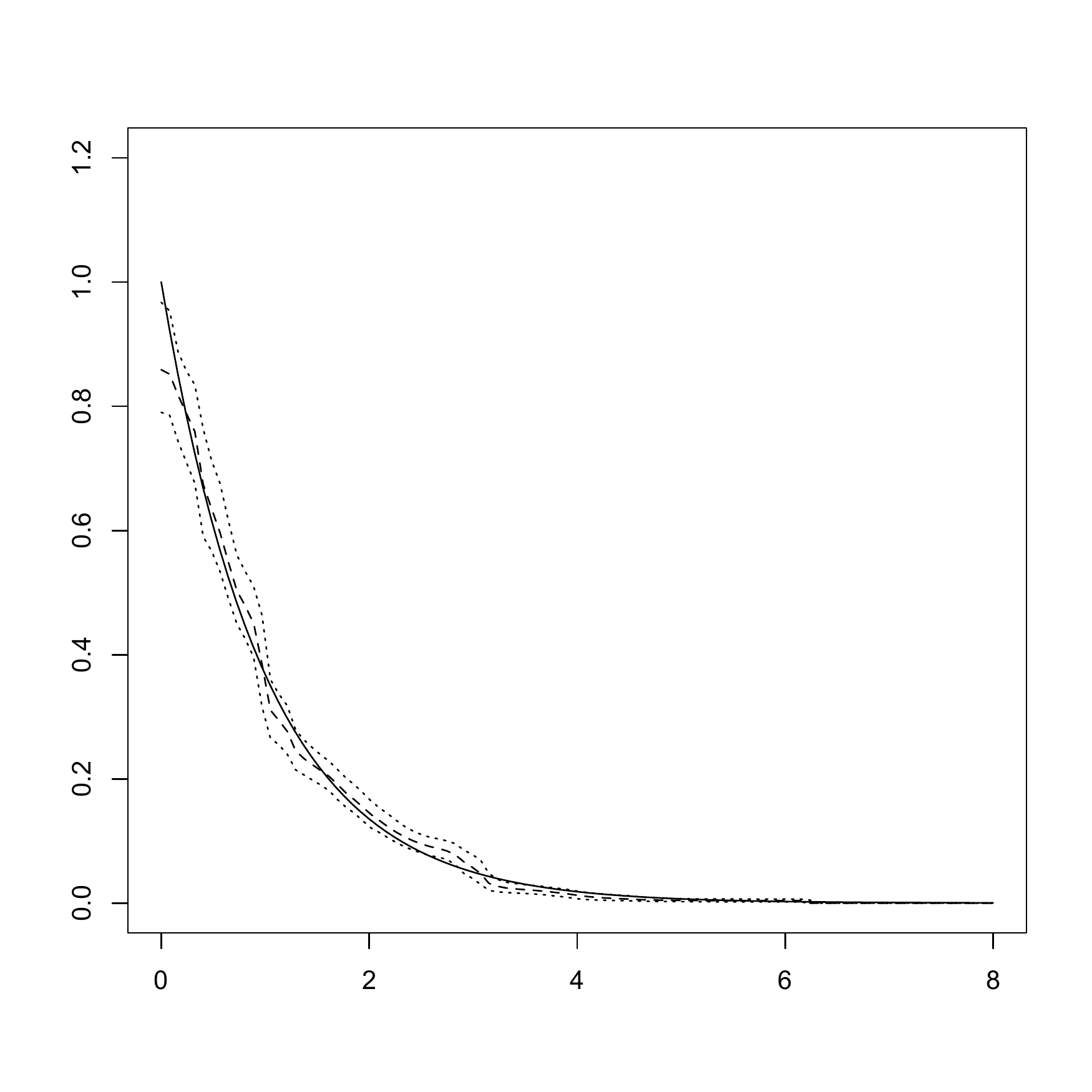} &
\includegraphics[width=0.33\textwidth]{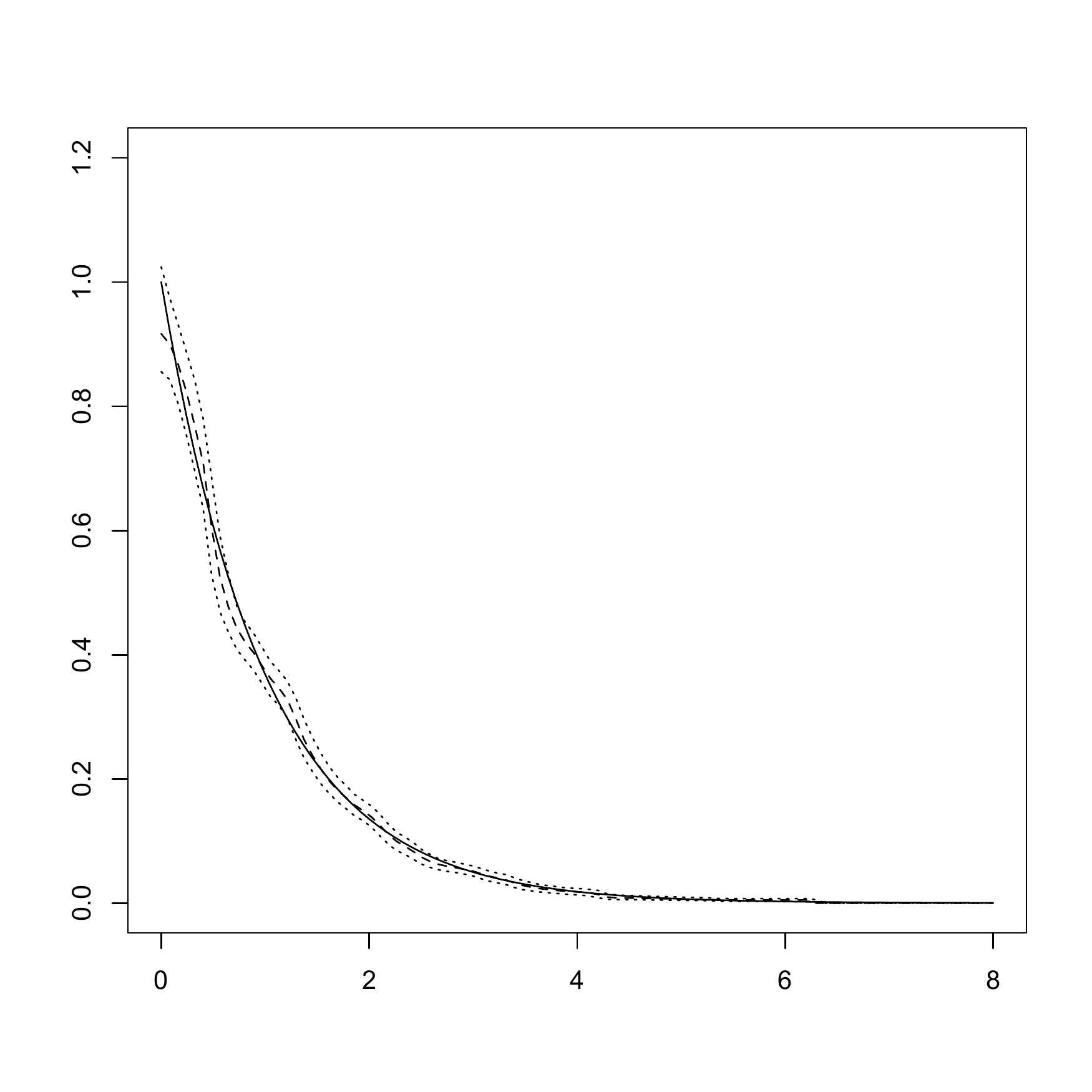}  
\end{tabular}
\caption{Estimation of $\overline{\lambda}_2$ from $D_{500}^2$ (first column), $D_{1000}^2$ (second column) and $D_{2000}^2$ (third column) using the four strategies: empirical prior (line 1), fixe $\gamma$ (line 2), hierarchical empirical prior (line 3), concentrated hierarchical empirical prior (line 4) . True density (plain line), estimation (dashed line) and confidence band (dotted lines)     }
\label{fig:estim2}
\end{figure}

\begin{figure}
\centering
\begin{tabular}{ccc}
$D^3_{500}$ & $D^3_{1000}$& $D^3_{2000}$ \\
\includegraphics[width=0.33\textwidth]{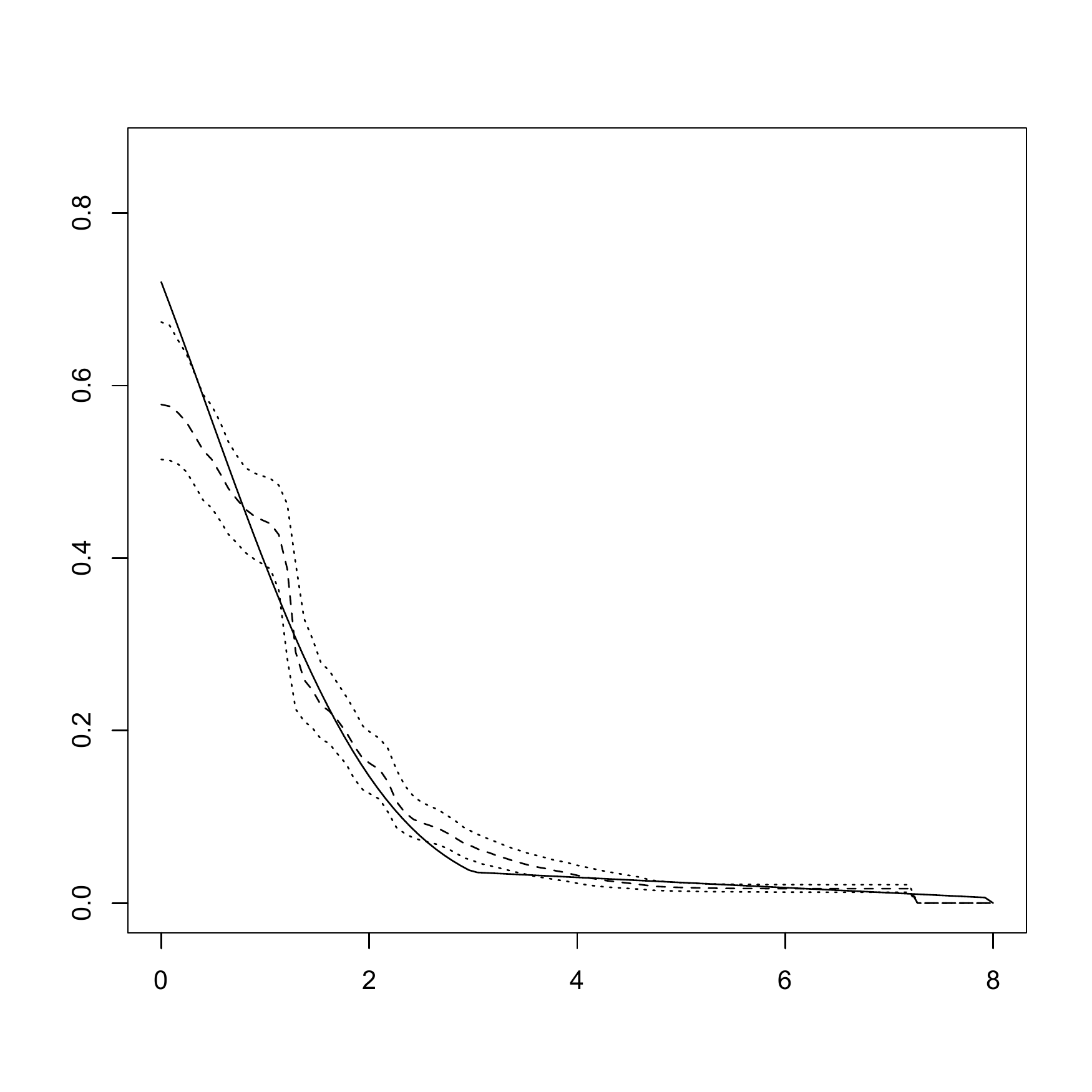} &
\includegraphics[width=0.33\textwidth]{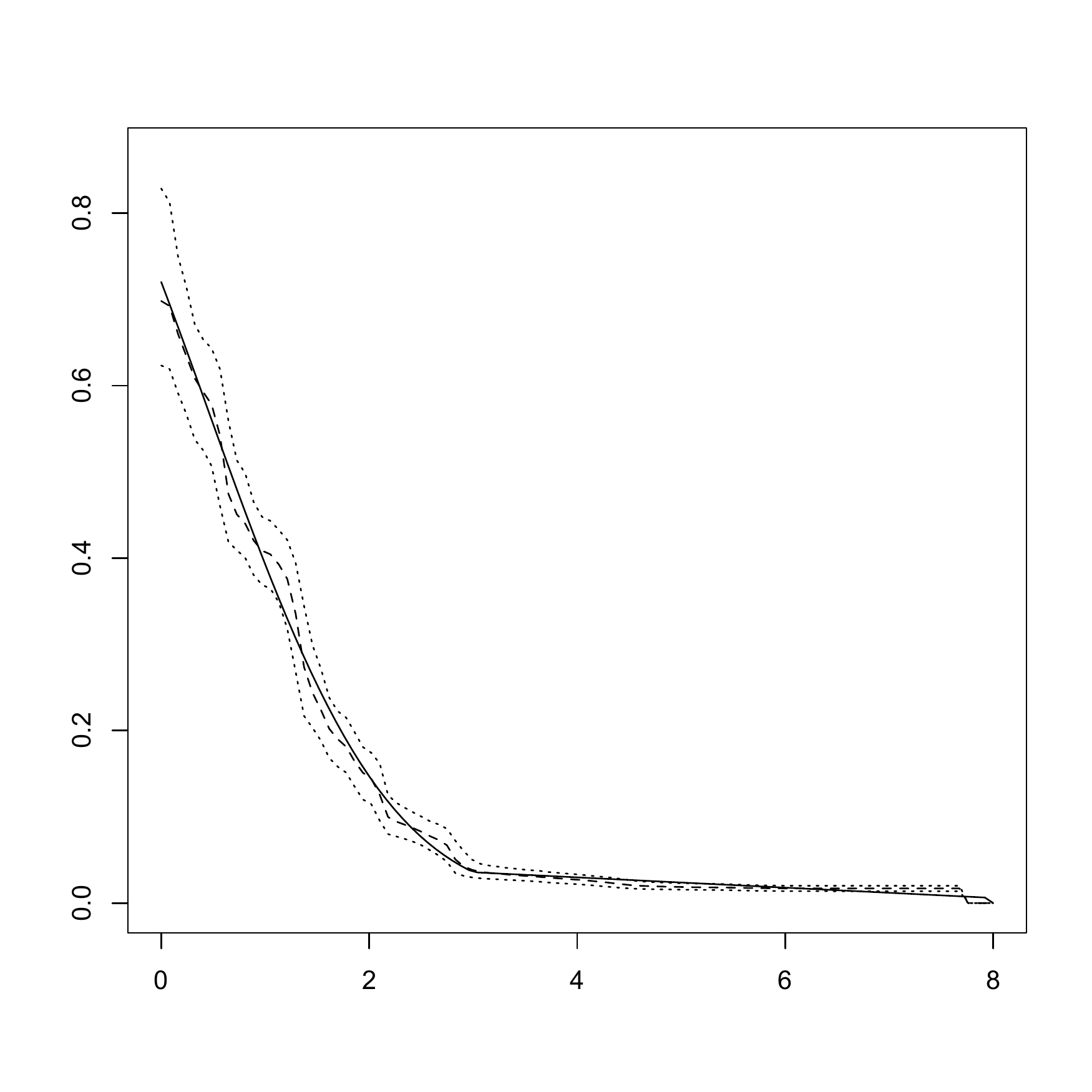} &
\includegraphics[width=0.33\textwidth]{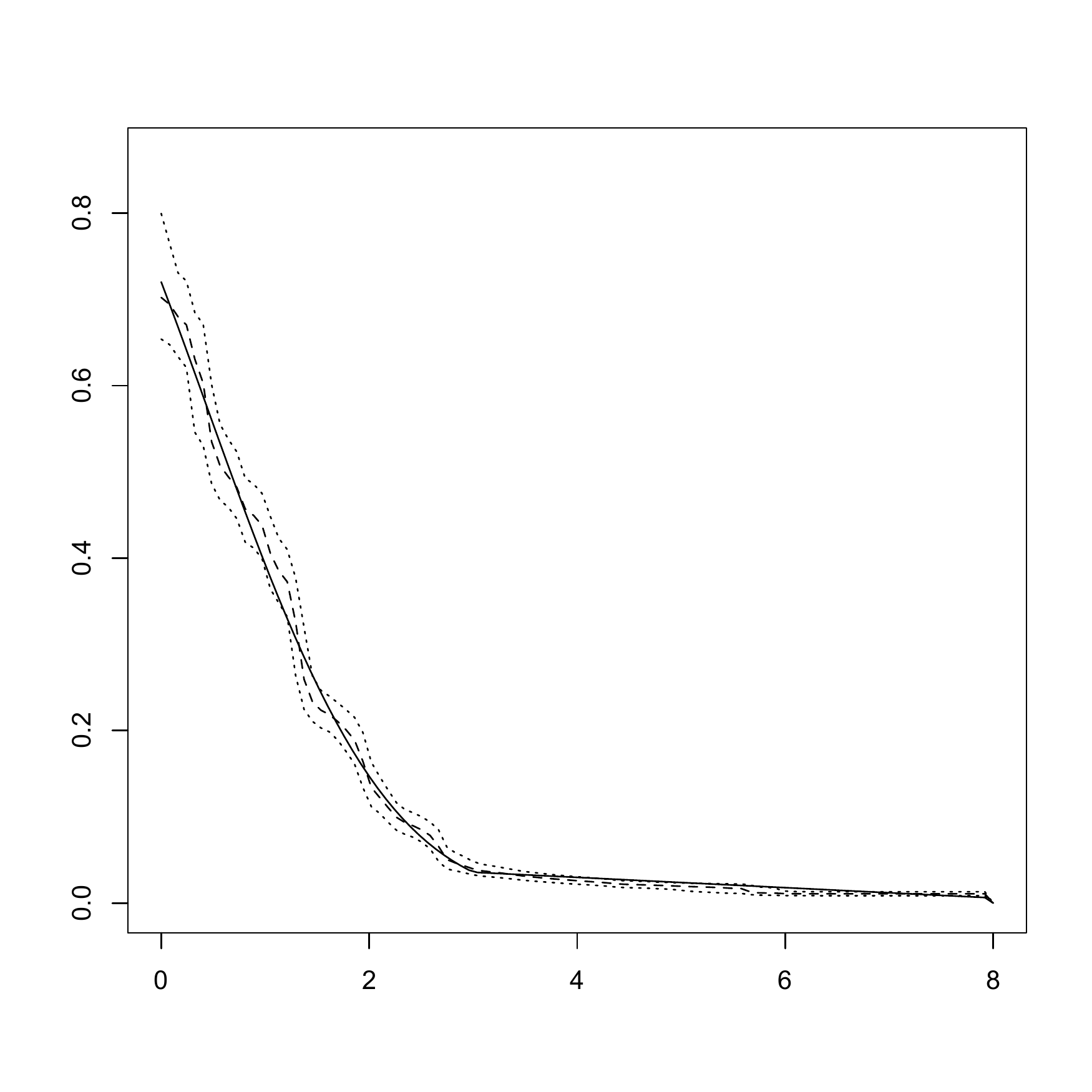} \\\includegraphics[width=0.33\textwidth]{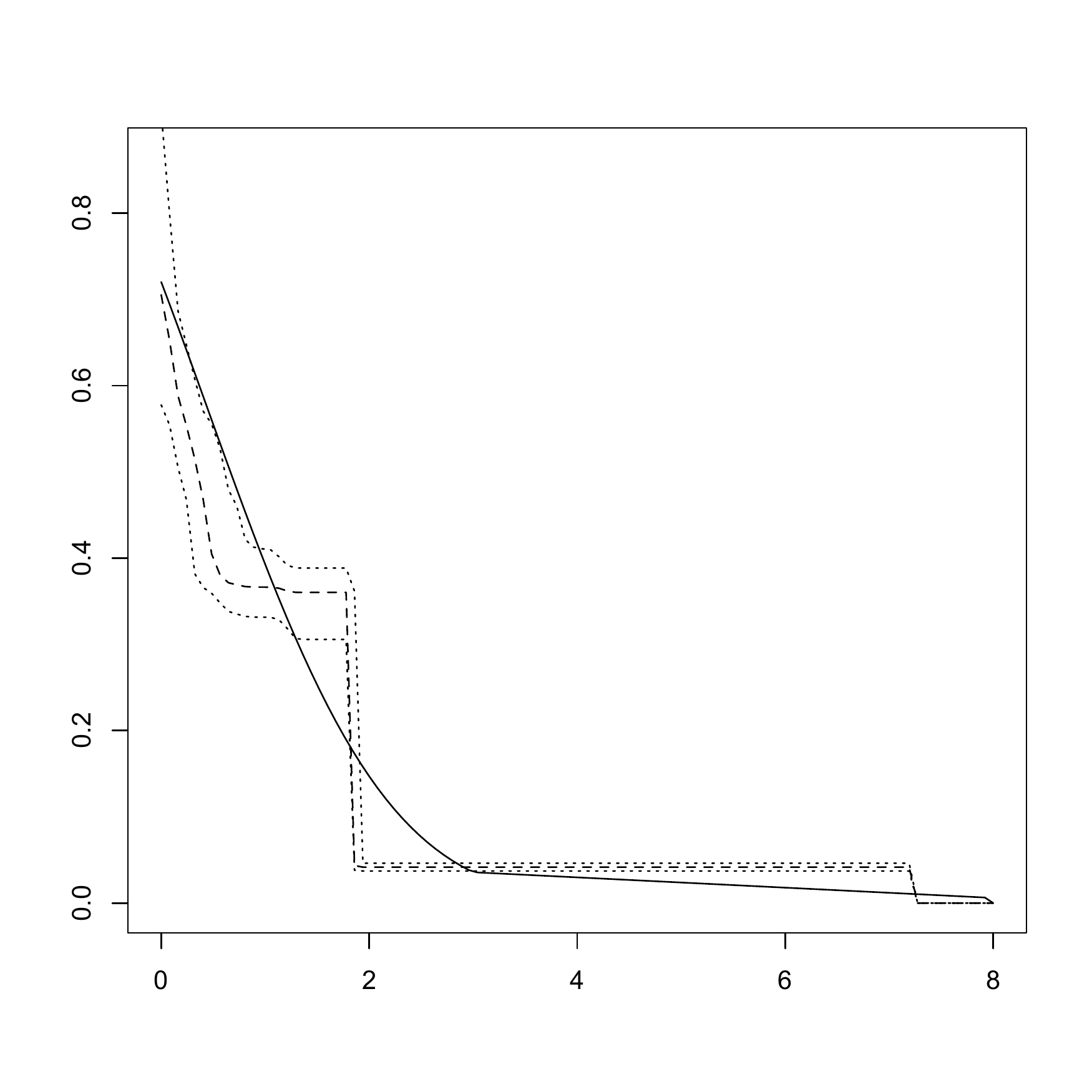} &
\includegraphics[width=0.33\textwidth]{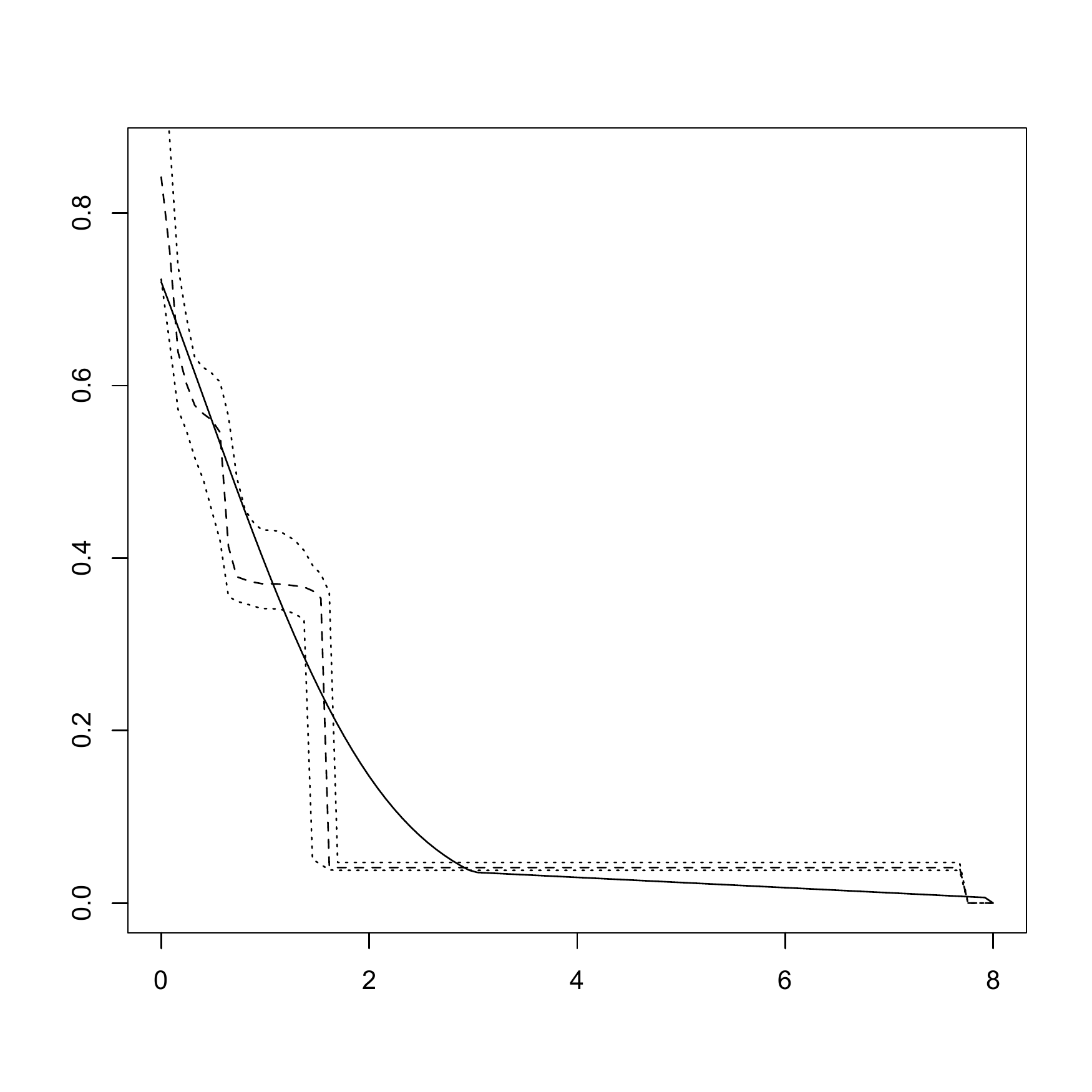} &
\includegraphics[width=0.33\textwidth]{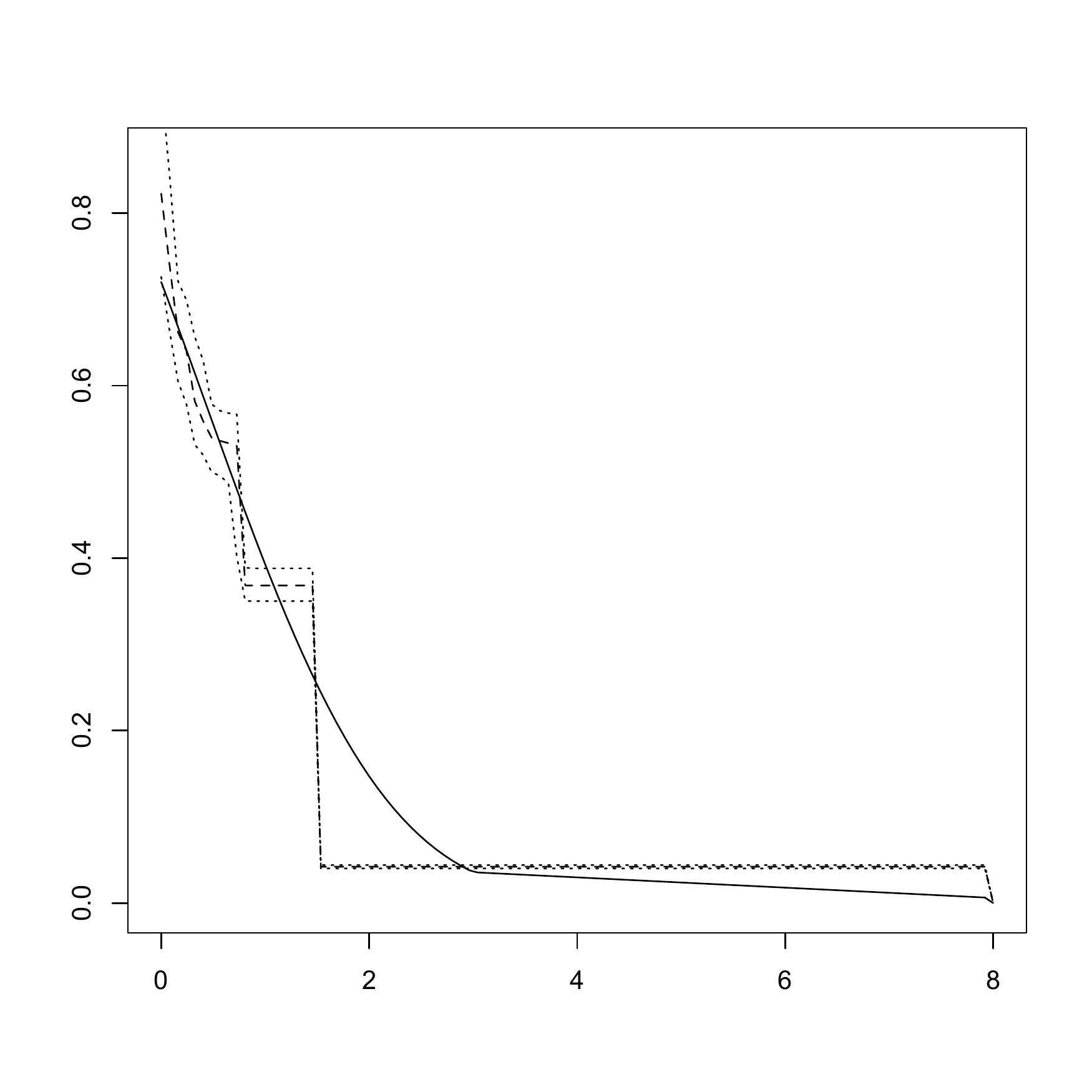} \\
\includegraphics[width=0.33\textwidth]{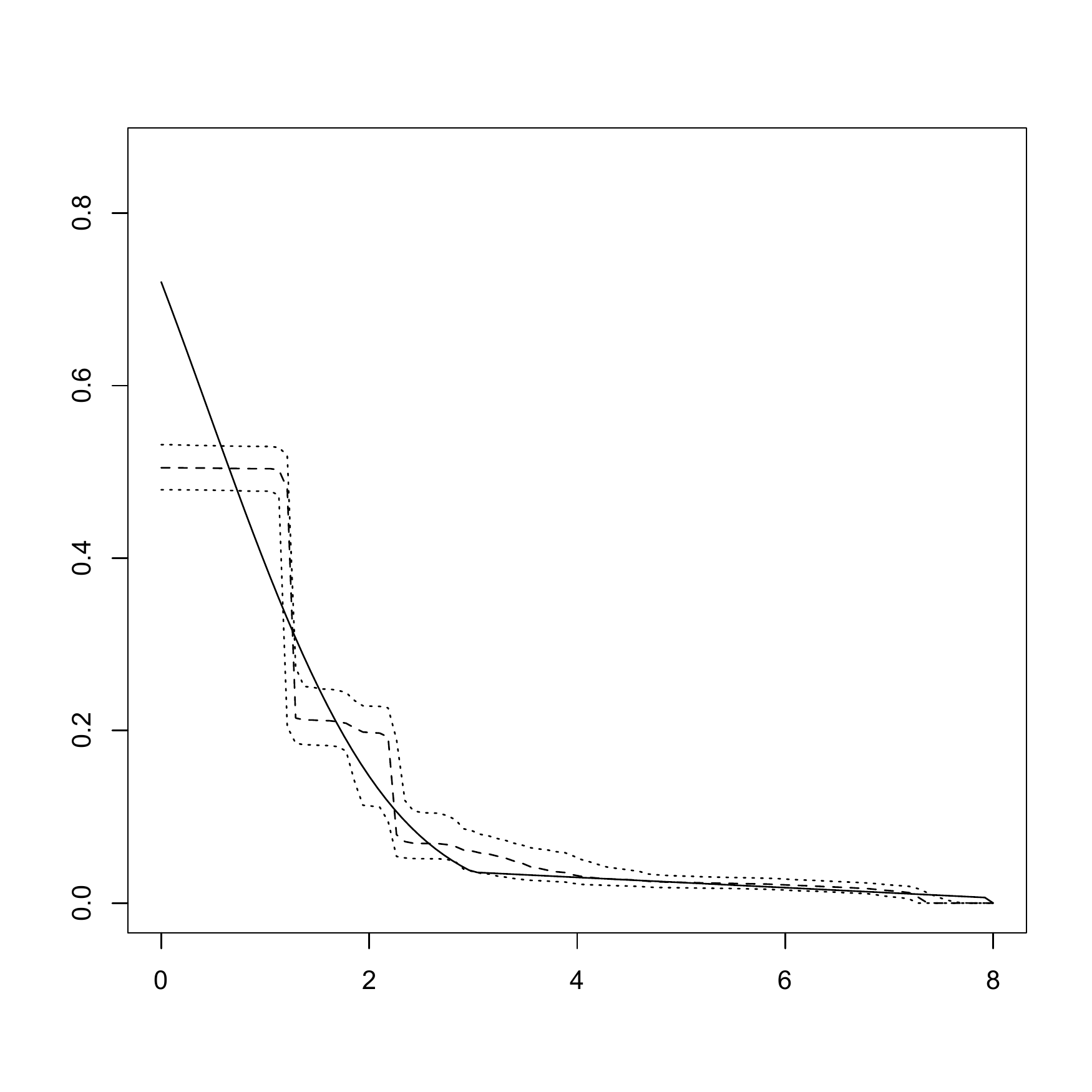} &
\includegraphics[width=0.33\textwidth]{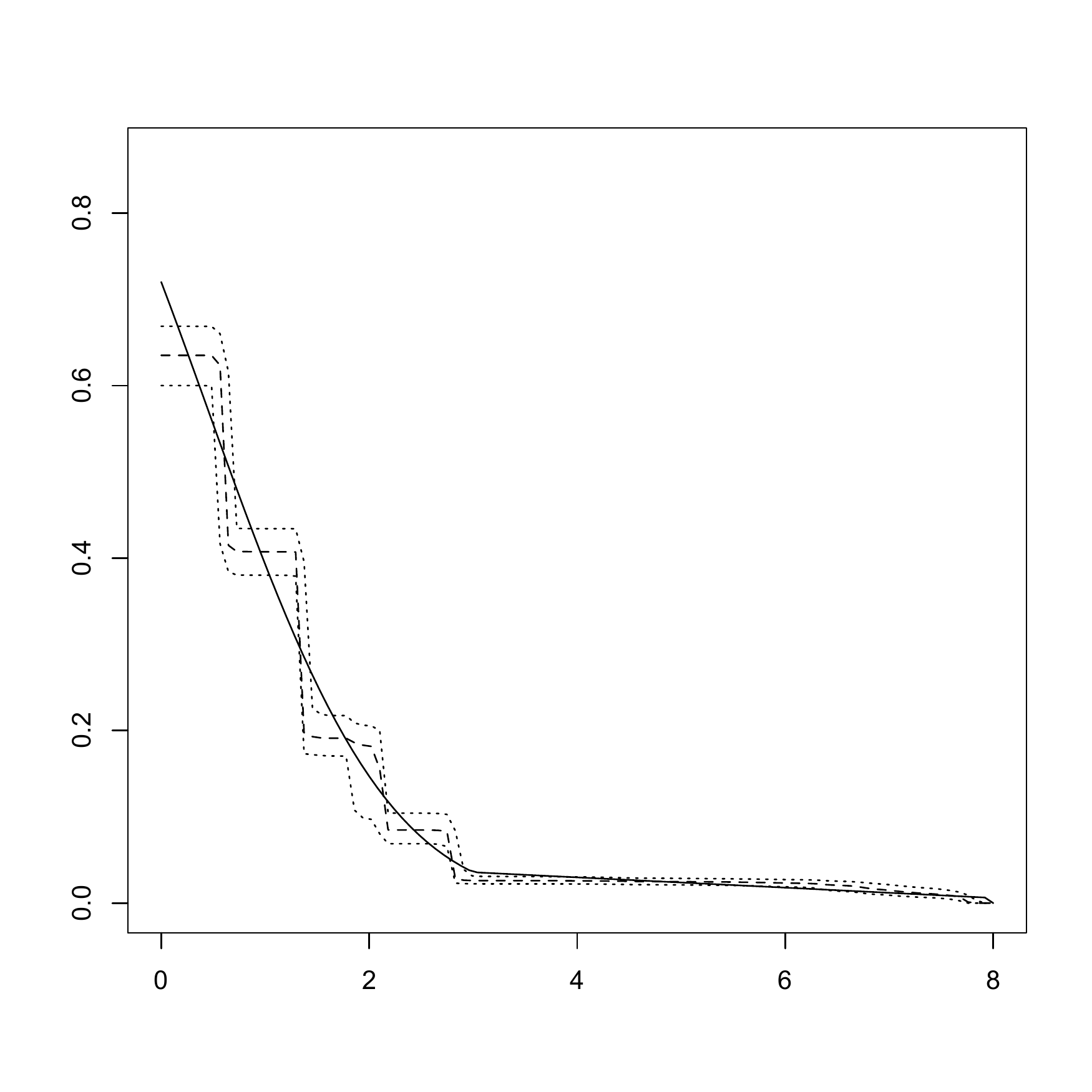} &
\includegraphics[width=0.33\textwidth]{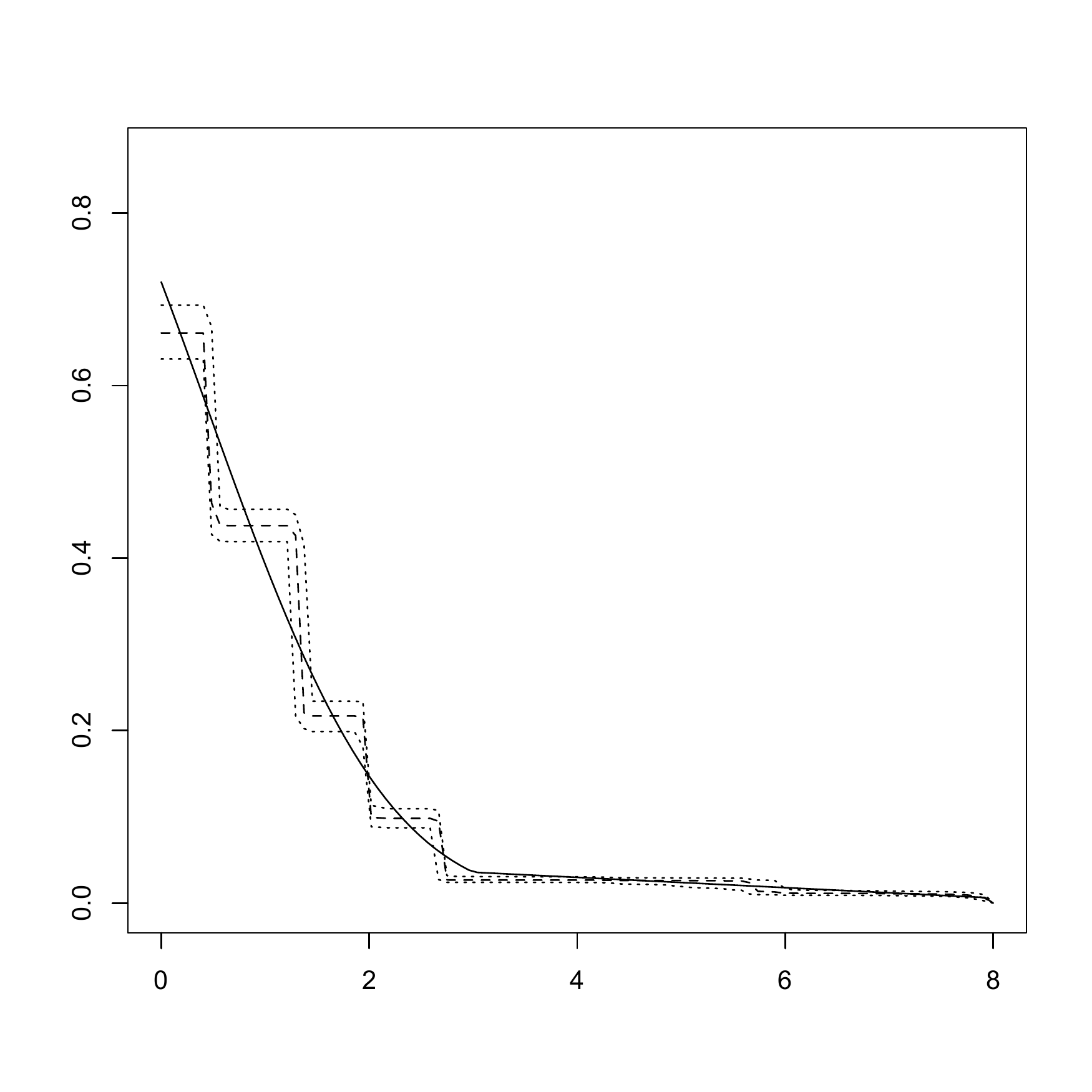} \\
\includegraphics[width=0.33\textwidth]{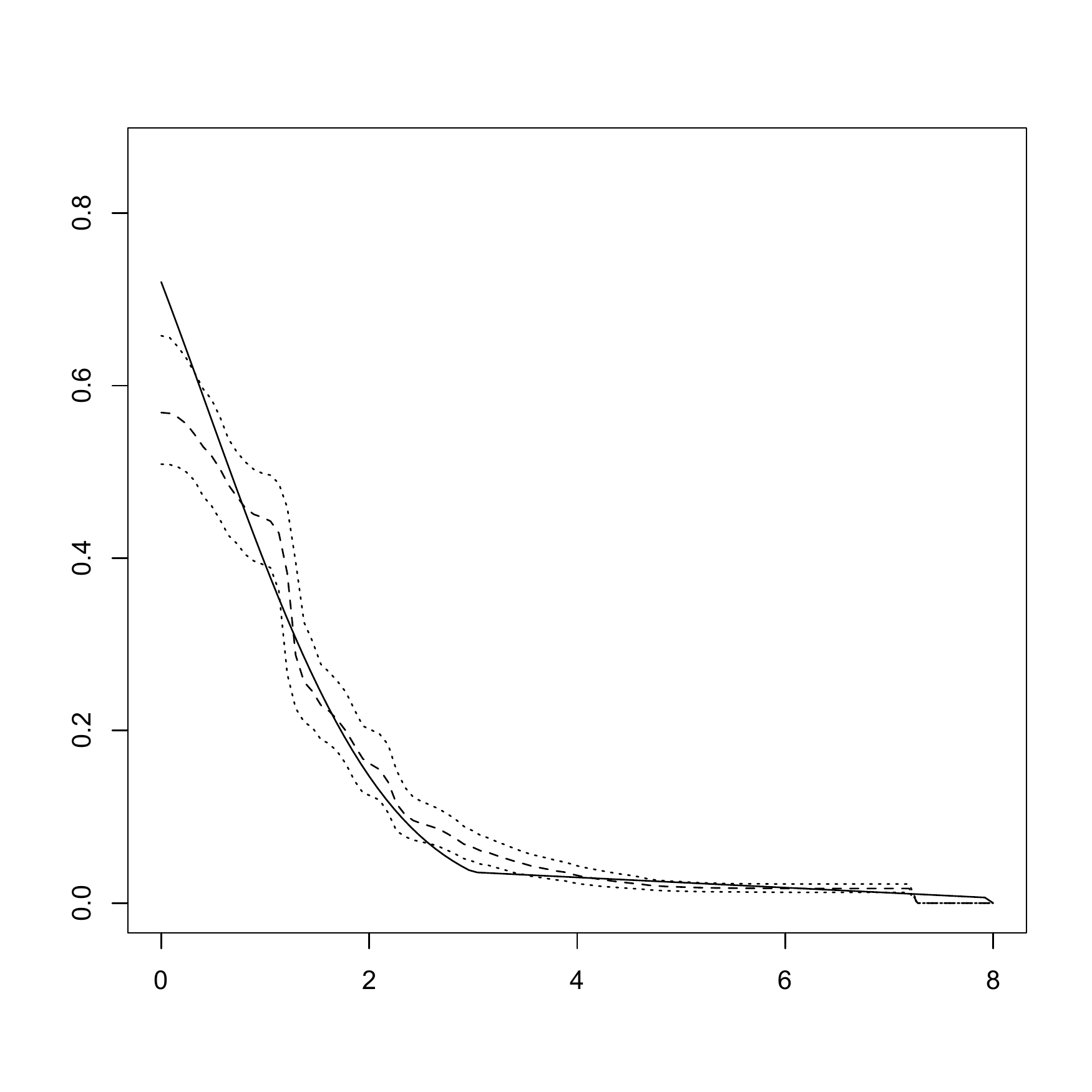} &
\includegraphics[width=0.33\textwidth]{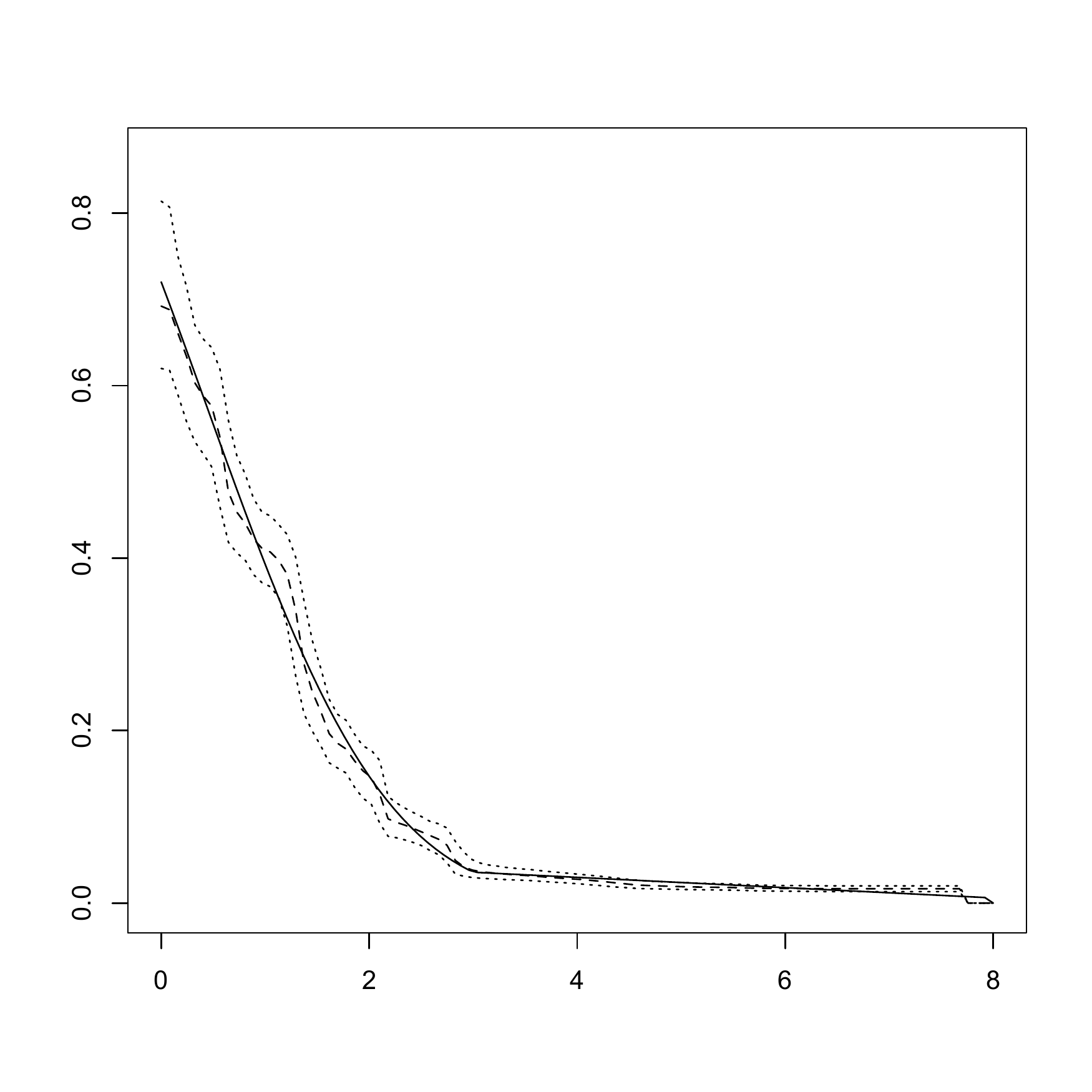} &
\includegraphics[width=0.33\textwidth]{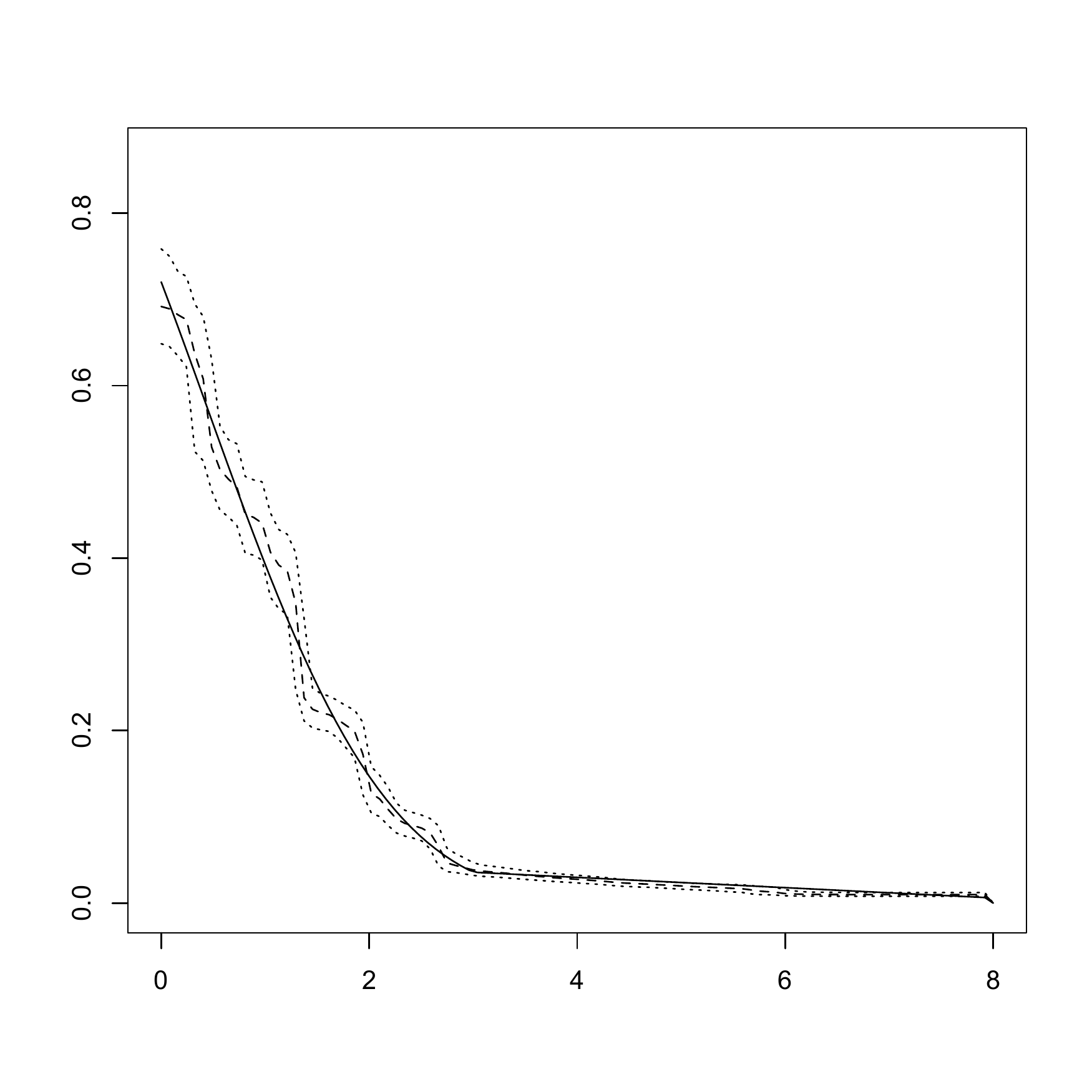}  
\end{tabular}
\caption{Estimation of $\overline{\lambda}_3$ from $D_{500}^3$ (first column), $D_{1000}^3$ (second column) and $D_{2000}^3$ (third column) using the four strategies: empirical prior (line 1), fixe $\gamma$ (line 2), hierarchical empirical prior (line 3), concentrated hierarchical empirical prior (line 4) . True density (plain line), estimation (dashed line) and confidence band (dotted lines)     }
\label{fig:estim3}
\end{figure}

To compare the three  different strategies  used to calibrate $\gamma$, several criteria are taken into account: tuning of the hyperparameters,  quality of the estimation, convergence of the MCMC, computational time. In terms of tuning effort on $\gamma$, the empirical Bayes  (strategy 1) together with the fixed $\gamma$ approach are comparable and significantly simpler than the hierarchical one, which increases the space to be explored by the MCMC algorithm and consequently slows down its convergence. Moreover,  setting an hyper-prior distribution on $\gamma$ requires to choose the parameters of this additional distribution ($a_{\gamma}$  and $b_{\gamma}$) and thus postpone the problem, even though these second order hyperparameters are supposed to be less influential. 
 In our particular example, computational time, for a fixed number of iterations here equal to $N_{iter} = 30000$, turned out to be also a key point. Indeed,  the simulation of the $\bar \lambda$  conditionally to the other variables involves an  accept-reject (AR) step. 
 For some particular values of $\gamma$ (small values of $\gamma$), we observe that the acceptance rate of the AR step can be dramatically low, automatically inflating the execution time of the algorithm, this phenomenon occurring randomly. The computational times (CpT) are summarized in  Table \ref{tab:num_comput_times}, which also provides for each of the $9$ datasets   the number of points  ($N_T$), the   $\hat{\gamma}$ computed using formula (\ref{eq:formula gamma emp}) and to be compared with the targeted value $\Psi^{-1}(E_{theo})$, the perturbation factor $\rho$ used in the fixed $\gamma$ strategy and the standard deviation of the prior  distribution of $\gamma$ $\sigma_{\gamma}$ (the prior mean being $\hat{\gamma}$) used in the two hierarchical approaches. The second hierarchical prior distribution (last columns of Table \ref{tab:num_comput_times})  corresponds to a prior distribution more concentrated around the empirical value  $\hat{\gamma}$. \\
 \noindent On Figure \ref{fig:estim1} (respectively \ref{fig:estim2} and \ref{fig:estim3}), we plot for each strategy and each dataset the posterior median of $\overline{\lambda}_1$ (respectively $\overline{\lambda}_2$ and $\overline{\lambda}_3$)  together with a   pointwise credible interval using the $0.1$\% and $0.9$\% empirical quantiles obtained from the posterior simulation. Table \ref{tab:dist} gives the distances between the estimated normalized intensities $\widehat{\overline{\lambda}}_i$ and the true $\overline{\lambda}_i$, for each dataset and each prior setting.

\noindent  On $\lambda_{0,1}$ (which a simple 2-steps function), the 4 strategies lead to the same quality of estimation (in term of losses / distances between the functions of interest).  In this case, it is thus interesting to have a look at the computational time in Table \ref{tab:num_comput_times}. We notice that for a small $\gamma_0$ or for a diffuse prior distribution on $\gamma$ (so possibly generating small values of $\gamma$) the computational time explodes. In practice, this phenomenon can be so critical that the user may have to stop the execution and re-launch the algorithm. Moreover, interestingly the  posterior mean of the number of non empty components in the mixture (computed  over the last 10000 iterations)  is equal in the case  $n= 500$ to $4.21$ in the empirical strategy, $11.42$ when $\gamma$ is fixed arbitrarily,  $6.98$ under the hierarchical large  prior and $3.77$ with the narrow hierarchical prior.  In this case choosing a small value of $\gamma$ leads to a posterior distribution on mixtures with too many non empty components.  This phenomena tends to disappear when $n$  increases. 
 
 For $\lambda_{0,2}$ and $\lambda_{0,3}$, a bad choice of $\gamma$ - here $\gamma$  too large in strategy 2 - or a not enough informative prior on $\gamma$ (hierarchical prior with large variance)  has a significant negative impact on the behavior of the posterior distribution. Contrariwise, the  medians of the empirical and the informative hierarchical posterior distribution of $\lambda$ have similar losses, as seen in Table  \ref{tab:dist}.

\begin{table}
\centering
\begin{tabular}{|ccc||cc|cc|cc|cc|}
\hline
&&& \multicolumn{2}{c|}{Empirical} & \multicolumn{2}{c|}{$\gamma$ fixed } &\multicolumn{2}{c|}{Hierarchical}& \multicolumn{2}{c|}{Hierarchical 2} \\
&&$N_T$&$\widehat{\gamma}$&CpT& $\rho   \Psi^{-1}(E_{theo})$ & CpT & $\sigma_{\gamma}$ & CpT&$\sigma_{\gamma}$ & CpT\\
\hline
\multirow{3}{*}{$\lambda_{0,1}$} & $D^1_{500}$ &499 &0.0386 & 523.57 &\multirow{3}{*}{$0.01\times  0.0323$} & \underline{2085.03} &\multirow{3}{*}{ 0.005} & \underline{12051.22} & \multirow{3}{*}{ 0.001}  & 447.75 \\ 
& $D^1_{1000}$   & 1036 & 0.0372 & 783.53 &  & \underline{1009.58} & & 791.28 &  & 773.33 \\ 
 &$D^1_{2000}$   & 2007 & 0.0372 & 1457.40 & & 1561.64 & & 1477.50 & & 1456.03 \\ 
 \hline
\multirow{3}{*}{$\lambda_{0,2}$} & $D^2_{500}$ & 505& 0.6605 & 1021.73 &\multirow{3}{*}{$100  \times   0.6667$}  & 1022.59 &\multirow{3}{*}{0.1} & 663.54 & \multirow{3}{*}{0.01}  & 1047.42 \\ 
 &$D^2_{1000}$  & 978 & 0.6857 & 1873.05 & & 1416.40 & & 1207.07 &  & 2018.89 \\ 
 &$D^2_{2000}$  & 2034 & 0.6827 & 4849.80 &  & 2236.02 && 2533.62 &  & 4644.55 \\ 
 \hline
  \multirow{3}{*}{$\lambda_{0,3}$} & $D^3_{500}$&483 & 0.4094 & 782.19 &  \multirow{3}{*}{$30 \times   0.4302$}  & 822.12 & \multirow{3}{*}{0.1}  & 788.14 & \multirow{3}{*}{0.01}  & 788.00 \\ 
  &$D^3_{1000}$& 1058 & 0.4398 & 1610.47 & & 2012.96 & & 1559.17 &  & 1494.75 \\ 
 &$D^3_{2000}$ & 2055 & 0.4677 & 3546.57 & & \underline{9256.71} &  & 3179.96 & & 2770.83 \\ 
   \hline
\end{tabular}
\caption{Values of $\gamma$ (or  the hyperparameters of its prior distribution)
 and computational time required for the Gibbs algorithm with $30000$ iterations. $N_T$: number of observations for each dataset $D^i_{n}$.  CpT: Computational Time. $\sigma_{\gamma}$: standard deviation of the prior distribution on $\gamma$}
\label{tab:num_comput_times}
\end{table}


\bigskip

\begin{table}
\centering
\begin{tabular}{|cc| lll | lll | lll | }
\hline
&& \multicolumn{3}{c|}{$\lambda_{0,1}$} & \multicolumn{3}{c|}{$\lambda_{0,2}$}  & \multicolumn{3}{c|}{$\lambda_{0,3}$} \\
&&$D^1_{500}$ & $D^1_{1000}$ &$D^1_{2000}$ &$D^2_{500}$ & $D^2_{1000}$& $D^2_{2000}$& $D^3_{500}$ &$D^3_{1000}$ & $D^3_{2000}$ \\
\hline

\multirow{4}{*}{$d_{L_1}$} &  Empir& 0.0246 & 0.0238 & 0.0207 & 0.0921$^*$ & 0.0817$^*$ & 0.0549$^*$ & 0.1382 & 0.0596 & 0.0606 \\ 
& Fixed  & 0.0161 & 0.0219 & 0.0211 & 0.5381 & 0.7221 & 0.6356 & 0.3114 & 0.2852 & 0.2885 \\ 
&Hierar & 0.0132 & 0.0233 & 0.0317 & 0.1082 & 0.1280 & 0.0969 & 0.2154 & 0.1378 & 0.1405 \\ 
& Hiera 2 & 0.0191 & 0.0240 & 0.0208 & 0.0925$^*$ & 0.0815$^*$ & 0.055$^*$2 & 0.1383 & 0.0607 & 0.0724 \\ 
\hline 
\multirow{4}{*}{$d_{L_2}$} &  Empir& 0.0014 & 0.0006 & 0.0006 & 0.0010 & 0.0010 & 0.0008 & 0.0008 & 0.0005 & 0.0006 \\ 
& Fixed & 0.0014 & 0.0006 & 0.0006 & 0.0251 & 0.0376 & 0.0378 & 0.0084 & 0.0095 & 0.0104 \\ 
&Hierar   & 0.0005 & 0.0006 & 0.0020 & 0.0033 & 0.0050 & 0.0034 & 0.0066 & 0.0070 & 0.0072 \\ 
& Hiera 2 & 0.0009 & 0.0006 & 0.0006 & 0.0011 & 0.0011 & 0.0009 & 0.0009 & 0.0005 & 0.0013 \\ 
\hline
\multirow{4}{*}{$\| \cdot \|_{\infty}$}  &Empir & 0.0909 & 0.0909 & 0.0909 & 0.1828 & 0.1369 & 0.0833 & 0.1421 & 0.0462 & 0.0428 \\ 
  & Fixed &  0.0909 & 0.0909 & 0.0909 & 0.3261 & 0.5022 & 0.4271 & 0.1707 & 0.1818 & 0.2020 \\ 
 &Hierar & 0.0909 & 0.0909 & 0.0909 & 0.2330 & 0.2007 & 0.1389 & 0.2154 & 0.1000 & 0.1190 \\ 
  &Hierar2 & 0.0909 & 0.0909 & 0.0909 & 0.1886 & 0.1413 & 0.0835 & 0.1515 & 0.0513 & 0.0555 \\
\hline
\end{tabular}
\caption{Distances between the estimated   and the true density  for all datasets  and all strategies. Distance $L_1$ in horizontal block one, distance $L_2$ in horizontal  block two, distance $\|\cdot \|_{\infty}$ in horizontal  block three.    }
\label{tab:dist}
\end{table}

\section{Proofs}\label{sec:proofs}
The notation  $\lesssim$ will be used to denote inequality up to a constant that is fixed throughout.
We denote by $C$ a constant depending on $m_1$, $m_2$, $k$ and so forth, which may change from line to line.

\subsection{Proof of Theorem \ref{Th:DPMGaussian1}
on empirical Bayes Dirichlet process mixtures of Gaussian densities for the ordinary smooth case
}
\label{pr:catia}

To prove Theorem \ref{Th:DPMGaussian1}, we must verify assumptions [A1] and [A2]. We first define the change of parameter $\psi_{\gamma, \gamma'}(p_{F,\sigma})$.
Under the Dirichlet process mixture prior we can write, when the base measure is associated to the parameter $\gamma = (m,\,s^2)$, that $p_{F,\sigma}(\cdot)=\sum_{j\geq1}p_j\phi_\sigma(\cdot-\theta_j)$ almost surely, with the $\theta_j\sim \Nor_{(m,s^2)}$ independently and independently of the $(p_j)_{j\geq1}$. Without loss of generality, we assume that $\mathcal K_n = [m_1,\,m_2]\times [s_1^2,\, s_n^2]$, with $-\infty <m_1\leq \,m_2<\infty$, $s_1^2>0$ and $s_n^2$ possibly going to infinity as a power of $\log n$.
Consider a $u_n$-covering of $[m_1,\,m_2]$ with intervals $I_k$, for $k=1,\,\ldots,\,L_{mn}$, where $L_{mn}=\lfloor (m_2-m_1)/u_n\rfloor$, and a covering of $[s_1^2,\,s_n^2]$ with intervals of the form $J_l= [s_1^2(1+u_n)^{l-1},\, s_1^2(1+u_n)^{l}]$, for $l=1,\,\ldots,\, L_{sn} $, where $L_{sn} \leq 2u_n^{-1}\log (s_n/s_1)$ and such that $L_{sn}\rightarrow\infty$ as $n\rightarrow\infty$.
We suppose that $u_n\rightarrow 0$.

For $s^2\in J_l$, $l=1,\,\ldots,\, L_{sn}$, let $\rho_l=(s^2/s_l^2)^{1/2} \leq (1+ u_n)^{1/2}$, with $s_l^2 = s_1^2 (1 + u_n)^{l-1} $, and, for $m \in I_k$, $ k=1,\,\ldots,\,L_{mn}$, writing $m_k = m_1 + (k-1)u_n$, if, for every $j\in\mathbb{N}$,
$\theta'_j \sim\textrm{N}_{(m_k,\,s_l^2)}=\textrm{N}_{\gamma'}$ and $\theta_j=\rho_\ell(\theta'_j-m_k)+m$, then $\theta_j \sim
\textrm{N}_{(m,\,s^2)}=\textrm{N}_{\gamma}$. Therefore, conditionally on $\sigma$, for $F\sim \mathrm{DP}(\alpha_\mathbb{R} \Nor_{\gamma})$,
\begin{equation}\label{eq:II}
\psi_{\gamma,\gamma'}(p_{F,\sigma})(x) =
\sum_{j\geq1}p_j\phi_{\sigma}(x - \theta'_j - [\theta'_j(\rho_l-1) + m-m_k\rho_l])
\end{equation}
is distributed according to a Dirichlet process location mixture of Gaussian distributions, with base measure $\Nor_{\gamma'}$.
The following  inequalities are used repeatedly in the sequel,
 $$ (1 + u_n)^{1/2} \geq \rho_l \geq 1, \qquad  - m_k u_n  \leq m-\rho_lm_k \leq u_n.$$
With abuse of notation, we also denote by $\psi_{\gamma, \gamma'} (\theta_j) = \rho_\ell(\theta'_j-m_k)+m$ for any $\gamma' = (m_k,\, s_l^2)$ and $\gamma= (m,\,s^2)$.
We first verify assumption [A1]. From condition \eqref{Holder}, let $\sigma \in (\sigma_n /2,\, 2 \sigma_n)$, with $\sigma_n = \epsilon_n^{1/\beta}$,
and $F^* $ be as defined in \eqref{Holder}, so that $F^* = \sum_{j=1}^{N_\sigma} p_j^* \delta_{\theta_j^*}$ with the minimal distance between any two $\theta_j$'s bounded from below by $\delta_j = \sigma \epsilon_n^{2b}$ for some $b>0$. Construct a partition $(U_j)_{j=1}^M$ of $\R$ following  the proof of Theorem 4 of \citet{ghosal:shen:tokdar}. Let $a_\sigma = a_0 |\log \sigma|^{1/\tau}$. The partition is such that $(U_j)_{j=1}^K$ is a partition of  $[-a_\sigma, a_\sigma]$ composed of intervals in the form $[\theta_j^* - \delta_j/2,\, \theta_j^* + \delta_j/2]$, $j=1,\, \ldots,\, N_\sigma$, and of intervals with diameter smaller than or equal to $\sigma$ to complete $[-a_\sigma,\, a_\sigma]$. Then, construct a partition of $(-\infty,\, -a_\sigma)$ and $(a_\sigma, \,\infty)$ with intervals $U_j$, $j> K$, such that $1 \geq \Nor_\gamma(U_j) \geq \sigma \epsilon_n^{2b}$. Note that, as in \citet{ghosal:shen:tokdar}, $M \lesssim \sigma^{-1} (\log n)^{1+ 1/\tau}$ and that, for all $j \leq K$,
$\Nor_\gamma (U_j) \gtrsim \frac{ \delta_j}{s} e^{-2 a_\sigma^2 /s^2 } \gtrsim \epsilon_n^{b'}$
for some $b'> b$ uniformly in $\gamma \in \mathcal K_n$. As in \citet{ghosal:shen:tokdar},
define $B_n^{\gamma}$ as the set  of $(F,\, \sigma)$ such that $\sigma \in (\sigma_n/2, 2\sigma_n)$ and
$$\sum_{j=1}^M |F(U_j)- p_j^*| \leq 2 \epsilon_n^{2b'}, \qquad 
\min_{1\leq j\leq M} \Nor_\gamma(U_j) \geq \epsilon_n^{4b'}/2.$$
Following Lemma 10 of \citet{ghosal:vdv:mixture:07}, we obtain that for some constant $c$ independent of $\gamma \in \mathcal K_n$,
\begin{equation}\label{piBn:DPM}
\inf_{\gamma\in\mathcal{K}_n}\pi\left( B_n^\gamma\mid\gamma\right)\gtrsim e^{-c \sigma^{-1} (\log n)^{2+1/\tau}}.
\end{equation}
Moreover, for all $(F,\, \sigma) \in B_n^\gamma$,
with $\gamma' = (m_k,\,s_l^2)$ and any $\gamma \in I_k\times J_l$,
\begin{equation*}
\begin{split}
p_{F,\sigma}(x) & =\sum_{j\geq1}p_j\phi_{\sigma}(x - \psi_{\gamma, \gamma'}(\theta_j))
 \geq \sum_{j\geq 1} \1_{|\theta_j'|\leq a_\sigma} p_j\phi_{\sigma}(x - \psi_{\gamma, \gamma'}(\theta_j)) \\
&\geq  \sum_{j\geq 1} \1_{|\theta_j'|\leq a_\sigma} p_j\phi_{\sigma}(x - \theta_j') e^{-4\frac{ |x-\theta_j'|(a_\sigma+1)u_n +(a_\sigma^2 +1) u_n^2}{ \sigma_n^2}}.
\end{split}
\end{equation*}
Note that $n^{-1}\sigma_n^{-1}=\epsilon^2_n$. Choosing $u_n \lesssim n^{-2}\sigma_n^{-1} (\log n)^{-1/\tau}=n^{-1}\epsilon_n^2(\log n)^{-1/\tau}$, we obtain that on the event $A_n = \{ \sum_{i=1}^n |X_i-m_0| \leq n \tau_0^2k_n\}$, for $k_n\geq (\log n)^{1/\tau}$, using the inequality
$\log x\geq (x-1)/x$ valid for $x>0$,
\[
\begin{split}
\ell_n(p_{F,\sigma}) - \ell_n(p_0)  &\geq \ell_n(p_{F_n, \sigma}) - \ell_n(p_0) + n\log c_\sigma-4n\sigma_n^{-2}[(a_\sigma^2+1)u_n^2+(a_\sigma +1)u_n(2a_\sigma+\tau_0^2k_n)]\\
&\geq
\ell_n(p_{F_n, \sigma}) - \ell_n(p_0) +n\log c_\sigma-C'n\epsilon_n^2\\
&\geq \ell_n(p_{F_n, \sigma}) - \ell_n(p_0) +n(c_\sigma-1)-C'n\epsilon_n^2\\
&\geq \ell_n(p_{F_n, \sigma}) - \ell_n(p_0) -2n\epsilon_n^{2b'}-C'n\epsilon_n^2
\geq \ell_n(p_{F_n, \sigma}) - \ell_n(p_0) -n\epsilon_n^2-C'n\epsilon_n^2,
\end{split}
\]
with a constant $C'$ large enough, where $p_{F_n,\sigma} (x) =  c_\sigma^{-1}\sum_{j\geq 1} \1_{|\theta_j'|\leq a_\sigma} p_j\phi_{\sigma}(x - \theta_j') $ and $c_\sigma=\sum_{j:\,|\theta_j'|\leq a_\sigma} p_j\geq (1-2\epsilon_n^{2b'})\geq (1-\epsilon_n^2)$ because $b'>1$. The proof of Theorem 4 of \citet{ghosal:shen:tokdar}, together with condition \eqref{Holder}, implies that [A1] is satisfied. We now prove [A2].
As in Proposition 2 of \citet{ghosal:shen:tokdar}, consider
\begin{equation}\label{eq:sieve}
\mathcal F_n = \pg{ (F,\, \sigma):\, \underline \sigma_n \leq \sigma \leq \bar \sigma_n,\,\,\, F = \sum_{j=1}^{\infty} p_j \delta_{\theta_j}, \,\,\, |\theta_j| \leq  n^{1/2} \,\,\, \forall\, j \leq H_n,\,\,\, \sum_{j>H_n}p_j \leq \epsilon_n},
\end{equation}
with $\underline{\sigma}_n = \epsilon_n^{1/\beta}$, $\bar \sigma_n = \exp( t n \epsilon_n^2)$ for some $t >0$ depending on the parameters $(\nu_1,\, \nu_2)$  of the inverse-gamma distribution for $\sigma$, and $H_n =\lfloor  n \epsilon_n^2 /\log n\rfloor$.
For some $x_0>0$, let $a_n=2x_0(\log n)^{1/\tau}$, $\gamma'= (m_k,\,s_l^2) \in \mathcal K_n$,
$ \gamma \in I_k \times J_l$ and $|x|\leq a_n/2 $. If $|\theta| \geq a_n$ then
$|x-\theta| \geq |\theta|/2$ and, for $ u_n \leq n^{-2}$, we can bound $p_{F,\sigma}$ as follows:
\begin{equation*}
\begin{split}
p_{F,\sigma}(x)&=\int_\R \phi_\sigma ( x - \psi_{\gamma, \gamma'}(\theta) )\mathrm{d}F'(\theta) \leq \int_\R \phi_\sigma (x-\theta) e^{\frac{|x-\theta|  u_n (|\theta | +1) }{\sigma^2 } } \d F'(\theta)\\
&\leq e^{3a_n( a_n+1) u_n\sigma^{-2}/2} \int_{|\theta|\leq a_n} \phi_\sigma (x-\theta) \d F'(\theta) +
\int_{|\theta|>a_n} \phi_{\sigma}( (x-\theta)(1 - 8u_n))\d F'(\theta)\\
 &\leq e^{o(1/n)} \int_{|\theta|\leq a_n} \phi_\sigma (x-\theta)  \d F'(\theta) +(1+ O(1/n))\int_{|\theta|>a_n} \phi_{\tilde \sigma_n}(x-\theta)  \d F'(\theta)\\
 &= e^{o(1/n)} \left( \int_{|\theta|\leq a_n} \phi_\sigma (x-\theta) \, \mathrm{d}F'(\theta) +\int_{|\theta|>a_n} \phi_{\tilde \sigma_n}(x-\theta) \d F'(\theta) \right),
\end{split}
\end{equation*}
where $F'\sim\mathrm{DP}(\alpha_{\mathbb{R}}\mathrm{N}_{\gamma'})$, $\tilde \sigma_n = \sigma (1 -1/n^2)^{-1/2}$. Now, since $\P_{p_0}^{(n)}\left(\Omega_n \right) = 1 + o(1) $, with $ \Omega_n =\{- a_n/2 \leq \min_i X_i \leq \max_i X_i \leq a_n/2 \} $, we can replace $\mathbb{R}^n$ with $\Omega_n$ and define for all $(F,\,\sigma) \in \mathcal F_n$,
\begin{equation*}\label{qgamma:DPM}
q_{\gamma}^{F,\sigma}(x)= \1_{[-a_n/2,\, a_n/2]}(x) e^{o(1/n)} \left( \int_{|\theta|\leq a_n} \phi_\sigma (x-\theta) \d F(\theta)  +\int_{|\theta|>a_n} \phi_{\tilde \sigma_n}(x-\theta) \d F(\theta)  \right).
\end{equation*}
 Similarly, we can replace $p_{F,\sigma}$ by $p_{F,\sigma}\1_{[-a_n/2,\, a_n/2]}$. We then have that $q_\gamma^{F,\sigma}$ and $p_{F,\sigma}$ are contiguous and
 $\| q_\gamma^{F,\sigma} - p_{F,\sigma} \|_1 = o(1/n)$. Therefore,  we can consider the same tests as in Lemma 1 of \citet{ghosal:vdv:mixture:07} and \eqref{test} is verified, together with \eqref{eq:entropy} using Proposition 2 of \citet{ghosal:shen:tokdar}. Since \eqref{piBn:DPM} implies condition \eqref{ratio:mass}, there only remains to verify assumption \eqref{Thetanc}. The difficulty here is to control $q_\gamma^{F,\sigma}$ as $\sigma\rightarrow0$. Indeed, $q_\gamma^{F,\sigma}$ can be used as an upper bound on $\psi_{\gamma, \gamma'}(p_{F,\sigma})$ for all $(F,\,\sigma)$ with $\sigma > \underline \sigma_n$, when $|x|\leq a_n/2$, which we are allowed to consider since we can restrict ourselves to the event $\Omega_n$. Thus,
$\int_{\sigma>\bar \sigma_n} Q_\gamma^{F,\sigma}([-a_n/2,\,a_n/2]^n)\,\d\pi(F\mid\gamma)\d\pi(\sigma) = o(e^{-2n\epsilon_n^2})\pi(B_n^\gamma\mid\gamma)$ uniformly in $\gamma \in \mathcal K_n$. We now study
\begin{equation*}\label{Thetanc2}
\int_{\sigma< \underline \sigma_n} Q_\gamma^{F,\sigma}([-a_n/2,\,a_n/2]^n)\,\d\pi(F\mid\gamma)\d\pi(\sigma).
\end{equation*}
We split $(0,\,\underline \sigma_n) = \cup_{j=0}^\infty [\underline \sigma_n2^{-(j+1)},\,\underline \sigma_n2^{-j})$. Then, for all $\sigma \in [\underline \sigma_n2^{-(j+1)},\,\underline \sigma_n2^{-j})$, $j\geq 0$, define $u_{n,j} = n^{-1}e_n (\underline \sigma_n 2^{-j})^2$,
with $e_n =o(1)$ and $\beta>1/2$ to guarantee that $u_{n,j}\leq n^{-2}$, so that, similarly to before, for all $\gamma \in \mathcal K_n$,
\begin{equation*}
\sup_{\|\gamma-\gamma'\| \leq u_{n,j}}\psi_{\gamma, \gamma'}(p_{F,\sigma}(x))\1_{[-a_n/2,\,a_n/2]}(x) \leq q_{\gamma}^{F,\sigma}(x),
\end{equation*}
with the distance $\|\gamma-\gamma'\|=|m-m_k|+|s/s_\ell- 1|$.
For all $\gamma \in \mathcal K_n$, using a $u_{n,j}$-covering of $\{\gamma':\, \|\gamma -\gamma'\|\leq u_n\}$
with
centering points $\gamma_i$, $i=1,\,\ldots,\, N_j$, where $N_j\leq (u_n/u_{n,j})^2 $, we have
 \begin{equation*}
 \begin{split}
q_\gamma^{F,\sigma}(x) &\leq \max_{1\leq i \leq N_j} \int_\R \sup_{\| \gamma' - \gamma_i\| \leq u_{n,j} } \phi_{\sigma}(x - \psi_{\gamma, \gamma'}(\theta) )\d F(\theta) \\
  & \leq \max_{1\leq i \leq N_j}  \int_\R \phi_\sigma ( x - \psi_{\gamma, \gamma_i}(\theta) ) e^{\frac{|x - \psi_{\gamma, \gamma_i}(\theta)| u_{n,j}(1 +|\psi_{\gamma,\gamma_i}(\theta)| ) }{\sigma^2} } \d F(\theta) \leq \max_{1\leq i \leq N_j}  c_{n,i} g_{\sigma, i }(x),
  \end{split}
  \end{equation*}
 with $g_{\sigma, i }$ the probability density over $[-a_n/2,\, a_n/2]$ proportional to
$$\int_{|\theta|\leq a_n}  \phi_\sigma ( x - \psi_{\gamma, \gamma_i}(\theta) ) e^{\frac{2a_n u_{n,j}(1 +a_n)}{\sigma^2}} \d F(\theta)
+ \int_{|\theta|> a_n}  \phi_\sigma ( (x - \psi_{\gamma, \gamma_j}(\theta))( 1 - o(1/n))) \d F(\theta),$$
so that
$$c_{n,i} \leq F([-a_n,\, a_n])  e^{16x_0 e_n (\log n)^{1/\tau}/n}+ (1-F([-a_n,\, a_n]))(1 + o(1/n^{1/2})) = 1+ O(e_n (\log n)^{1/\tau}/n).$$
This implies that, for $e_n \leq (\log n)^{-1/\tau}$,
\begin{equation*}
\begin{split}
\int_{\underline \sigma_n 2^{-(j+1)}}^{\underline \sigma_n 2^{-j} } Q_\gamma^{F, \sigma}([-a_n/2,\, a_n/2]^n)\d\pi(F\mid \gamma)\d\pi(\sigma)
& \lesssim N_j  \pi([\underline \sigma_n2^{-(j+1)},\, \underline \sigma_n2^{-j}))\\ & \lesssim u_n^2n^{2} \underline \sigma_n^{-4}e_n^{-2} 2^{4j}
e^{-2^{j-1}/\underline \sigma_n},
\end{split}
\end{equation*}
whence $
 \int_{\sigma< \underline \sigma_n} Q_\gamma^{F,\sigma}([-a_n/2,\,a_n/2]^n) \d\pi(F\mid \gamma)\d\pi(\sigma)\lesssim e^{-\underline \sigma_n^{-1}/2}$ and \eqref{Thetanc} is verified, which completes the proof of Theorem \ref{Th:DPMGaussian1}.


\subsection{Proof of Theorem \ref{Th:DPMGaussian2} on empirical Bayes Dirichlet process mixtures
of Gaussian densities for the super-smooth case}
\label{pr:catia2}
The main difference with the ordinary smooth case lies in the fact that, since the rate $\epsilon_n$ is almost parametric, we have $n\epsilon_n^2=O((\log n)^\kappa)$ for a suitable finite constant $\kappa>0$ so that, in order to compensate the number $N_n(u_n)$ of points, we need that, for some set $\tilde B_n$,
\begin{equation}\label{eq:expbound}
\sup_{\gamma\in\mathcal K_n}\sup_{p_{F,\sigma}\in \tilde B_n}\mathbb{P}_{p_0}^{(n)}\pt{\inf_{\|\gamma-\gamma'\|\leq u_n}\ell_n(\psi_{\gamma,\gamma'}(p_{F,\sigma}))-\ell_n(p_0)<-n\epsilon_n^2}\lesssim e^{-Cn\epsilon_n^2}
\end{equation}
for some constant $C>0$. It is known from Lemma 2 of \citet{shen:wasserman:01} that if $p_{F,\sigma}\in S_n=\{p_{F,\sigma}:\,\rho_\alpha(p_0;\,p_{F,\sigma})\leq\epsilon_n^2\}$ then
$$\mathbb{P}_{p_0}^{(n)}\pt{\ell_n(p_{F,\sigma})-\ell_n(p_0)<-n(1+C)\epsilon_n^2}\lesssim e^{-\alpha Cn\epsilon_n^2}.$$
Consider the same set $\mathcal K_n$ and the same $u_n$-covering as in the proof of Theorem \ref{Th:DPMGaussian1}.
Take $\sigma_n =O((\log n)^{-1/r})$. Let $\sigma \in (\sigma_n,\,\sigma_n+ e^{-d_1(1/\sigma_n)^{r}})$ for some positive constant $d_1$.
It is known from Lemma 6.3 of \citet{scricciolo:12} that there exists a distribution $F^*=\sum_{j=1}^{N_\sigma} p_j^* \delta_{\theta_j^*}$, with $N_\sigma=O((a_\sigma/\sigma)^2)$ points in $[-a_\sigma,\, a_\sigma]$, where
$a_\sigma=O(\sigma^{-r/(\tau\wedge2)})$, such that, for some constant $c>0$,
$n^{-1}\max\{\textrm{KL}(p_0;\, p_{F^*, \sigma}) ,\,V_2 (p_0;\, p_{F^*, \sigma})\}\lesssim e^{-c(1/\sigma)^{r}}$.
Inspection of the proof of Lemma 6.3 reveals that all arguments remain valid to
bound above any $\rho_\alpha$-divergence $\rho_\alpha(p_0;\,p_{F^*, \sigma})$ for $\alpha\in(0,\,1]$.
In fact, using the inequality valid for all $a,\,b>0$, $|a^\alpha-b^\alpha|\leq|a^\beta-b^\beta|^{\alpha/\beta}$, with $0\leq\alpha\leq\beta$, having set in our case $\beta=1$, we have $\rho_\alpha(p_0;\,p_{F^*,\sigma})\leq \alpha^{-1}\mathbb{E}_{p_0}\pg{\pq{|p_0(X_1)-p_{F^*, \sigma}(X_1)|/p_{F^*, \sigma}(X_1)}^\alpha}$.
All bounds used in the proof of Lemma 6.3 for the various pieces in which
$n^{-1}\textrm{KL}(p_0;\, p_{F^*, \sigma})$ is split can be used here to bound above $\rho_\alpha(p_0;\,p_{F^*, \sigma})$.
Thus, $\rho_\alpha(p_0;\,p_{F^*, \sigma})\lesssim e^{-c(1/\sigma)^{r}}$. Construct a partition $(U_j)_{j=0}^{N_\sigma}$ of $\R$, with $U_0:=(\bigcup_{j=1}^{N_{\sigma}}U_j)^c$,
$U_j\ni\theta_j^*$ and $\lambda(U_j)=O(e^{-c_1(1/\sigma)^{r}})$, $j=1,\,\ldots,\,N_{\sigma}$. Then,
$\inf_{\gamma\in\mathcal K_n}\min_{1\leq j\leq N_\sigma}\textrm{N}_\gamma(U_j)\gtrsim
e^{-c_1(1/\sigma)^{r}}$. Defined the set
$$B_n:=\pg{(F,\,\sigma):\,\sigma\in(\sigma_n,\,\sigma_n+ e^{-d_1(1/\sigma_n)^{r}}), \,\,\,
\sum_{j=1}^{N_\sigma} |F(U_j)- p_j^*| \leq e^{-c_1(1/\sigma)^{r}}},$$
we have $\inf_{\gamma\in\mathcal{K}_n}\pi(B_n\mid\gamma)\gtrsim \exp{\{-c_2N_{\sigma_n}(1/\sigma_n)^{r}\}}=\exp{\{-c_3(\log n)^{5+6/{r}}\}}=e^{-c_3n\epsilon_n^2}$
and, for every $(F,\,\sigma)\in B_n$, it results $\rho_\alpha(p_0;\,p_{F,\sigma})\lesssim e^{-c_2(1/\sigma_n)^{r}}\lesssim \epsilon_n^2$. Moreover, reasoning as in the proof of Theorem \ref{Th:DPMGaussian1}, for $u_n\lesssim k_n^{-1}\sigma_n^2\epsilon_n^2(\log n)^{-1/(\tau\wedge 2)}$, on the event $A_n=\{ \sum_{i=1}^n |X_i-m_0| \leq n\tau_0^2k_n\}$, for $k_n=O(n)$,
\[
\begin{split}
\ell_n(p_{F,\sigma}) - \ell_n(p_0)  &\geq \ell_n(p_{F_n, \sigma}) - \ell_n(p_0) + n(c_\sigma-1)\\
&\quad -4n\sigma_n^{-2}[(a_\sigma^2+1)u_n^2+(a_\sigma +1)u_n(2a_\sigma+\tau_0^2k_n)]\\ 
&\geq\ell_n(p_{F_n, \sigma}) - \ell_n(p_0) -n\epsilon_n^2-C'n\epsilon_n^2
\end{split}
\]
for some positive constant $C'$, with
$p_{F_n,\sigma} (\cdot ) =  c_\sigma^{-1}\sum_{j\geq 1} \1_{|\theta_j'|\leq a_\sigma} p_j\phi_{\sigma}(x - \theta_j') $, where $c_\sigma=\sum_{j:\,|\theta_j'|\leq a_\sigma} p_j\geq (1-e^{-c_1(1/\sigma_n)^{r}})\geq 1-\epsilon_n^2$.
Hence, the proof of Lemma 2 of \citet{shen:wasserman:01} implies that \eqref{eq:expbound} is satisfied.
The other parts of the proof of Theorem \ref{Th:DPMGaussian1} go through to this case.
We only need to verify that both $\mathbb{P}_{p_0}^{(n)}(A_n^c)$  and $\mathbb{P}_{p_0}^{(n)}(\Omega_n^c)$
go to $0$. Indeed, by Markov's inequality, $\mathbb{P}_{p_0}^{(n)}(A_n^c)\leq n^{-1}\leq e^{-c_3n\epsilon_n^2}$. Also, $\mathbb{P}_{p_0}^{(n)}(\Omega_n^c)\lesssim e^{-c_4n\epsilon_n^2}$ because of the assumption \eqref{cond:tail1} that $p_0$ has exponentially small tails.

\subsection{Proof of Corollary \ref{Mixing} on empirical Bayes posterior contraction rates for mixing distributions in Wasserstein
metrics}
We appeal to Corollary~\ref{Mixing} in \cite{scricciolo:12} and the following remark which gives an indication on
how to remove the condition that $\Theta$ is bounded. In particular, we need that, for every $1\leq q<\infty$, there exist
$q<u<\infty$ and $0<B<\infty$ such that $\mathbb{E}_F[|X|^u]<B$ with
$[\mathrm{DP}(\alpha_{\mathbb{R}}\mathrm{N}_{\hat{\gamma}_n})]$-probability one, for almost every sample path when sampling from $\mathbb{P}_{p_0}^{(\infty)}$.
This can be proved appealing to the properties of the tails of the distribution functions sampled
from a Dirichlet process as in \citet{doss:sellke:82}.

\subsection{Proof of Corollary~\ref{Mixingdensity} on adaptive empirical Bayes density deconvolution}
The result is based on the following inversion inequalities which relate the $\mathbb{L}_2$-distance between
the true mixing density and the random approximating mixing density to the $\mathbb{L}_2$/$\mathbb{L}_1$-distance between the corresponding mixed densities:
\begin{equation*}\label{Ka}
\|p_Y-p_{0Y}\|_2\lesssim
\left\{
\begin{array}{ll}
\|K\ast p_Y-K\ast p_{0Y}\|_2^{\beta/(\beta+\eta)}, & \hbox{ ordinary smooth case,}\\[5pt]
(-\log \|K\ast p_Y-K\ast p_{0Y}\|_1)^{-\beta/r}, & \hbox{ super-smooth case.}
\end{array}
\right.
\end{equation*}
In what follows, we use \vir{os} and \vir{ss} as short-hands for \vir{ordinary smooth} and \vir{super-smooth}, respectively.
To prove the preceding inequalities, we instrumentally use the \emph{sinc} kernel
to characterize regular densities in terms of their approximation properties.
We recall that the \emph{sinc} kernel
\[\operatorname{sinc}(x)=
\left\{
  \begin{array}{cl}
  (\sin x)/(\pi x), & \hbox{\mbox{ if}\,\,\, $x\neq0$,}\\
  1/\pi, & \hbox{\mbox{ if}\,\,\, $x=0$,}
  \end{array}
\right.
\]
has Fourier transform $\widehat{\operatorname{sinc}}$ identically equal to $1$
on $[-1,\,1]$ and vanishing outside it.
For $\delta>0$, let $\operatorname{sinc}_\delta(\cdot)=\operatorname{sinc}(\cdot/\delta)$ and
define $g_\delta$ as the inverse Fourier transform of $\widehat{\operatorname{sinc}_\delta}/\widehat{K}$,
\[g_\delta(x)=\frac{1}{2\pi}\int e^{-itx}\frac{\widehat{\operatorname{sinc}_\delta}(t)}{\widehat{K}(t)}\d t,\quad x\in\mathbb{R}.\]
Let $\widehat{g_\delta}=\widehat{\operatorname{sinc}_\delta}/\widehat{K}$
be the Fourier transform of $g_\delta$. So, $\operatorname{sinc}_\delta=K\ast g_\delta$ and
$p_Y\ast \operatorname{sinc}_\delta=(p_Y\ast K)\ast g_\delta=(K\ast p_Y)\ast g_\delta$.
We then have
\[
\begin{split}
\|p_Y-p_{0Y}\|_2^2&\leq
\|p_Y\ast\operatorname{sinc}_\delta-p_{0Y}\ast \operatorname{sinc}_\delta\|_2^2
+
\|p_Y-p_Y\ast \operatorname{sinc}_\delta\|_2^2
+
\|p_{0Y}-p_{0Y}\ast\operatorname{sinc}_\delta\|_2^2\\
&\lesssim
\|p_Y\ast\operatorname{sinc}_\delta-p_{0Y}\ast \operatorname{sinc}_\delta\|_2^2+\|p_Y-p_Y\ast \operatorname{sinc}_\delta\|_2^2+ \delta^{2\beta}
\end{split}
\]
because
$
\|p_{0Y}-p_{0Y}\ast \operatorname{sinc}_\delta\|_2^2=\int_{-\infty}^\infty|\widehat{p_{0Y}}(t)|^2|1-\operatorname{sinc}_\delta(t)|^2\d t<
\delta^{2\beta}\int_{|t|>1/\delta}(1+t^2)^\beta|\widehat{p_{0Y}}(t)|^2\d t\lesssim \delta^{2\beta}$ by assumption \eqref{Fourierpoly}.
Now, recall that $p_Y= p_{F,\sigma}=F\ast\phi_\sigma$. For $(F,\,\sigma)\in\mathcal{F}_n$ defined as in \eqref{eq:sieve}, with
$\sigma\geq\underline{\sigma}_n\propto C\delta(\log n)^{\kappa_2}$, where $2\kappa_2\geq1$,
$C^2/2\geq (2\beta+1)/[2(\beta+\eta)+1]$ and $\delta\gtrsim n^{-1/[2(\beta+\eta)+1]}$, we have
\[\begin{split}
\|p_Y-p_Y\ast \operatorname{sinc}_\delta\|_2^2
&=
\int_{|t|>1/\delta}|\widehat{p_Y}(t)|^2\d t  =
\int_{|t|>1/\delta}|\widehat{F}(t)|^2|\widehat{\phi_\sigma}(t)|^2\d t
\leq
\int_{|t|>1/\delta}|\widehat{\phi_\sigma}(t)|^2\d t \\
&\lesssim (\sigma^2/\delta)^{-1}e^{-(\sigma/\delta)^2/2}\lesssim
[\delta(\log n)^{2\kappa_2}]^{-1}e^{-C^2(\log n)^{2\kappa_2}/2}\lesssim \delta^{2\beta}.\end{split}\]
In the ordinary smooth case,
\[\begin{split}
\|p_Y\ast\operatorname{sinc}_\delta-p_{0Y}\ast \operatorname{sinc}_\delta\|_2^2&=
\|(K\ast p_Y)\ast g_\delta-(K\ast p_{0Y})\ast g_\delta\|_2^2\\&\leq \delta^{-2\eta}
\int_{|t|\leq1/\delta}|\widehat{K}(t)|^2|\widehat{p_Y}(t)-\widehat{p_{0Y}}(t)|^2\d t\\
&\leq\delta^{-2\eta}
\|K\ast p_Y-K\ast p_{0Y}\|_2^2.
\end{split}
\]
In the super-smooth case,
\[
\|p_Y\ast\operatorname{sinc}_\delta-p_{0Y}\ast \operatorname{sinc}_\delta\|_2^2=
\|(K\ast p_Y)\ast g_\delta-(K\ast p_{0Y})\ast g_\delta\|_2^2\leq\|K\ast p_Y-K\ast p_{0Y}\|_1^2\|g_\delta\|_2^2,
\]
where
\[
\|g_\delta\|_2^2=\frac{1}{2\pi}\int_{-\infty}^\infty \frac{|\widehat{\operatorname{sinc}_\delta}(t)|^2}{|\widehat{K}(t)|^2}\d t
=\frac{1}{2\pi}\int_{\delta|t|\leq1}|\widehat{K}(t)|^{-2}\d t\lesssim e^{2\varrho\delta^{-r}}.
\]
Combining pieces, for $(F,\,\sigma)\in\mathcal F_n$,
\[\|p_Y-p_{0Y}\|_2^2\lesssim \delta^{2\beta}+
\left\{
\begin{array}{ll}
\|K\ast p_Y-K\ast p_{0Y}\|_2^2\times\delta^{-2\eta}, & \hbox{ os case,}\\[5pt]
\|K\ast p_Y-K\ast p_{0Y}\|_1^2 \times e^{2\varrho\delta^{-r}}, & \hbox{ ss case,}
\end{array}
\right.
\]
so that the optimal choice for $\delta$ is
\[\delta=
\left\{
\begin{array}{ll}
O(\|K\ast p_Y-K\ast p_{0Y}\|_2^{1/(\beta+\eta)}), & \hbox{ os case,}\\[5pt]
O\pt{(-\log \|K\ast p_Y-K\ast p_{0Y}\|_1)^{-1/r}}, & \hbox{ ss case.}
\end{array}
\right.
\]
Taking into account that $p_Y= p_{F,\sigma}=F\ast\phi_\sigma$, for $1\leq q\leq\infty$, we have
\[\|K\ast p_Y-K\ast p_{0Y}\|_q=\|K\ast F\ast\phi_\sigma-K\ast p_{0Y})\|_q=
\|(K\ast F)\ast\phi_\sigma-K\ast p_{0Y}\|_q=\|p_{F\ast K,\sigma}-K\ast p_{0Y}\|_q.\]
Then,
\[\|p_Y-p_{0Y}\|_2=\|p_{F,\sigma}-p_{0Y}\|_2\lesssim
\left\{
\begin{array}{ll}
\|p_{F\ast K,\sigma}-K\ast p_{0Y}\|_2^{\beta/(\beta+\eta)}, & \hbox{ os case,}\\[5pt]
(-\log \|p_{F\ast K,\sigma}-K\ast p_{0Y}\|_1)^{-\beta/r}, & \hbox{ ss case.}
\end{array}
\right.
\]
For suitable constants $\tau_1,\,\kappa_1>0$, let
\[\psi_n=
\left\{
\begin{array}{ll}
n^{-(\beta+\eta)/[2(\beta+\eta)+1]}(\log n)^{\kappa_1}, & \hbox{ os case,}\\[5pt]
n^{-1/2}(\log n)^{\tau_1}, & \hbox{ ss case,}
\end{array}
\right.
\]
and let $\epsilon_n$ be as in the statement. Hence, for all $(F,\,\sigma)\in\mathcal F_n$, the following inclusions hold:
\[\left\{
\begin{array}{ll}
\{(F,\,\sigma):\,\|p_{F\ast K,\sigma}-K\ast p_{0Y}\|_2\lesssim \psi_n\}\subseteq
\{(F,\,\sigma):\,\|p_Y-p_{0Y}\|_2\lesssim\epsilon_n\}, & \hbox{ os case,}\\[5pt]
\{(F,\,\sigma):\,\|p_{F\ast K,\sigma}-K\ast p_{0Y}\|_1\lesssim \psi_n\}\subseteq
\{(F,\,\sigma):\,\|p_Y-p_{0Y}\|_2\lesssim\epsilon_n\}, & \hbox{ ss case.}
\end{array}
\right.
\]
For $q=2$ in the ordinary smooth case and $q=1$ in the super-smooth case, if
$\pi(\{(F,\,\sigma)\in\mathcal F_n:\,\|p_{F\ast K,\sigma}-K\ast p_{0Y}\|_q\lesssim \psi_n\}\mid \Data,\,\hga)\rightarrow1$ in
$\P_{p_{0X}}^{(n)}$-probability or $\P_{p_{0X}}^{(\infty)}$-almost surely,
then also $\pi(\{(F,\,\sigma):\,\|p_Y-p_{0Y}\|_2\lesssim \epsilon_n\}\mid \Data,\,\hga)\rightarrow1$ in the same mode of convergence
and the proof is complete.
\quad$\square$
\subsection{Proofs for Aalen models}
For any intensity $\lambda$, we still denote $M_\lambda=\int_\Omega\lambda(t)\d t$ and $\bar\lambda=M_\lambda^{-1}\times\lambda\in{\mathcal F}_1.$
To prove Theorem \ref{cor:EB} and Theorem \ref{th:gene:aalen}, we need the following intermediate results that are based on classical tools subsequently defined. Recall that $\bar\epsilon_n = (\log n/n)^{1/3} $ and set $\epsilon_n = J_1 \bar \epsilon_n$ for some $J_1>0$ large enough.  The first result controls  the Kullback-Leibler divergence and moments of $\ell_n(\lambda_0)-\ell_n(\lambda)$, where $\ell_n(\lambda)$ is the log-likelihood evaluated at $\lambda$, whose expression is given in \eqref{loglik:aalen}.
\begin{Prop} \label{prop:KL:Aalen}
Let $v_n$ be a positive  sequence converging to 0  such that $nv_n^2 \rightarrow \infty$.
For $ H>0$, let  $$B_{k,n}(H,\, \lambda_0,\, v_n) = \{ \lambda:  \ \bar \lambda \in \bar B_{k,n}(H, \,\bar \lambda_0,\, v_n), \,\,\, |M_\lambda - M_{\lambda_0}| \leq v_n\}.$$
Then, for all  $k\geq 2$ and $\lambda \in B_{k,n}(H, \, \lambda_0,\, v_n) $, under \eqref{ass:Y1},
\begin{equation*}
\mathrm{KL}(\lambda_0;\,\lambda)\leq \kappa_0n v_n^2, \quad V_k(\lambda_0;\, \lambda)\leq \kappa( n v_n^2)^\frac{k}{2}
\end{equation*}
where $\kappa, \kappa_0$ depend only on $k$, $C_{1k}$, $H$, $\lambda_0$, $m_1$ and $m_2$; an expression of $\kappa_0$ is given in \eqref{kappa0}.
\end{Prop}
The second result establishes the existence of tests that control the numerator of posterior distributions.
\begin{Prop} \label{prop:test:Aalen}
Assume that conditions (i) and (ii) of Theorem \ref{th:gene:aalen} are satisfied. For any positive integer $j$, define
$$S_{n,j}(v_n) = \{ \lambda: \ \bar\lambda\in{\mathcal F}_n\mbox{ and } jv_n \leq \|\lambda -\lambda_0\|_1 \leq (j+1)v_n \}.$$
Then, under \eqref{ass:Y1},
 there exist $J_0 ,\, \rho,\, c>0$ such that, for all $j\geq J_0$, there exists a  test  $\phi_{n,j} \in [0,1]$  such that
$$ \E_{\lambda_0}^{(n)}[1_{\Gamma_n} \phi_{n,j}] \leq Ce^{ - c n j^2 v_n^2 }, \quad
\sup_{\lambda \in S_{n,j}(v_n) }\E_{\lambda}[1_{\Gamma_n} (1- \phi_{n,j}) ] \leq Ce^{ - c n j^2 v_n^2 }\quad\mbox{if }  j \leq \rho /v_n$$
$$ \E_{\lambda_0}^{(n)}[1_{\Gamma_n} \phi_{n,j}] \leq Ce^{ - c n j v_n }, \quad
\sup_{\lambda \in S_{n,j}(v_n) }\E_{\lambda}[1_{\Gamma_n} (1- \phi_{n,j}) ] \leq Ce^{ - c n j v_n }\quad\mbox{if }  j \geq \rho /v_n$$
for $C$ a constant.
\end{Prop}
\subsubsection{ Proof of Theorem \ref{cor:EB}}  \label{sec:pr:corEB}
\begin{proof}
Without loss of generality, we assume that $\Omega=[0,\,T]$. To apply Theorem \ref{th:gene:EB}, we must first define the transformation $\psi_{\gamma, \gamma'}$. Note that the parameter $\gamma$ only influences the prior on $\bar \lambda$ and has no impact on $M_\lambda$. As explained in Section \ref{sec:gene},  we can consider the following transformation:
for all $\gamma,\,\gamma' \in \R_+^*$, we set, for any $x$,
 $$ \bar \lambda (x) =\sum_{j=1}^\infty p_j \frac{\1_{(0,\,\theta_j)}(x)}{\theta_j} ,\quad \psi_{\gamma, \gamma'}(\bar \lambda)(x) =\sum_{j=1}^\infty p_j \frac{\1_{(0,\,G_{\gamma'}^{-1} (G_\gamma(\theta_j)))}(x)}{G_{\gamma'}^{-1}( G_\gamma(\theta_j))}, $$
 with
 $$p_j = V_j \prod_{l <j} (1-V_l), \quad V_j \sim {\rm Beta}( 1,\, A), \quad \theta_j \sim G_\gamma \quad \mbox{independently}.$$
So, if  $\bar \lambda$ is distributed according to a Dirichlet process mixture of uniform distributions with base measure indexed by $\gamma$, then
$\psi_{\gamma, \gamma'}(\bar \lambda)$ is distributed according to a Dirichlet process mixture of uniform distributions with base measure indexed by $\gamma'$. We prove Theorem  \ref{cor:EB} for both types of base measure introduced in \eqref{base:mono}.
 Let $G$ denote the cumulative distribution function of a $\Gamma(a,\,1)$ random variable and $g$ its density. Then, for the first type of base measure we have
$$G_{\gamma'}^{-1} (G_\gamma(\theta)) = \frac{ G^{-1} \left( G(\gamma \theta ) G(\gamma' T ) / G(\gamma T) \right) }{ \gamma'}\quad \mbox{ if $\theta \leq T$,} $$
and $G_{\gamma'}^{-1} (G_\gamma(\theta)) = T$ if $\theta \geq T$.   Therefore, for any $\theta \in [0,\, T]$, if $\gamma' \geq \gamma$ then \begin{equation}\label{ineq:psi}
G_{\gamma'}^{-1} (G_\gamma(\theta)) \leq \theta ,\quad
G_{\gamma'}^{-1} (G_\gamma(\theta)) \geq \frac{ \gamma \theta }{ \gamma' }.
\end{equation}
The second inequality of \eqref{ineq:psi} is straightforward. The first inequality of \eqref{ineq:psi} is equivalent to
$G(\gamma \theta ) G(\gamma' T )  \leq G(\gamma' \theta ) G(\gamma T )$ and  is deduced from the following argument. Let $$\Delta(\theta)=G(\gamma \theta ) G(\gamma' T ) - G(\gamma' \theta ) G(\gamma T ).$$ Then, $\Delta (0)=0$ and $\Delta ( T)=0$.
By Rolle's Theorem there  exists $c\in(0,\,T)$ such that $\Delta' ( c ) =0$.
We have
$$
\Delta'(\theta) = \gamma g(\gamma \theta ) G(\gamma' T ) - \gamma' g(\gamma' \theta ) G(\gamma T )$$ so is proportional to $\theta^{a-1}e^{-\gamma' \theta} ( \gamma^a e^{(\gamma'-\gamma )\theta} G(\gamma' T )- (\gamma')^a G(\gamma T ))$.
The function inside the brackets is increasing so that $\Delta'(\theta)\leq 0$ for $\theta\leq c$ and $\Delta'(\theta)\geq 0$ for $\theta\geq c$. Therefore, $\Delta $ is first decreasing and then increasing. Since $\Delta(0)=\Delta(T)=0$, $\Delta$ is negative on $(0,\,T)$, which achieves the proof of \eqref{ineq:psi}. For the second type of base measure, we have for $\theta\leq T$,
 $$ \forall\, \gamma,\, \gamma' >0,\quad  G_{\gamma'}^{-1}(G_\gamma ( \theta))  = \frac{ T \gamma }{\gamma'} \frac{ \theta }{ ( T - \theta + \theta\gamma/\gamma' )}$$ and  \eqref{ineq:psi} is straightforward.

We first verify assumption {[A1]}. At several places, by using \eqref{informally} and \eqref{ass:Y3}, we use that under $\P_\lambda^{(n)}(\cdot\mid \Gamma_n)$, for any interval $I$, the number of points of $N$ falling in $I$ is controlled by the number of points of a Poisson process with intensity $n(1+\alpha)m_2\lambda$ falling in $I$. Let $u_n = 1/(n\log n)$, so that $u_n=o(\bar\epsilon_n^2 )$ and choose $k \geq 6 $ so that $u_n^{-1}=  o((n\bar\epsilon_n^2)^{k/2})$ and \eqref{Nn} holds.  For $\kappa_0$ given in Proposition~\ref{prop:KL:Aalen}, we control $\P_{\lambda_0}^{(n)}\left( \ell_n(\lambda) - \ell_n( \lambda_0) \leq - (\kappa_0+1)n \bar \epsilon_n^2 \right) $. We follow most of the computations of \citet{salomond:13}. Let $e_n =  (n \bar\epsilon_n^2)^{-k/2}$ and set
$$\bar \lambda_{0n} (x) =  \frac{ \lambda_0(x)  \1_{ x \leq \theta_n} }{ \int_0^{\theta_n} \lambda_0(x) \d x }  \quad \mbox{ with }\quad \theta_n = \inf \left\{ \theta:\, \int_{0}^\theta \bar \lambda_0(x) \d x \geq 1 - \frac{e_n  }{n } \right\}$$ and $\lambda_{0n} = M_{\lambda_0} \bar \lambda_{0n}$.
Define the event $A_n = \{ \forall\, X\in N , \, X \leq \theta_n \}$. We shall need the following result. Given $\tilde N$  a Poisson process with intensity $n(1+\alpha)m_2\lambda_0$. If $\tilde N[0,\,T]=k$, we denote $\tilde N=\{X_1,\,\ldots,\,X_k\}$ and  conditionally on $\tilde N[0,\,T]=k$, $X_1,\,\ldots,\,X_k$ are i.i.d. with density $\bar\lambda_0$. So,
\begin{eqnarray*}
\P_{\lambda_0}^{(n)}(A_n^c\mid \Gamma_n)&\leq&\sum_{k=1}^{\infty}\P_{\lambda_0}^{(n)}(\exists\; X_i>\theta_n\mid \tilde N[0,\,T]=k)\P_{\lambda_0}^{(n)}(\tilde N[0,\,T]=k)\\
&\leq&\sum_{k=1}^{\infty}\left(1-\left(\int_0^{\theta_n}\bar\lambda_0(t)\d t\right)^k\right)\P_{\lambda_0}^{(n)}(\tilde N[0,\,T]=k)\\
&\leq&\sum_{k=1}^{\infty}\left(1-\left(1-\frac{e_n}{n}\right)^k\right)\P_{\lambda_0}^{(n)}(\tilde N[0,\,T]=k)\\
&=&O\left(\frac{e_n}{n}\E_{\lambda_0}^{(n)}[\tilde N[0,\,T]]\right)=O(e_n)= O((n \bar\epsilon_n^2)^{-k/2}).
\end{eqnarray*}
Now, we have
$$\P_{\lambda_0}^{(n)}\left( \ell_n(\lambda) - \ell_n( \lambda_0) \leq - (\kappa_0+2)n \bar \epsilon_n^2\mid \Gamma_n \right)\leq \P_{\lambda_0}^{(n)}\left( \ell_n(\lambda) - \ell_n( \lambda_0) \leq - (\kappa_0+2)n\bar  \epsilon_n^2\mid A_n,\,\Gamma_n \right)+ \P_{\lambda_0}^{(n)}(A_n^c\mid \Gamma_n).$$
We now deal with the first term.
On $\Gamma_n\cap A_n$,
\begin{equation*}
\begin{split}
\ell_n (\lambda_0) &= \ell_n( \lambda_{0n}) + \int_0^{\theta_n}\log\left(\frac{\lambda_0(t)}{\lambda_{0n}(t)}\right)\d N_t-\int_0^T(\lambda_0(t)-\lambda_{0n}(t))Y_t\d t
\\
&=\ell_n( \lambda_{0n}) +N[0,\,T] \log \left( \int_0^{\theta_n} \bar \lambda_0(t) \d t \right)  - M_{\lambda_0} \int_0^T \bar \lambda_0(t)Y_t  \d t + M_{\lambda_0}\frac{ \int_0^{\theta_n} \bar \lambda_0(t)Y_t  \d t}{ \int_0^{\theta_n} \bar \lambda_0(t) \d t }  \\
&\leq\ell_n( \lambda_{0n}) + M_{\lambda_0}\frac{ \int_{\theta_n}^T  \bar \lambda_0(t) \d t \int_0^{\theta_n} \bar \lambda_0(t)Y_t  \d t}{ \int_0^{\theta_n} \bar \lambda_0(t) \d t}
 -  M_{\lambda_0} \int_{\theta_n}^T \bar \lambda_0(t)Y_t  \d t \\
 &\leq \ell_n( \lambda_{0n}) + M_{\lambda_0} \frac{e_n (1+\alpha)m_2 }{ 1 - e_n/n }.
\end{split}
\end{equation*}
So, for $n$ large enough, for any $\lambda$,
\begin{eqnarray*}
\P_{\lambda_0}^{(n)}\left( \ell_n(\lambda) - \ell_n( \lambda_0) \leq - (\kappa_0+2)n \bar \epsilon_n^2\mid A_n,\,\Gamma_n \right)&\leq &
\P_{\lambda_0}^{(n)}\left( \ell_n(\lambda) - \ell_n( \lambda_{0n}) \leq - (\kappa_0+1)n \bar  \epsilon_n^2\mid A_n,\,\Gamma_n \right)\\
&=&\P_{\lambda_{0n}}^{(n)} \left( \ell_n(\lambda) - \ell_n( \lambda_{0n})  \leq - (\kappa_0+1)n \bar \epsilon_n^2\mid \Gamma_n \right)
\end{eqnarray*}
because $\P_{\lambda_0}^{(n)}(\cdot\mid A_n)=\P_{\lambda_{0n}}^{(n)}(\cdot).$
For all $\lambda = M_\lambda \bar \lambda \in B_{k,n}(H,\, \lambda_{0n},\,\bar\epsilon_n)$, using Proposition \ref{prop:KL:Aalen}, we obtain
\begin{equation} \label{likeratio:aalen}
\P_{\lambda_{0n}}^{(n)} \left( \ell_n(\lambda) - \ell_n(\lambda_{0n}) \leq -(\kappa_0+1)n \bar \epsilon_n^2\mid \Gamma_n\right) = O((n\bar\epsilon_n^2)^{-k/2}).
\end{equation}
To prove \eqref{KLcond}, we need to control $ \inf_{\gamma' \in [\gamma,\,\gamma+u_n]} \ell_n(M_\lambda\psi_{\gamma, \gamma'}(\bar \lambda) ) $. 
Using \eqref{ineq:psi}, we have for any $\gamma' \in [\gamma,\, \gamma+u_n]$, on $\Gamma_n$,
\begin{equation*}
\frac{ \gamma }{ \gamma+ u_n } \psi_{\gamma, \gamma+ u_n}(\bar \lambda)(t)  \leq \psi_{\gamma, \gamma'}(\bar\lambda)(t) \leq \psi_{\gamma, \gamma+ u_n}(\bar \lambda)(t)    + \sum_{j=1}^{\infty} p_j \frac{ \1_{\left(\frac{ \gamma }{ \gamma+ u_n }G_{\gamma+ u_n}^{-1}( G_\gamma( \theta_j) ),\, G_{\gamma+ u_n}^{-1}( G_\gamma( \theta_j) ) \right)} (t) }{ G_{\gamma+ u_n}^{-1}( G_\gamma( \theta_j) )}
\end{equation*}
so that for $n$ large enough,
\begin{equation*}
\begin{split}
 \inf_{\gamma' \in [\gamma,\,\gamma+u_n]} \ell_n(M_\lambda\psi_{\gamma, \gamma'}(\bar \lambda) )   & \geq
  -M_\lambda \int_0^T \psi_{\gamma, \gamma+ u_n}(\bar \lambda)(t)  Y_t \d t  -  \frac{M_\lambda u_n n(1+\alpha)m_2}{ \gamma+ u_n } \\
  & \quad \quad  + \int_0^T\log \left( M_\lambda \psi_{\gamma, \gamma+ u_n}(\bar \lambda)(t) \right) \d N_t + N[0,\,T] \log \left( \frac{ \gamma }{ \gamma+ u_n } \right) \\
  & \geq \ell_n(M_\lambda\psi_{\gamma, \gamma+u_n}(\bar \lambda) ) -  \frac{ u_n }{ \gamma+ u_n }\left( M_\lambda n(1+\alpha)m_2+  2N[0,\,T]\right)
  \end{split}
\end{equation*}
and, on the event $\{N[0,\,T] \leq 2 M_{\lambda_0} m_2 n \}$ which has probability going to 1,
 \begin{equation} \label{lb:like:aalen}
 \inf_{\gamma' \in [\gamma,\,\gamma+u_n]} \ell_n(M_\lambda\psi_{\gamma, \gamma'}(\bar \lambda) )
\geq \ell_n(M_\lambda\psi_{\gamma, \gamma+ u_n} (\bar \lambda)  )- n\gamma^{-1} (M_{\lambda}(1+\alpha)m_2 + 4m_2M_{\lambda_0})u_n.
 \end{equation}
Combining this lower bound with \eqref{likeratio:aalen}, we obtain that, for all $\lambda = M_\lambda \psi_{\gamma, \gamma+ u_n} (\bar \lambda)$, with $ \psi_{\gamma, \gamma+ u_n} (\bar \lambda)\in \bar B_{k,n}(H,\, \bar \lambda_{0n},\, \bar \epsilon_n)$ and $|M_\lambda-M_{\lambda_0}|\leq \bar \epsilon_n$, \begin{equation*}
\P_{\lambda_0}^{(n)}\left(  \inf_{\gamma' \in [\gamma,\,\gamma+u_n]} \ell_n(M_\lambda\psi_{\gamma, \gamma'}(\bar \lambda) ) - \ell_n( \lambda_0) \leq - (\kappa_0+2)n \bar \epsilon_n^2\mid A_n,\,\Gamma_n \right)= O((n\bar\epsilon_n^2)^{-k/2})
\end{equation*}
and assumption [A1] is satisfied if $J_1^2 \geq \kappa_0 + 2$.  We now verify assumption [A2]. First,
mimicking the proof of Lemma 8 of \citet{salomond:13}, we have that over any compact subset $\mathcal K'$ of $(0, \,\infty)$,
\begin{equation}\label{salomond:lem8}
\inf_{\gamma\in \mathcal K'}  \pi_1 \left( \bar B_{k,n}(H,\, \bar \lambda_{0n}, \,\bar \epsilon_n)  \mid \gamma \right) \geq e^{ - C_k n \bar \epsilon_n^2 }
 \end{equation}
for some $C_k >0$, when $n$ is large enough.
Let $B_n^\gamma = \{\lambda : \, \psi_{\gamma, \gamma+u_n}(\bar\lambda) \in \bar B_{k,n}(H,\, \bar \lambda_{0n},\, \bar \epsilon_n),\, M_\lambda\in [M_{\lambda_0}  - \bar \epsilon_n,\, M_{\lambda_0} +\bar \epsilon_n]\}$. Then
$$\pi(B_n^\gamma\mid \gamma)  = \pi_1(\bar B_{k,n}(H, \,\bar \lambda_{0n},\, \bar \epsilon_n)\mid \gamma+u_n) \pi_M([M_{\lambda_0} - \bar\epsilon_n,\, M_{\lambda_0} + \bar \epsilon_n])$$
which, together with \eqref{salomond:lem8},  implies that
\begin{equation*}\label{KL:prior:aalen}
\inf_{\gamma \in \mathcal K} \pi(B_n^\gamma\mid \gamma) \geq e^{ - 2C_k n\bar \epsilon_n^2 }
\end{equation*}
when $n$ is large enough, so that \eqref{ratio:mass} is satisfied as soon as $J_1$ is large enough.  We now control the measure $Q_{\gamma , n}^\lambda$, with
$$\d Q_{\gamma, n}^\lambda=1_{\Gamma_n}\times \sup_{\| \gamma'  - \gamma\| \leq u_n} \exp(\ell_n(M_\lambda\psi_{\gamma, \gamma'}(\bar\lambda)))\d\mu$$
and $\mu$ is the measure such that under $\mu$ the process is an homogeneous Poisson process with intensity $1$.
Using \eqref{ineq:psi} and similarly to \eqref{lb:like:aalen}, we obtain that, for all $\gamma' \in [ \gamma,\, \gamma + u_n]$,
 $$ \bar \lambda (t)   - \sum_{j=1}^{\infty } p_j \frac{ \1_{\left( \frac{\gamma \theta_j}{ \gamma'},\, \theta_j\right)}(t) }{ \theta_j} \leq \psi_{\gamma, \gamma'}(\bar\lambda)(t) \leq \frac{ \gamma+ u_n}{ \gamma } \bar \lambda(t) $$
 and we have
 \begin{equation*}
 \begin{split}
Q_{\gamma, n}^\lambda( \Xn) &= \E_{\mu}^{(n)} \left[ 1_{\Gamma_n}\sup_{\gamma' \in [\gamma,\, \gamma+u_n]}\exp\left( T - M_\lambda\int_0^T \psi_{\gamma, \gamma'}(\bar \lambda) (t) Y_t \d t + \int_0^T \log \left( M_\lambda \psi_{\gamma, \gamma'}( \bar \lambda)(t)\right) \d N_t \right) \right] \\
&\leq \E_{\lambda}^{(n)} \left[1_{\Gamma_n} \exp( nm_2(1+\alpha)M_\lambda u_n/(\gamma+u_n)  + \log(1+u_n/\gamma) N[0,\,T]) \right] \\
&\leq \E_{\lambda}^{(n)} \left[1_{\Gamma_n} \exp( nm_2(1+\alpha)M_\lambda \gamma^{-1}u_n + u_n\gamma^{-1} N[0,\,T]) \right] \\
&\leq  \exp\left( nm_2(1+\alpha)M_\lambda \gamma^{-1}u_n + (1+ \alpha)nm_2M_\lambda(e^{u_n/\gamma}-1)\right)\\
&\leq  \exp\left( 3nm_2(1 + \alpha)\gamma^{-1}M_\lambda u_n\right)
\end{split}
 \end{equation*}
 when $n$ is large enough. Let $\phi_{n,j}$ be the tests defined in Proposition \ref{prop:test:Aalen} over $S_{n,j}(\bar\epsilon_n)$.  Using the previous computations, we have
 \begin{equation*}
 \begin{split}
 Q_{\gamma, n}^\lambda[1 - \phi_{n,j} ] &\leq  \E_{\lambda} ^{(n)}\left[(1 - \phi_{n,j})\exp( nm_2(1+\alpha)M_\lambda \gamma^{-1}u_n + u_n\gamma^{-1} N[0,\,T]) 1_{\Gamma_n}\right] \\
 &\leq e^{ nm_2(1+\alpha)M_\lambda \gamma^{-1}u_n} \left(\E_{\lambda} ^{(n)}\left[ (1 -  \phi_{n,j})1_{\Gamma_n}\right] \E_{\lambda} ^{(n)} \left[ e^{2\gamma^{-1}u_nN[0,\,T]}1_{\Gamma_n}\right]\right)^{1/2} \\
 &\leq e^{4nm_2(1 + \alpha)\gamma^{-1}M_\lambda u_n} \max \{ e^{ - cnj^2\bar \epsilon_n^2 /2},\, e^{ - cnj \bar\epsilon_n /2}\}.
 \end{split}
 \end{equation*}
As in  \citet{salomond:13}, we set $\mathcal F_n = \{ \bar \lambda: \ \bar \lambda (0 ) \leq M_n\}$ with $M_n = \exp \left( c_1 n \epsilon_n^2 \right)$ and $c_1$ is a positive constant. From Lemma 9 of \citet{salomond:13}, there exists $a>0$ such that
$\sup_{\gamma\in \mathcal K'}  \pi_1(\mathcal F_n^c|\gamma) \leq e^{ - c_1 (a+1) n\epsilon_n^2 } $, so that 	when $n$ is large enough,
\begin{equation*}
\begin{split}
 \sup_{\gamma\in \mathcal K'} \int_{\R^+} \int_{\mathcal F_n^c} Q_{\gamma, n}^\lambda( \Xn)d\pi_1( \bar \lambda | \gamma) \pi_M(M_\lambda) dM_\lambda
& \lesssim e^{ - c_1 (a+1) n\epsilon_n^2 } \int_{\R^+}  \exp\left( \delta M_\lambda \right) \pi_M(M_\lambda) dM_\lambda \\
&\lesssim e^{ - c_1 (a+1) n\epsilon_n^2 },
\end{split}
\end{equation*}
with $\delta$ that can be chosen as small as needed since $nu_n=o(1)$. This proves  \eqref{Thetanc} by choosing $c_1$ conveniently.
Combining the above upper bound with Proposition \ref{prop:test:Aalen}, together with Remark 1, achieves the proof of Theorem \ref{cor:EB}.
\end{proof}
\subsubsection{Proof of Theorem \ref{th:gene:aalen}}
The proof of Theorem  \ref{th:gene:aalen} is similar to the proof of Theorem 1 of \citet{ghosal:vdv:07} once Proposition~\ref{prop:KL:Aalen} and Proposition \ref{prop:test:Aalen} are proved.
Let $U_n = \{\lambda: \ \| \lambda - \lambda_0\|_1 \geq  J_1v_n \}$
and recall that
\begin{equation*}
\pi(U_n\mid N) = \frac{ \int_{U_n} e^{\ell_n(\lambda) - \ell_n(\lambda_0) }\d\pi(\lambda) }{ \int_{\mathcal F} e^{\ell_n(\lambda) - \ell_n(\lambda_0) }\d\pi(\lambda) },
\end{equation*}
where $\ell_n(\lambda)$ is the log-likelihood evaluated at $\lambda$. We use notations of Proposition~\ref{prop:KL:Aalen} and Proposition~\ref{prop:test:Aalen}. Writing
 $$ D_n = \int_{\mathcal F} e^{\ell_n(\lambda) - \ell_n(\lambda_0) }\d\pi(\lambda),$$
 we have
\begin{equation*}
 \begin{split}
 &
 \P_{\lambda_0}^{(n)}\left( D_n \leq e^{ - (\kappa_0+1) n v_n^2 }\pi_1( \bar B_{k,n}(H, \bar \lambda_0,  v_n) ) \right) \\
& \leq \P_{\lambda_0}^{(n)}\left(\int_{B_{k,n}(H,  \lambda_0, v_n)}\frac{ \exp[\ell_n(\lambda) - \ell_n(\lambda_0)]d\pi(\lambda)}{\pi(B_{k,n}(H, \lambda_0, v_n) ) } \leq - (\kappa_0+1) n v_n^2+\log\left(\frac{\pi_1( \bar B_{k,n}(H, \bar \lambda_0, v_n))}{\pi(B_{k,n}(H, \lambda_0,  v_n))}\right)\right).
\end{split}
\end{equation*}
 Note that since  $v_n^2 \geq \log n/n$,
 $$ \pi(B_{k,n}(H,  \lambda_0, v_n))\gtrsim \pi_1( \bar B_{k,n}(H, \bar \lambda_0,  v_n))v_n\gtrsim \pi_1( \bar B_{k,n}(H, \bar \lambda_0, v_n) )e^{-nv_n^2/2}  $$ so that
using Proposition \ref{prop:KL:Aalen} and the Markov inequality, we obtain that
 $$  \P_{\lambda_0}^{(n)}\left( D_n \leq e^{ - (\kappa_0 +1) n  v_n^2 }\pi_1( \bar B_{k,n}(H, \bar \lambda_0, v_n) ) \right) \lesssim (nv_n^2)^{-k/2}.$$
Moreover,  \eqref{norm:mino} implies that $\pi(S_{n,j}(v_n))\leq \pi_1(\bar S_{n,j})$
 and  using the tests $\phi_{n,j}$ of Proposition \ref{prop:test:Aalen}, we have for  $J_1 \geq J_0$, mimicking the proof of Theorem~1 of \cite{ghosal:vdv:07},
 \begin{equation*}
 \begin{split}
   &\E_{\lambda_0}^{(n)} \left[ 1_{\Gamma_n} \pi\left( \lambda : \ \|\lambda -\lambda_0\|_1\geq J_1 v_n |N\right) \right]
 \leq \sum_{j\geq J_1 }  \E_{\lambda_0}^{(n)} \left[1_{\Gamma_n}\phi_{n,j} \right]  + \sum_{j = J_1}^{\lfloor \rho/v_n\rfloor }  e^{  (\kappa_0 +1) n v_n^2 } \frac{ \pi_1(\bar S_{n,j}) e^{ - c n j^2 v_n^2 } }{ \pi_1( \bar B_{k,n}(H, \bar \lambda_0, v_n) )  } \\
 & \quad + \sum_{j > \rho/v_n  }   \frac{ e^{  (\kappa_0 +1) n v_n^2 }\pi_1(\bar S_{n,j} ) e^{ - c n j v_n } }{ \pi_1( \bar B_{k,n}(H, \bar \lambda_0, v_n) )  }  + \frac{e^{  (\kappa_0 +1) n  v_n^2 }  \pi_1(\mathcal F_n^c )  }{ \pi_1( \bar B_{k,n}(H, \bar \lambda_0,  v_n) )  } + \P_{\lambda_0}^{(n)}\left( D_n \leq e^{ - (\kappa_0 +1) n v_n^2 }\pi_1( \bar B_{k,n}(H, \bar \lambda_0, v_n) ) \right)  \\
     & \lesssim (nv_n^2)^{-k/2},
 \end{split}
 \end{equation*}
which proves the result since $\P_{\lambda_0}^{(n)}(\Gamma_n^c)=o(1)$.

\appendix
\section{Additional proofs}\label{app:Aalen}
We use in the sequel that for any densities $f$ and $g$,
$\|f-g\|_1\leq 2h(f,\,g).$
\subsection{Proof of Proposition \ref{prop:KL:Aalen}}
We recall that the log-likelihood evaluated at $\lambda$ is given by
$$\ell_n(\lambda)=\int_0^T\log(\lambda(t))\d N_t-\int_0^T\lambda(t)Y_t\d t,$$ see \cite{ABGK}.  Since on $\Omega^c$, $N$ is empty and $Y_t\equiv0$ almost surely, in the sequel we assume, without loss of generality, that $\Omega=[0,\,T]$.
 We denote by $$M_n(\lambda) = \int_0^T \lambda(x) \mu_n(x) \d x, \quad M_n(\lambda_0) = \int_0^T \lambda_0(x) \mu_n(x) \d x,$$
$$\bar \lambda_n(x) = \frac{\lambda(x) \mu_n(x) }{M_n(\lambda)}= \frac{ \bar \lambda(x) \tilde\mu_n(x) }{ \int_0^T \bar \lambda(t) \tilde\mu_n(t)\d t } \quad \mbox{and } \quad  \bar \lambda_{0,n}(x) = \frac{ \lambda_0(x) \mu_n(x) }{M_n(\lambda_0) } = \frac{ \bar \lambda_0(x)\tilde\mu_n(x) }{ \int_0^T \bar \lambda_0(t)\tilde \mu_n(t)\d t}.$$
Then, by using straightforward computations,
\begin{equation}\label{compil1}
\begin{split}
\textrm{KL}(\lambda_0;\,\lambda)&=\E_{\lambda_0}^{(n)}[\ell_n(\lambda_0)-\ell_n(\lambda)]=M_{n}(\lambda_0)\left(\textrm{KL}
(\bar{\lambda}_{0,n};\,\bar{\lambda}_{n})+\frac{M_n(\lambda)}{M_{n}(\lambda_0)} -1-\log\left(\frac{M_n(\lambda)}{M_{n}(\lambda_0)}\right)  \right)\\
&=M_{n}(\lambda_0)\pq{\textrm{KL}(\bar{\lambda}_{0,n};\,\bar{\lambda}_{n})+\phi\left(\frac{M_n(\lambda)}{M_{n}(\lambda_0)}\right)}\\
&\leq nm_2M_{\lambda_0}\pq{\textrm{KL}(\bar{\lambda}_{0,n};\,\bar{\lambda}_{n})+\phi\left(\frac{M_n(\lambda)}{M_{n}(\lambda_0)}\right)},
\end{split}
\end{equation}
where $\phi(x) = x- 1 - \log x$ and
$$\textrm{KL}(\bar{\lambda}_{0,n};\,\bar{\lambda}_{n})=
\int_0^T\log\left(\frac{\bar\lambda_{0,n}(t)}{\bar\lambda_{n}(t)}\right)\bar\lambda_{0,n}(t)\d t.$$
Now, we control $\textrm{KL}(\bar{\lambda}_{0,n};\,\bar{\lambda}_{n})$ for $\lambda \in B_{k,n}(H, \lambda_0,v_n)$. By using Lemma~8.2 of \cite{ghosal:ghosh:vdv:00}, we have
\begin{eqnarray}\label{compil2}
\textrm{KL}(\bar{\lambda}_{0,n};\,\bar{\lambda}_{n})&\leq& 2h^2(\bar{\lambda}_{0,n},\,\bar{\lambda}_{n})\left(1+\log\left\|\frac{\bar{\lambda}_{0,n}}{\bar{\lambda}_{n}}\right\|_\infty\right)\nonumber\\
&\leq&2h^2(\bar{\lambda}_{0,n},\,\bar{\lambda}_{n})\left(1+\log\left(\frac{m_2}{m_1}\right)+\log\left\|\frac{\bar{\lambda}_0}{\bar{\lambda}}\right\|_\infty\right)\nonumber\\
&\leq&2\left(1+\log\left(\frac{m_2}{m_1}\right)\right)h^2(\bar{\lambda}_{0,n},\,\bar{\lambda}_{n})\left(1+\log\left\|\frac{\bar{\lambda}_0}{\bar{\lambda}}\right\|_\infty\right)
\end{eqnarray}
since $1+\log(m_2/m_1)\geq 1$. We now deal with $h^2(\bar{\lambda}_{0,n},\,\bar{\lambda}_{n})$. We have
\begin{eqnarray*}
h^2(\bar{\lambda}_{0,n},\,\bar{\lambda}_{n})&=&\int\left(\sqrt{\bar \lambda_{0,n}(x)}-\sqrt{ \bar \lambda_{n}(x)}\right)^2\d x
=\int\left(\sqrt{\frac{\bar \lambda_0(x)\tilde\mu_n(x)}{\int\bar\lambda_0(t)\tilde\mu_n(t)\d t}}-\sqrt{\frac{\bar \lambda(x)\tilde\mu_n(x)}{\int\bar\lambda(t)\tilde\mu_n(t)\d t}}\right)^2\d x\\
&\leq&2m_2\int\left(\sqrt{\frac{\bar \lambda_0(x)}{\int\bar\lambda_0(t)\tilde\mu_n(t)\d t}}-\sqrt{\frac{\bar \lambda_0(x)}{\int\bar\lambda(t)\tilde\mu_n(t)\d t}}\right)^2dx\\&&+2m_2\int\left(\sqrt{\frac{\bar \lambda_0(x)}{\int\bar\lambda(t)\tilde\mu_n(t)\d t}}-\sqrt{\frac{\bar \lambda(x)}{\int\bar\lambda(t)\tilde\mu_n(t)\d t}}\right)^2\d x\\
&\leq&2m_2U_n+\frac{2m_2}{m_1}h^2(\bar\lambda_0,\,\bar\lambda),
\end{eqnarray*}
with $$U_n=\left(\sqrt{\frac{1}{\int\bar\lambda_0(x)\tilde\mu_n(x)\d x}}-\sqrt{\frac{1}{\int\bar\lambda(x)\tilde\mu_n(x)\d x}}\right)^2.$$
We denote by
$$\tilde\epsilon_n:=\frac{1}{\int\bar\lambda_0(x)\tilde\mu_n(x)\d x}\int(\bar\lambda(x)-\bar\lambda_0(x))\tilde\mu_n(x)\d x,$$
so that
$$|\tilde\epsilon_n|\leq\frac{1}{m_1}\int|\bar\lambda(x)-\bar\lambda_0(x)|\tilde\mu_n(x)\d x\leq\frac{2m_2}{m_1}h(\bar\lambda_0,\,\bar\lambda).$$
Then,
$$U_n=\frac{1}{\int\bar\lambda_0(x)\tilde\mu_n(x)\d x}\left(1-\frac{1}{\sqrt{1+\tilde\epsilon_n}}\right)^2\leq \frac{\tilde\epsilon_n^2}{4m_1}\leq \frac{m_2^2}{m_1^3}h^2(\bar\lambda_0,\,\bar\lambda).$$
Finally,
\begin{equation}\label{compil3}
h^2(\bar{\lambda}_{0,n},\,\bar{\lambda}_{n})\leq \frac{2m_2}{m_1}\left(1+\frac{m_2^2}{m_1^2}\right)h^2(\bar\lambda_0,\,\bar\lambda).
\end{equation}
It remains to bound $\phi\left(\frac{M_n(\lambda)}{M_{n}(\lambda_0)}\right)$. We have
\begin{equation*}
\begin{split}
|M_n(\lambda_0)-M_n(\lambda)|&\leq\int_0^T|\lambda(t)-\lambda_0(t)|\mu_n(t)\d t  \leq nm_2\int_0^T|\lambda(t)-\lambda_0(t)|\d t\\
&\leq\frac{m_2}{m_1M_{\lambda_0}}M_n(\lambda_0)\left(M_{\lambda_0}\|\bar\lambda-\bar\lambda_0\|_1+|M_\lambda-M_{\lambda_0}|\right)\\
&\leq\frac{m_2M_n(\lambda_0)}{m_1M_{\lambda_0}}\left(2M_{\lambda_0}h(\bar\lambda,\bar\lambda_0)+|M_\lambda-M_{\lambda_0}|\right)
\leq\frac{m_2(2M_{\lambda_0}+1)M_n(\lambda_0)v_n}{m_1M_{\lambda_0}}.
\end{split}
\end{equation*}
Since $\phi(u+1)\leq u^2$ if $|u|\leq\frac{1}{2}$, we have for $n$ large enough,
\begin{equation}\label{compil4}
\phi\left(\frac{M_n(\lambda)}{M_{n}(\lambda_0)}\right)\leq \frac{m_2^2(2M_{\lambda_0}+1)^2}{m_1^2M_{\lambda_0}^2}v_n^2.
\end{equation}
Combining \eqref{compil1}, \eqref{compil2}, \eqref{compil3} and \eqref{compil4}, we have $\textrm{KL}(\lambda_0;\,\lambda)\leq \kappa_0 nv_n^2$  for $n$ large enough,
with
\begin{equation}\label{kappa0}
\kappa_0=m_2^2M_{\lambda_0}\left(\frac{4}{m_1}\left(1+\log\left(\frac{m_2}{m1}\right)\right)\left(1+\frac{m_2^2}{m_1^2}\right)+\frac{m_2(2M_{\lambda_0}+1)^2}{m_1^2M_{\lambda_0}^2}\right).
\end{equation}
We now deal with
$$ V_{2k}(\lambda_0;\,\lambda) = \E_{\lambda_0}^{(n)}[|\ell_n(\lambda_0)-\ell_n(\lambda)-\E_{\lambda_0}^{(n)}[\ell_n(\lambda_0)-\ell_n(\lambda)]|^{2k}],$$
with $k>1$. In the sequel, we denote by $C$ a constant that may change from line to line. Straightforward computations lead to
\begin{eqnarray*}
V_{2k}(\lambda_0;\, \lambda)&=&\E_{\lambda_0}^{(n)}\left[\left|-\int_0^T\left(\lambda_0(t)-\lambda(t)-\lambda_0(t)\log\left(\frac{\lambda_0(t)}{\lambda(t)}\right)\right)( Y_t-\mu_n(t))\d t \right. \right.\\
& &\hspace{5cm} \left. \left.+\int_0^T\log\left(\frac{\lambda_0(t)}{\lambda(t)}\right)(\d N_t-Y_t\lambda_0(t)\d t)\right|^{2k}\right]\\
&\leq&2^{2k-1}(A_{2k}+B_{2k}),
\end{eqnarray*}
with
$$B_{2k}:=\E_{\lambda_0}^{(n)}\left[ \left|\int_0^T\log\left(\frac{\lambda_0(t)}{\lambda(t)}\right)(\d N_t-Y_t\lambda_0(t)\d t)\right|^{2k}\right]$$
and, by using \eqref{moment},
\begin{eqnarray*}
A_{2k}&:=&\E_{\lambda_0}^{(n)}\left[\left|\int_0^T\left(\lambda_0(t)-\lambda(t)-\lambda_0(t)\log\left(\frac{\lambda_0(t)}{\lambda(t)}\right)\right)( Y_t-\mu_n(t))\d t \right|^{2k} \right]\\
&\leq&\left(\int_0^T \left(\lambda_0(t)-\lambda(t)-\lambda_0(t)\log\left(\frac{\lambda_0(t)}{\lambda(t)}\right)\right)^2 \d t \right)^k\times
\E_{\lambda_0}^{(n)}\left[ \left(  \int_0^T ( Y_t-\mu_n(t))^2 \d t\right)^k  \right] \\
&\leq& 2^{2k-1}C_{1k}n^k\left(A_{2k,1}+A_{2k,2}\right),
\end{eqnarray*}
with, for $\lambda\in B_{k,n}(H,\, \lambda_0,v_n)$,
\begin{eqnarray*}
A_{2k,1}&:=&\left(\int_0^T \lambda_0^2(t)\log^2\left(\frac{\lambda_0(t)}{\lambda(t)}\right) \d t \right)^k
\leq M_{\lambda_0}^{2k}\|\bar\lambda_0\|_\infty^k\left(\int_0^T\bar \lambda_0(t)\log^2\left(\frac{M_{\lambda_0}\bar\lambda_0(t)}{M_{\lambda}\bar\lambda(t)}\right) \d t\right)^k\\
&\leq&2^{2k-1} M_{\lambda_0}^{2k}\|\bar\lambda_0\|_\infty^k\left(E_2^k(\bar\lambda_0,\,\bar\lambda)+\left|\log\left(\frac{M_\lambda}{M_{\lambda_0}}\right)\right|^{2k}\right)\\
&\leq&C\left(E_2^k(\bar\lambda_0,\,\bar\lambda)+\left|M_\lambda-M_{\lambda_0}\right|^{2k}\right)
\leq Cv_n^{2k}
\end{eqnarray*}
and
\begin{eqnarray*}
A_{2k,2}&:=&\left(\int_0^T \left(\lambda_0(t)-\lambda(t)\right)^2 \d t \right)^k\\
&=&\left(\int_0^T \left((M_{\lambda_0}-M_\lambda)\bar\lambda_0(t)-M_\lambda(\bar\lambda(t)-\bar\lambda_0(t))\right)^2 \d t \right)^k\\
&\leq&2^{2k-1}\left(\|\bar\lambda_0\|_\infty^k(M_{\lambda_0}-M_\lambda)^{2k}+M_{\lambda}^{2k}\left(\int_0^T \left(\sqrt{\bar\lambda_0(t)}-\sqrt{\bar\lambda(t)}\right)^2 \left(\sqrt{\bar\lambda_0(t)}+\sqrt{\bar\lambda(t)}\right)^2 \d t \right)^k\right)\\
&\leq&2^{2k-1}\left(\|\bar\lambda_0\|_\infty^k(M_{\lambda_0}-M_\lambda)^{2k}+2^kM_{\lambda}^{2k}(\|\bar\lambda_0\|_\infty+\|\bar\lambda\|_\infty)^kh^{2k}(\bar\lambda_0,\,\bar\lambda)\right)
\leq Cv_n^{2k}.
\end{eqnarray*}
Therefore,
$$A_{2k}\leq C(nv_n^2)^k.$$
To deal with $B_{2k}$, we set for any $T>0$,
$$M_T:=\int_0^T\log\left(\frac{\lambda_0(t)}{\lambda(t)}\right)(\d N_t-Y_t\lambda_0(t)\d t),$$
so $(M_T)_T$ is a martingale. Using the Burkholder-Davis-Gundy Inequality (see Theorem B.15 in \cite{Karr}), there exists a constant $C(k)$ only depending on $k$ such that, since $2k> 1$,
$$\E_{\lambda_0}^{(n)}[|M_T|^{2k}]\leq C(k)\E_{\lambda_0}^{(n)}\left[ \left|\int_0^T\log^2\left(\frac{\lambda_0(t)}{\lambda(t)}\right)\d N_t\right|^k\right].$$
Therefore, for $k>1$,
\begin{eqnarray*}
B_{2k}&=&\E_{\lambda_0}^{(n)}[|M_T|^{2k}]\\
&\leq&3^{k-1}C(k)\left(\E_{\lambda_0}^{(n)}\left[ \left|\int_0^T\log^2\left(\frac{\lambda_0(t)}{\lambda(t)}\right)(\d N_t-Y_t\lambda_0(t)\d t)\right|^k+\left|\int_0^T\log^2\left(\frac{\lambda_0(t)}{\lambda(t)}\right)(Y_t-\mu_n(t))\lambda_0(t)\d t\right|^k\right.\right.\\
&&\hspace{3cm}+\left.\left.\left|\int_0^T\log^2\left(\frac{\lambda_0(t)}{\lambda(t)}\right)\mu_n(t)\lambda_0(t)\d t\right|^k\right]\right)\\
&\leq&3^{k-1}C(k)(B_{k,2}^{(0)}+B_{k,2}^{(1)}+B_{k,2}^{(2)}),
\end{eqnarray*}
with
\begin{eqnarray*}
B_{k,2}^{(0)}&=&\E_{\lambda_0}^{(n)}\left[\left|\int_0^T\log^2\left(\frac{\lambda_0(t)}{\lambda(t)}\right)(\d N_t-Y_t\lambda_0(t)\d t)\right|^k\right],\\
B_{k,2}^{(1)}&=&\E_{\lambda_0}^{(n)}\left[\left|\int_0^T\log^2\left(\frac{\lambda_0(t)}{\lambda(t)}\right)(Y_t-\mu_n(t))\lambda_0(t)\d t\right|^k\right]\\
B_{k,2}^{(2)}&=&\left|\int_0^T\log^2\left(\frac{\lambda_0(t)}{\lambda(t)}\right)\mu_n(t)\lambda_0(t)\d t\right|^k.
\end{eqnarray*}
This can be iterated: we set $J=\min\{j\in\mathbb{N}: \ 2^j\geq k\}$ so that  $1<k2^{1-J}\leq 2$. There exists a constant $C_k$, only depending on $k$, such that for
$$B_{k2^{1-j},2^j}^{(1)}=\E_{\lambda_0}^{(n)}\left[\left|\int_0^T\log^{2^j}\left(\frac{\lambda_0(t)}{\lambda(t)}\right)(Y_t-\mu_n(t))\lambda_0(t)\d t\right|^{k2^{1-j}}\right]$$
and
$$B_{k2^{1-j},2^j}^{(2)}=\left|\int_0^T\log^{2^j}\left(\frac{\lambda_0(t)}{\lambda(t)}\right)\mu_n(t)\lambda_0(t)\d t\right|^{k2^{1-j}},$$
\begin{eqnarray*}
B_{2k}&\leq& C_k\left(\E_{\lambda_0}^{(n)}\left[\left|\int_0^T\log^{2^J}\left(\frac{\lambda_0(t)}{\lambda(t)}\right)(\d N_t-Y_t\lambda_0(t) \d t)\right|^{k2^{1-J}}\right]+\sum_{j=1}^J(B_{k2^{1-j},2^j}^{(1)}+B_{k2^{1-j},2^j}^{(2)}) \right)\\
&\leq&C_k\left(\left(\E_{\lambda_0}^{(n)}\left[\left|\int_0^T\log^{2^J}\left(\frac{\lambda_0(t)}{\lambda(t)}\right)(\d N_t-Y_t\lambda_0(t)\d t)\right|^2\right]\right)^{k2^{-J}}+\sum_{j=1}^J(B_{k2^{1-j},2^j}^{(1)}+B_{k2^{1-j},2^j}^{(2)}) \right)\\
&=&C_k\left(\left(\E_{\lambda_0}^{(n)}\left[\int_0^T\log^{2^{J+1}}\left(\frac{\lambda_0(t)}{\lambda(t)}\right)Y_t\lambda_0(t)\d t\right]\right)^{k2^{-J}}+\sum_{j=1}^J(B_{k2^{1-j},2^j}^{(1)}+B_{k2^{1-j},2^j}^{(2)}) \right)\\
&=&C_k\left(\left(\int_0^T\log^{2^{J+1}}\left(\frac{\lambda_0(t)}{\lambda(t)}\right)\mu_n(t)\lambda_0(t)\d t\right)^{k2^{-J}}+\sum_{j=1}^J(B_{k2^{1-j},2^j}^{(1)}+B_{k2^{1-j},2^j}^{(2)}) \right)\\
&=&C_k\left(B_{k2^{-J},2^{J+1}}^{(2)}+\sum_{j=1}^J(B_{k2^{1-j},2^j}^{(1)}+B_{k2^{1-j},2^j}^{(2)}) \right).
\end{eqnarray*}
Note that for any $1\leq j\leq J$,
\begin{eqnarray*}
B_{k2^{1-j},2^j}^{(1)}&\leq&\left(\int_0^T\log^{2^{j+1}}\left(\frac{\lambda_0(t)}{\lambda(t)}\right)\lambda_0^2(t)\d t\right)^{k2^{-j}}\times \E_{\lambda_0}^{(n)}\left[\left(\int_0^T(Y_t-\mu_n(t))^2\d t\right)^{k2^{-j}}\right]\\
&\leq&C(M_{\lambda_0}^2\|\bar\lambda_0\|_\infty)^{k2^{-j}}\left(\int_0^T\log^{2^{j+1}}\left(\frac{M_{\lambda_0}\bar\lambda_0(t)}{M_{\lambda}\bar\lambda(t)}\right)\bar\lambda_0(t)\d t\right)^{k2^{-j}}\times n^{k2^{-j}}\\
&\leq& C\left(\log^{2^{j+1}}\left(\frac{M_{\lambda_0}}{M_\lambda}\right)+E_{2^{j+1}}(\bar\lambda_0,\,\bar\lambda)\right)^{k2^{-j}}\times n^{k2^{-j}}\\
&\leq&C(nv_n^2)^{k2^{-j}}\leq C(nv_n^2)^k,
\end{eqnarray*}
where we have used \eqref{moment}. Similarly,  for any $j\geq 1$,
\begin{eqnarray*}
B_{k2^{1-j},2^j}^{(2)}&\leq& (nm_2M_{\lambda_0})^{k2^{1-j}}
\left(\int_0^T\log^{2^j}\left(\frac{M_{\lambda_0}\bar\lambda_0(t)}{M_{\lambda}\bar\lambda(t)}\right)\bar\lambda_0(t)\d t\right)^{k2^{1-j}}\\
&\leq&C\left(\log^{2^j}\left(\frac{M_{\lambda_0}}{M_\lambda}\right)+E_{2^j}(\bar\lambda_0,\,\bar\lambda)\right)^{k2^{1-j}} \times n^{k2^{1-j}}\\
&\leq&C(nv_n^2)^{k2^{1-j}} \leq C(nv_n^2)^k.
\end{eqnarray*}
Therefore, for any $k>1$,
$$V_{2k}(\lambda_0;\, \lambda)\leq \kappa(nv_n^2)^k,$$
where $\kappa$ depends on $C_{1k}$, $k$, $H$, $\lambda_0$, $m_1$ and $m_2$. Using previous computations, the case $k=1$ is straightforward.
So, we obtain the result for $V_k(\lambda_0;\, \lambda)$ with $k\geq 2$.
\subsection{Proof of Proposition \ref{prop:test:Aalen}}
We shall use the following lemma whose proof is given in Section~\ref{sec:proof1}.
 \begin{Lemma} \label{lem:tests}
There exist constants $\xi,\,K>0$, only depending on $M_{\lambda_0}$, $\alpha,$ $m_1$ and $m_2$, such that, for any $\lambda_1$,
there exists a test $\phi_{\lambda_1}$ so that
$$\E_{\lambda_0}^{(n)}[1_{\Gamma_n}\phi_{\lambda_1}]\leq 2\exp\left(- K n\| \lambda_1 - \lambda_0 \|_{1} \times\min\{\| \lambda_1 - \lambda_0 \|_1,\,m_1\} \right)$$
and
$$\sup_{\lambda: \ \| \lambda - \lambda_1 \|_{1} <\xi\| \lambda_1 - \lambda_0 \|_1 } \E_{\lambda}[1_{\Gamma_n}(1-\phi_{\lambda_1})]\leq 2\exp\left(-Kn\| \lambda_1 - \lambda_0 \|_1\times\min\{\| \lambda_1 - \lambda_0 \|_{1},\,m_1\} \right).$$
\end{Lemma}

\medskip

We consider the setting of Lemma~\ref{lem:tests} and a covering of $S_{n,j}(v_n)$ with $\L_1$-balls with radius $\xi j v_n$ and centers $(\lambda_{l,j})_{l=1,\,\ldots,\, D_j }$, where $D_j$ is the covering number of $S_{n,j}(v_n)$ by such balls. We then set $\phi_{n,j} = \max_{l=1,\,\ldots,\, D_j}\phi_{\lambda_{l,j}}$, where the $\phi_{\lambda_{l,j}}$'s are defined in Lemma~\ref{lem:tests}. So, there exists a constant $\rho>0$ such that
 $$ \E_{\lambda_0}^{(n)}[ 1_{\Gamma_n} \phi_{n,j} ] \leq 2D_j e^{-K n j^2v_n^2}, \quad  \sup_{\lambda \in S_{n,j}(v_n)}
 \E_{\lambda}^{(n)}[ 1_{\Gamma_n} ( 1 -\phi_{n,j}) ] \leq 2e^{-K n j^2v_n^2 },\quad\mbox{if }  j \leq \rho /v_n$$
and
 $$ \E_{\lambda_0}^{(n)}[ 1_{\Gamma_n} \phi_{n,j} ] \leq 2D_j e^{-K n jv_n } , \quad  \sup_{\lambda \in S_{n,j}(v_n)} \E_{\lambda}^{(n)}[ 1_{\Gamma_n} ( 1 -\phi_{n,j}) ] \leq 2e^{-K n jv_n},\quad\mbox{if }  j \geq \rho /v_n,$$
where $K$ is a constant (see Lemma~\ref{lem:tests}). We now bound $D_j$.
First note that for any $\lambda=M_\lambda\bar \lambda$ and $\lambda' = M_{\lambda'}\bar \lambda'$,
\begin{equation}\label{llbar}
\|\lambda -\lambda'\|_1 \leq M_\lambda\|\bar \lambda - \bar \lambda'\|_1+ |M_\lambda-M_{\lambda'}|.
\end{equation}
Assume that $M_\lambda \geq M_{\lambda_0}$. Then, we have
\begin{equation*}
\begin{split}
 \|\lambda - \lambda_0\|_1 & \geq \int_{\bar \lambda >\bar \lambda_0} (M_\lambda \bar \lambda (x) - M_{\lambda_0} \bar \lambda_0 (x))\d x \\
 &= M_\lambda\int_{\bar \lambda >\bar \lambda_0} ( \bar \lambda (x) - \bar \lambda_0 (x))\d x + (M_\lambda - M_{\lambda_0}) \int_{\bar \lambda >\bar \lambda_0}\bar \lambda_0 (x)\d x\\
 &\geq  M_\lambda\int_{\bar \lambda >\bar \lambda_0} ( \bar \lambda (x) - \bar \lambda_0 (x))\d x = \frac{ M_\lambda\|\bar \lambda -\bar \lambda_0\|_1}{ 2}.
 \end{split}
 \end{equation*}
Conversely, if $M_\lambda < M_{\lambda_0}$,
\begin{equation*}
\begin{split}
 \|\lambda - \lambda_0\|_1 & \geq \int_{\bar \lambda_0 >\bar \lambda} (M_{\lambda_0} \bar \lambda_0 (x) - M_\lambda \bar \lambda (x))\d x \\
 &\geq M_{\lambda_0}\int_{\bar \lambda_0 >\bar \lambda} ( \bar \lambda_0 (x) -  \bar \lambda (x))\d x = \frac{ M_{\lambda_0}\|\bar \lambda -\bar \lambda_0\|_1}{ 2 }.
 \end{split}
 \end{equation*}
 So,
$2  \| \lambda - \lambda_0\|_1 \geq  (M_\lambda \vee M_{\lambda_0})  \|\bar \lambda - \bar \lambda_0\|_1$
 and  we finally have
 \begin{equation} \label{norm:mino}
 \| \lambda - \lambda_0\|_1 \geq \max \left\{ \frac{(M_\lambda \vee M_{\lambda_0})  \|\bar \lambda - \bar \lambda_0\|_1 }{ 2 } , \,| M_\lambda - M_{\lambda_0}|  \right\}
 \end{equation}
So, for all $\lambda = M_\lambda \bar \lambda  \in S_{n,j}(v_n)$,
\begin{equation}\label{7.9}
\| \bar \lambda -\bar \lambda_0 \|_1 \leq \frac{2(j+1) v_n}{M_{\lambda_0}}, \quad |M_\lambda-M_{\lambda_0}| \leq (j+1) v_n.
 \end{equation}
 Therefore, $S_{n,j}(v_n)\subset  (\bar S_{n,j} \cap \mathcal F_n)\times \{M:\, \ |M-M_{\lambda_0}| \leq  (j+1) v_n\}$ and any covering of $(\bar S_{n,j} \cap \mathcal F_n)\times \{M:\, \ |M-M_{\lambda_0}| \leq  (j+1) v_n\}$ will give a covering of $S_{n,j}(v_n)$. So, to bound $D_j$, we have to build a convenient covering of $( \bar S_{n,j} \cap \mathcal F_n)\times \{M:\, \ |M-M_{\lambda_0}| \leq  (j+1) v_n\}$. We distinguish two cases.
\begin{itemize}
\item We assume that $(j+1)v_n \leq 2M_{\lambda_0}$. Then, \eqref{7.9} implies that $M_\lambda \leq 3 M_{\lambda_0}$.
 Moreover, if
 $$\|\bar \lambda -\bar  \lambda'\|_1 \leq \frac{\xi j v_n}{3M_{\lambda_0}+1} \quad \mbox{ and } \quad |M_\lambda-M_{\lambda'}|\leq \frac{\xi j v_n}{3M_{\lambda_0}+1},$$
 then, by using  \eqref{llbar},
 $$\| \lambda - \lambda'\|_1\leq \frac{(M_\lambda +1) \xi j v_n}{3M_{\lambda_0}+1} \leq \xi j v_n.$$
By Assumption $(ii)$ of Theorem \ref{th:gene:aalen}, this implies that, for any $\delta>0$, there exists $J_0$ such that for $j\geq J_0$,
  $$D_j\leq D( (3 M_{\lambda_0}+1)^{-1}  \xi j v_n,\, \bar S_{n,j} \cap \mathcal F_n,\, \| \cdot \|_1 ) \times \left(2(j+1)v_n\times\frac{(3 M_{\lambda_0}+1) }{ \xi jv_n}+\frac{1}{2}\right) \lesssim \exp( \delta n(j+1)^2 v_n^2).$$
\item We assume that $(j+1) v_n \geq 2 M_{\lambda_0}$.   If $$ \| \bar \lambda - \bar \lambda' \|_1 \leq \frac{\xi}{4} \quad  \mbox{ and } \quad |M_\lambda-M_{\lambda'}| \leq \frac{\xi (M_\lambda\vee M_{\lambda_0})}{ 4}, $$
   then using again \eqref{llbar} and \eqref{7.9},
  $$
  \| \lambda - \lambda'\|_1\leq \frac{\xi M_\lambda }{4}+\frac{\xi (M_\lambda+M_{\lambda_0})}{4}
  \leq\frac{3\xi M_{\lambda_0}}{4}+\frac{\xi (j+1)v_n}{2}\leq \frac{7\xi (j+1)v_n}{8}\leq \xi jv_n,$$
for $n$ large enough. By Assumption $(i)$ of Theorem \ref{th:gene:aalen}, this implies that, for any $\delta>0$,
  $$D_j \lesssim D( \xi /4,\, \mathcal F_n,\, \| \cdot \|_1 ) \times \log((j+1)v_n)\lesssim \log(jv_n)\exp (\delta n).$$
 \end{itemize}
It is enough to choose $\delta$ small enough to obtain the result of Proposition \ref{prop:test:Aalen}.
\subsection{Proof of Lemma \ref{lem:tests} }\label{sec:proof1}
For any $\lambda$, we denote $\E^{(n)}_{\lambda,\Gamma_n}[\cdot]=\E^{(n)}_{\lambda}[1_{\Gamma_n}\times\cdot]$. We set for any $\lambda,\,\lambda'$,
$$ \|\lambda - \lambda'\|_{\tilde \mu_n}:=\int_\Omega |\lambda(t)-\lambda'(t)|\tilde \mu_n(t)\d t.$$
On $\Gamma_n$ we have
\begin{equation}\label{comp:norm}
m_1 \| \lambda - \lambda_0\|_1 \leq \| \lambda - \lambda_0\|_{\tilde \mu_n} \leq m_2 \| \lambda - \lambda_0\|_1.
\end{equation}
The main tool for building convenient tests is Theorem~3 of \cite{HRR} (and  its proof) applied in the univariate setting. By mimicking the proof of this theorem from Inequality (7.5) to Inequality (7.7), if $H$ is a deterministic function bounded by $b$, we have that, for any $u\geq 0$,
\begin{equation}\label{concentration}
\P_\lambda^{(n)}\left(\left|\int_0^T H_t(\d N_t-\d \Lambda_t)\right|\geq \sqrt{2v u } + \frac{bu}{3} \mbox{ and } \Gamma_n\right)\leq 2e^{-u},
\end{equation}
where we recall that  $ \Lambda_t = \int_0^t Y_s\lambda(s)\d s$ and $v $ is a deterministic constant such that on $\Gamma_n$, almost surely,
$$\int_0^T H_t^2 Y_t\lambda(t)\d t \leq  v .$$
For any non-negative function $\lambda_1$, we define the sets
$$A:=\{t\in\Omega:\ \lambda_1(t)\geq \lambda_0(t)\}, \quad A^c:=\{t\in\Omega:\ \lambda_1(t)< \lambda_0(t)\}$$ and the following pseudo-metrics
$$d_A(\lambda_1,\,\lambda_0):=\int_A[\lambda_1(t)-\lambda_0(t)]\tilde\mu_n(t)\d t \quad d_{A^c}(\lambda_1,\,\lambda_0):=\int_{A^c}[\lambda_0(t)-\lambda_1(t)]\tilde\mu_n(t)\d t.$$
Note that
$ \| \lambda_1 - \lambda_0 \|_{\tilde\mu_n}  = d_A(\lambda_1,\,\lambda_0) +  d_{A^c}(\lambda_1,\,\lambda_0).$
For $u>0$, if $d_A(\lambda_1, \lambda_0) \geq d_{A^c}(\lambda_1, \lambda_0)$, define the test
\[\phi_{\lambda_1,A}(u):=1\left\{N(A)-\int_{A} \lambda_0(t)Y_t\d t \geq  \rho_n(u) \right\} , \mbox{ with  } \rho_n(u) := \sqrt{2n v(\lambda_0) u} + \frac{u}{3},\]
where, for any non-negative function $\lambda$,
\begin{equation}\label{vl}
v(\lambda) := (1+\alpha)\int_\Omega \lambda(t) \tilde\mu_n(t)\d t.
\end{equation}
Similarly, if
$d_A(\lambda_1,\, \lambda_0) < d_{A^c}(\lambda_1,\, \lambda_0)$, define
\[\phi_{\lambda_1,A^c}(u):=1\left\{N(A^c)-\int_{A^c}\lambda_0(t)Y_t\d t \leq - \rho_n(u) \right\}. \]
Since for any non-negative function $\lambda$,  on $\Gamma_n$, by using \eqref{ass:Y3},
\begin{equation}
\label{Yn:mun}
(1-\alpha) \int_\Omega \lambda(t) \tilde\mu_n(t)\d t\leq  \int_\Omega \lambda(t)\frac{Y_t}{n} \d t \leq (1+\alpha)\int_\Omega\lambda(t) \tilde\mu_n(t)\d t,
\end{equation}
then inequality \eqref{concentration} applied with $H=1_A$ or $H=1_{A^c}$,  $b=1$ and $v=n v(\lambda_0)$ implies that for any $u>0$,
\begin{equation}\label{err1}
\E^{(n)}_{\lambda_0,\Gamma_n}[\phi_{\lambda_1,A}(u) ]\leq 2e^{-u },\quad \E^{(n)}_{\lambda_0,\Gamma_n}[\phi_{\lambda_1,A^c}(u) ]\leq 2e^{-u }.
\end{equation}
Now, we state the following lemma whose proof is given in Section~\ref{sec:proof2}.
\begin{Lemma} \label{lem:tests:aalen:error2}
Let $\lambda$ be a non-negative function. Assume that $$\| \lambda - \lambda_1\|_{\tilde\mu_n} \leq \frac{1-\alpha}{ 4(1+\alpha)} \| \lambda_1 -\lambda_0\|_{\tilde\mu_n}.$$
 We set $\tilde M_n(\lambda_0)=\int_\Omega\lambda_0(t)\tilde\mu_n(t)\d t$
 and we distinguish two cases.
\begin{enumerate}
\item Assume that $d_A(\lambda_1,\, \lambda_0) \geq d_{A^c}(\lambda_1,\, \lambda_0)$. Then,
$$\E^{(n)}_{\lambda,\Gamma_n}[1-\phi_{\lambda_1,A}(u_A)]\leq 2\exp(-u_A),$$
where
$$
u_A=\left\{\begin{array}{lc}
u_{0A}n d_A^2(\lambda_1,\,\lambda_0),&\mbox{ if }  \| \lambda_1 -\lambda_0\|_{\tilde\mu_n}\leq 2\tilde M_n(\lambda_0),\\[2pt]
u_{1A}nd_A(\lambda_1,\lambda_0),&\mbox{ if }  \| \lambda_1 -\lambda_0\|_{\tilde\mu_n}> 2\tilde M_n(\lambda_0),
\end{array}\right.
$$
and $u_{0A}$ and $u_{1A}$ are two constants only depending on $\alpha$, $M_{\lambda_0}$, $m_1$ and $m_2$.
\item Assume that $d_A(\lambda_1,\, \lambda_0) < d_{A^c}(\lambda_1,\, \lambda_0)$. Then,
$$\E^{(n)}_{\lambda,\Gamma_n}[1-\phi_{\lambda_1,A^c}(u_{A^c})]\leq 2\exp(-u_{A^c}),$$
where
$$
u_{A^c}=\left\{\begin{array}{lc}
u_{0A^c}nd_{A^c}^2(\lambda_1,\,\lambda_0),&\mbox{ if }  \| \lambda_1 -\lambda_0\|_{\tilde\mu_n}\leq 2\tilde M_n(\lambda_0),\\[2pt]
u_{1A^c}nd_{A^c}(\lambda_1,\,\lambda_0),&\mbox{ if }  \| \lambda_1 -\lambda_0\|_{\tilde\mu_n}> 2\tilde M_n(\lambda_0),
\end{array}\right.
$$
and $u_{0A^c}$ and $u_{1A^c}$ are two constants only depending on $\alpha$, $M_{\lambda_0}$, $m_1$ and $m_2$.
\end{enumerate}
\end{Lemma}
Note that, by using \eqref{comp:norm}, if $d_A(\lambda_1,\, \lambda_0)\geq d_{A^c}(\lambda_1,\, \lambda_0)$,
\begin{eqnarray*}
u_A&\geq&\min\{u_{0A}n d_A^2(\lambda_1,\,\lambda_0),\,u_{1A}nd_A(\lambda_1,\,\lambda_0)\}\\
&\geq&nd_A(\lambda_1,\,\lambda_0)\times\min\{u_{0A} d_A(\lambda_1,\,\lambda_0),\,u_{1A}\}\\
&\geq&\frac{nm_1\|\lambda_1-\lambda_0\|_1}{2}\times\min\left\{\frac{u_{0A}m_1\|\lambda_1-\lambda_0\|_1}{2},\,u_{1A}\right\}\\
&\geq& K_A n\| \lambda_1 - \lambda_0 \|_{1} \times\min\{\| \lambda_1 - \lambda_0 \|_1,\,m_1\},
\end{eqnarray*}
for $K_A$ a positive constant small enough only depending on $\alpha$, $M_{\lambda_0}$, $m_1$ and $m_2$.
Similarly, if $d_A(\lambda_1,\, \lambda_0) < d_{A^c}(\lambda_1,\, \lambda_0)$,
\begin{eqnarray*}
u_{A^c}&\geq&\frac{nm_1\|\lambda_1-\lambda_0\|_1}{2}\times\min\left\{\frac{u_{0A^c}m_1\|\lambda_1-\lambda_0\|_1}{2},\,u_{1A^c}\right\}\\
&\geq&K_{A^c} n\| \lambda_1 - \lambda_0 \|_{1} \times\min\{\| \lambda_1 - \lambda_0 \|_1,\,m_1\},
\end{eqnarray*}
for $K_{A^c}$  a positive constant small enough only depending on $\alpha$, $M_{\lambda_0}$, $m_1$ and $m_2$.
Now, we set
$$\phi_{\lambda_1}=\phi_{\lambda_1,A}(u_A)1_{\left\{d_A(\lambda_1,\, \lambda_0)\geq d_{A^c}(\lambda_1,\, \lambda_0)\right\}}+\phi_{\lambda_1,A^c}(u_{A^c})1_{\left\{d_A(\lambda_1,\, \lambda_0)< d_{A^c}(\lambda_1,\, \lambda_0)\right\}}$$
so that, with
$K=\min\{K_A,\,K_{A^c}\}$, by using \eqref{err1},
\begin{eqnarray*}
\E^{(n)}_{\lambda_0,\Gamma_n}[\phi_{\lambda_1}]&=&\E^{(n)}_{\lambda_0,\Gamma_n}[\phi_{\lambda_1,A}(u_A)]1_{\left\{d_A(\lambda_1,\, \lambda_0)\geq d_{A^c}(\lambda_1,\, \lambda_0)\right\}}+\E^{(n)}_{\lambda_0,\Gamma_n}[\phi_{\lambda_1,A^c}(u_{A^c})]1_{\left\{d_A(\lambda_1,\, \lambda_0)< d_{A^c}(\lambda_1,\, \lambda_0)\right\}}\\
&\leq&2e^{-u_A}1_{\left\{d_A(\lambda_1,\, \lambda_0)\geq d_{A^c}(\lambda_1,\, \lambda_0)\right\}}+2e^{-u_{A^c}}1_{\left\{d_A(\lambda_1, \,\lambda_0)< d_{A^c}(\lambda_1,\, \lambda_0)\right\}}\\
&\leq&2\exp\left(- K n\| \lambda_1 - \lambda_0 \|_{1} \times\min\{\| \lambda_1 - \lambda_0 \|_1,\,m_1\} \right).
\end{eqnarray*}
If $\|\lambda - \lambda_1 \|_{1} <\xi\| \lambda_1 - \lambda_0 \|_1$, $ \xi=m_1(1-\alpha)/( 4m_2(1+\alpha)) $,
then
$$\| \lambda - \lambda_1\|_{\tilde\mu_n} \leq \frac{1-\alpha}{ 4(1+\alpha)} \| \lambda_1 -\lambda_0\|_{\tilde\mu_n}$$
and Lemma \ref{lem:tests:aalen:error2} shows that
\begin{eqnarray*}
\E^{(n)}_{\lambda,\Gamma_n}[1-\phi_{\lambda_1}]&\leq&2e^{-u_A}1_{\left\{d_A(\lambda_1,\, \lambda_0)\geq d_{A^c}(\lambda_1, \, \lambda_0)\right\}}+2e^{-u_{A^c}}1_{\left\{d_A(\lambda_1,\, \lambda_0)< d_{A^c}(\lambda_1,\, \lambda_0)\right\}}\\
&\leq&2\exp\left(- K n\| \lambda_1 - \lambda_0 \|_{1} \times\min\{\| \lambda_1 - \lambda_0 \|_1,\,m_1\} \right),
\end{eqnarray*}
which ends the proof of Lemma \ref{lem:tests}.
\subsection{Proof of  Lemma \ref{lem:tests:aalen:error2} }\label{sec:proof2}
We only consider the case $d_A(\lambda_1,\, \lambda_0)\geq d_{A^c}(\lambda_1,\, \lambda_0)$.
The case $d_A(\lambda_1,\, \lambda_0)< d_{A^c}(\lambda_1,\, \lambda_0)$ can be dealt with by using similar arguments.
So, we assume that $d_A(\lambda_1,\, \lambda_0) \geq d_{A^c}(\lambda_1, \,\lambda_0)$. On $\Gamma_n$ we have
\begin{eqnarray*}
\int_A(\lambda_1(t)-\lambda_0(t))Y_t\d t&\geq&n(1-\alpha)\int_A(\lambda_1(t)-\lambda_0(t))\tilde\mu_n(t)\d t\\
&\geq&\frac{n(1-\alpha)}{2}\|\lambda_1-\lambda_0\|_{\tilde\mu_n} \geq2n(1+\alpha)\|\lambda-\lambda_1\|_{\tilde\mu_n}\\
&\geq&2n(1+\alpha)\int_A|\lambda(t)-\lambda_1(t)|\tilde\mu_n(t)\d t
\geq 2\int_A|\lambda(t)-\lambda_1(t)|Y_t\d t.
\end{eqnarray*}
Therefore, we have
\begin{eqnarray*}
\E^{(n)}_{\lambda,\Gamma_n}[1- \phi_{\lambda_1,A}(u_A)] &=& \P^{(n)}_{\lambda,\Gamma_n} \left\{N(A)- \int_{A} \lambda(t)Y_t\d t <  \rho_n(u_A)+ \int_{A} (\lambda_0-\lambda)(t)Y_t\d t  \right\}  \\
&=& \P^{(n)}_{\lambda,\Gamma_n}  \left\{ N(A)- \int_{A} \lambda(t)Y_t\d t <  \rho_n(u_A)- \int_{A} (\lambda_1-\lambda_0)(t)Y_t\d t\right. \\
&&\hspace{7cm}+\left.\int_{A} (\lambda_1-\lambda)(t)Y_t\d t\right\} \\
&\leq& \P^{(n)}_{\lambda,\Gamma_n}  \left\{ N(A)- \int_{A} \lambda(t)Y_t\d t <  \rho_n(u_A)- \frac{1}{2}\int_{A} (\lambda_1-\lambda_0)(t)Y_t\d t\right\}.
\end{eqnarray*}
\noindent\underline{Assume $\|\lambda_1-\lambda_0\|_{\tilde\mu_n}\leq 2\tilde M_n(\lambda_0)$}\\

\noindent This assumption implies that
$$d_A(\lambda_1,\,\lambda_0)\leq\|\lambda_1-\lambda_0\|_{\tilde\mu_n}\leq 2\tilde M_n(\lambda_0)\leq 2m_2M_{\lambda_0}.$$
Since
$v(\lambda_0)= (1+\alpha)\tilde M_n(\lambda_0),$ with
$u_A=u_{0A}n d_A^2(\lambda_1,\, \lambda_0),$
where $u_{0A}\leq 1$ is a constant depending on $\alpha$, $m_1$ and $m_2$ chosen later, we  have
$$ \rho_n(u_A) \leq n d_A(\lambda_1,\, \lambda_0)\sqrt{2u_{0A}(1+\alpha)\tilde M_n(\lambda_0)}+ \frac{ u_{0A} n d_A^2(\lambda_1,\, \lambda_0)}{ 3 }  \leq \sqrt{u_{0A}}  n d_A(\lambda_1, \,\lambda_0) K_1 $$
as soon as $K_1\geq\sqrt{2(1+\alpha)\tilde M_n(\lambda_0)}+ \frac{2\tilde M_n(\lambda_0)\sqrt{ u_{0A}}}{3}$. Observe that
the definition of $v(\lambda)$ in \eqref{vl} gives
\begin{eqnarray*}
v(\lambda) &=&(1+\alpha)\int_\Omega\lambda_0(t)\tilde\mu_n(t)\d t+ (1+\alpha)\int_\Omega(\lambda(t)-\lambda_0(t))\tilde\mu_n(t)\d t\\
&\leq& v(\lambda_0)+ (1+\alpha)\|\lambda -\lambda_0\|_{\tilde\mu_n}
\leq  v(\lambda_0) + (1+\alpha) \left( \|\lambda -\lambda_1\|_{\tilde\mu_n} + \|\lambda_1 -\lambda_0\|_{\tilde\mu_n}\right) \\
&\leq& v(\lambda_0) + \frac{5+3\alpha}{4}\|\lambda_1 -\lambda_0\|_{\tilde\mu_n} \leq C_1,
\end{eqnarray*}
where $C_1$ only depends on $\alpha,$ $M_{\lambda_0},$ $m_1$ and $m_2$.
Combined with \eqref{Yn:mun},  this implies that, on $\Gamma_n$, if $K_1\leq \frac{1-\alpha}{4\sqrt{u_{0A}}}$, which is true for $u_{0A}$ small enough,
\begin{equation*}
\begin{split}
 \frac{1}{2}\int_{A} (\lambda_1-\lambda_0)(t)Y_t\d t-  \rho_n(u_A)
&\geq \frac{(1-\alpha)n}{ 2 } d_A(\lambda_1,\, \lambda_0) \left[ 1 - \frac{2K_1\sqrt{u_{0A}}}{1-\alpha}\right] \\
&\geq
 \frac{ (1-\alpha)n}{4 } d_A(\lambda_1,\, \lambda_0) \geq \sqrt{ 2 nC_1r} + \frac{r}{3}\geq \sqrt{ 2 nv(\lambda) r} + \frac{r}{3},
\end{split}
\end{equation*}
with
$$ r = n\min \left\{
\frac{(1-\alpha)^2 d_A^2(\lambda_1,\, \lambda_0) }{ 128 C_1 }  , \,\frac{ 3 (1-\alpha) d_A(\lambda_1,\, \lambda_0)}{ 8} \right\}.$$
Inequality \eqref{concentration} then leads to
\begin{equation}\label{r}
\E^{(n)}_{\lambda,\Gamma_n}[1- \phi_{\lambda_1,A}(u_A) ] \leq 2e^{-r}.
\end{equation}
For $u_{0A}$ small enough only depending on $M_{\lambda_0}$, $\alpha$, $m_1$ and $m_2$, we have
$$\frac{(1-\alpha)}{4\sqrt{u_{0A}}}\geq \sqrt{2(1+\alpha)\tilde M_n(\lambda_0)}+ \frac{2\tilde M_n(\lambda_0)\sqrt{u_{0A}}}{3}$$
so (\ref{r}) is true. Since $r\geq u_A$ for $u_{0A}$ small enough, then
$$
\E^{(n)}_{\lambda,\Gamma_n}[1- \phi_{\lambda_1,A}(u) ] \leq 2e^{-u_A}.
$$

\noindent\underline{Assume $\|\lambda_1-\lambda_0\|_{\tilde\mu_n}\geq 2\tilde M_n(\lambda_0)$}\\

\noindent We take
$u_A=u_{1A}nd_A(\lambda_1,\,\lambda_0),$
where $u_{1A}\leq 1$ is a constant depending on $\alpha$ chosen later. We still consider the same test $\phi_{\lambda_1, A}(u_A)$. Observe now that, since $d_A(\lambda_1,\,\lambda_0)\geq \frac{1}{2}\|\lambda_1 -\lambda_0\|_{\tilde\mu_n}\geq \tilde M_n(\lambda_0)$,
\begin{eqnarray*}
\rho_n(u_A)&=&\sqrt{2nu_Av(\lambda_0)}+\frac{u_A}{3}\\
&\leq&n\sqrt{2(1+\alpha)u_{1A}\tilde M_n(\lambda_0)d_A(\lambda_1,\,\lambda_0)}+\frac{d_A(\lambda_1,\,\lambda_0)nu_{1A}}{3}\\
&\leq&\left(\sqrt{2(1+\alpha)}+\frac{1}{3}\right)\sqrt{u_{1A}}nd_A(\lambda_1,\,\lambda_0)
\end{eqnarray*}
and, under the assumptions of the lemma,
\begin{equation}\label{v}
v(\lambda)\leq(1+\alpha)\tilde M_n(\lambda_0) + (1+\alpha) \left( \|\lambda -\lambda_1\|_{\tilde\mu_n} + \|\lambda_1 -\lambda_0\|_{\tilde\mu_n}\right)\leq C_2d_A(\lambda_1,\,\lambda_0),
\end{equation}
where $C_2$ only depends on $\alpha.$ Therefore,
\begin{eqnarray*}
\frac{1}{2}\int_{A} (\lambda_1-\lambda_0)(t)Y_t\d t-\rho_n(u_A)&\geq&\frac{n(1-\alpha)}{2}\int_A(\lambda_1(t)-\lambda_0(t))\tilde\mu_n(t)\d t -\left(\sqrt{2(1+\alpha)}+\frac{1}{3}\right)\sqrt{u_{1A}}nd_A(\lambda_1,\,\lambda_0)\\
&\geq&\left(\frac{1-\alpha}{2}-\left(\sqrt{2(1+\alpha)}+\frac{1}{3}\right)\sqrt{u_{1A}}\right)nd_A(\lambda_1,\,\lambda_0)\\
&\geq&\frac{1-\alpha}{4}nd_A(\lambda_1,\,\lambda_0),
\end{eqnarray*}
where the last inequality is true for $u_{1A}$ small enough depending only on $\alpha$. Finally, using (\ref{v}) and since $u_A=u_{1A}nd_A(\lambda_1,\,\lambda_0)$, we have
\begin{eqnarray*}
\frac{1-\alpha}{4}nd_A(\lambda_1,\,\lambda_0)&\geq&\sqrt{2nC_2d_A(\lambda_1,\,\lambda_0)u_{1A}nd_A(\lambda_1,\,\lambda_0)}+\frac{u_{1A}nd_A(\lambda_1,\,\lambda_0)}{3}\\
&\geq&\sqrt{2nv(\lambda)u_A}+\frac{u_A}{3}
\end{eqnarray*}
for $u_{1A}$ small enough depending only on $\alpha$. We then obtain
$$\E^{(n)}_{\lambda,\Gamma_n}[1- \phi_{\lambda_1,A}(u_A)] \leq 2e^{-u_A},$$
which ends the proof of the lemma.

\section{Details of the Gibbs algorithm} \label{sec:details algo}

We detail the algorithm used to sample from the posterior distribution $(\overline{\lambda}, A, \gamma)| \Nb$ in the Poisson process context, in the most complete case (with a hierarchical level on $\gamma$).  In case where $\gamma$ is set to a fixed value, then the corresponding part in the algorithm is removed. 

\noindent As a standard Gibbs algorithm, the simulation is decomposed into  three steps: 
$$[1] \quad  \overline{\lambda} | A, \gamma, \Nb \quad \quad [2] \quad    A | \overline{\lambda}, \gamma, \Nb   \quad\quad [3]  \quad \gamma | A, \overline{\lambda}, \Nb. $$
where $N$ is the observed Poisson process over $[0,T]$,  namely a number of jumps $N(T)$ and jump instants $(T_1,\dots, T_{N(T)})$. 
\noindent In order to avoid an artificial truncation in $\overline{\lambda}$, we use the slice sampler strategy proposed by \cite{fall2012}.  
\noindent More precisely, we consider  the stick breaking representation of $\overline{\lambda}$. Let $c_i$ be the affectation variable of data $W_i$. The DPM model is written as: 
$$W_i  | c_i, \btheta^* \sim  h_{\theta^{\star}_{c_i}}, \quad  P(c_i = k) = w_k,  \forall k \in \mathbb{N}^{*}\quad   (w_k)_{k \in \mathbb{N}^{\star}} \sim  Stick(A),  \quad (\theta^{*}_{k})_{k \in \mathbb{N}^{\star}} \sim_{i.i.d} G_{\gamma}.$$ 
The  slice sampler strategy consists in introducing a latent variable $u_i$ such that the joint distribution of $(W_i,u_i)$ is 
$
 p(W_i,u_i|\bw,\btheta^{*}) = \sum_{k=1}^\infty w_k h_{\theta^{*}_{k}}(W_i) \frac{1}{\xi_k}\ind_{[0,\xi_k]}(u_i)$ with  $\xi_k = \min(w_k,\zeta)$, which  can be reformulated as:  
 \begin{equation}\label{eq:u x jointe2}
 p(W_i,u_i |\bw,\btheta^*)=  \frac{1}{\zeta}\ind_{[0,\zeta]}(u_i)  \sum_{k=1,w_k >\zeta}^\infty w_k h_{\theta^{*}_{k}}(W_i) +   \sum_{k=1,  u_i \leq w_k \leq \zeta}^\infty h_{\theta^{*}_{k}}(W_i)\ind_{[0,w_k]}(u_i)
 \end{equation}
 $(w_k)_{k \geq 1} $ verifying  $\sum_{k \geq 1} w_k=1$  (implying $\lim_{k \rightarrow \infty} w_k =0$), the cardinal of $\{k, w_k > \varepsilon\}$ is finite for   every $\varepsilon >0$,  and  the sum in  (\ref{eq:u x jointe2}) is finite.  
 
 \begin{remark}
 Note that the number of non-null terms in (\ref{eq:u x jointe2})  is rigorously  dependent of  the observation index $i$ (denoted $K_i^{*}$). But in the algorithm we will deal with the maximum of the $K^{*}_i$:  $$K^{*} = \max\{K_i^{*},i=1\dots N(T)\}$$ 
 \end{remark}


\noindent The Gibbs algorithm with the Slice sampler strategy now takes into account the latent variable $\bu=(u_1,\dots,u_{N(T)})$ which is sampled conjointly with $\overline{\lambda}$, resulting into the following steps: 
$$[1^{\star}] \quad  \overline{\lambda},\bu | A, \gamma, \Nb \quad \quad [2^{\star}] \quad    A | \overline{\lambda}, \bu,\gamma, \Nb   \quad\quad [3^{\star}]  \quad \gamma | A, \overline{\lambda}, \bu,\Nb. $$
We now detail steps $[1^{\star}]$ , $[2^{\star}]$ and $[3^{\star}]$.

\subsection{Details of the algorithm}

\noindent \emph{Initialisation} The Gibbs algorithms are initialized on $ A^{(0)}=10$ ($0$ referring to the iteration number of the Gibbs algorithm). We set: 
$K^{(0)}=  N_T$. for $k=1\dots K$, $(\theta^*)^{(0)}_k \sim G_{\gamma}$. 
 When the hierarchical approach is considered on $\gamma$ ($\gamma \sim \Gamma(a_\gamma,b_\gamma)$), we initialize $\gamma$ on its prior mean value: $\gamma^{(0)}= \frac{a_{\gamma}}{b_{\gamma}}$.

\paragraph{$[1^{\star}]$ Sampling from  $ \overline{\lambda},\bu | A, \gamma, \Nb$}$\;$ \\
 
\noindent Step $[1^{\star}]$ of the Gibbs algorithm is  decomposed into $5$ steps which are detailed below. Let  $\bu = (u_1,\dots,u_{N(T)})$, $\bc=(c_1,\dots,c_{N(T)})$, $\btheta^* = (\theta_1,\dots,\theta_{K^{\star}})$, $\bw = (w_1,\dots,w_{K^{\star}})$ --$\bc$, $\btheta^*$ and $\bw$ representing $\overline{\lambda}$-- be the current object. We denote by $K_{N(T)}$ the number of non-empty classes: 
$$ K_{N(T)} = \# \{k\in \{1\dots K^{\star}\} | \exists j \in \{1,\dots,N(T)\} \mbox{ such that }  c_j = k\} $$ 
$\btheta^*$ and $\bw$ are ordered such that the elements indexed from $K_{N(T)}+1$ to $K^{\star}$ correspond to empty classes. 
$\bu$, $\bc$, $\btheta^*$ and $\bw$ are iteratively sampled as follows:

 \begin{enumerate}
\item[[$1^{\star}$.a]]  First we sample $\bu |  \bw,\btheta^*,\bc,\Nb,A, \gamma$ using the following identities: 
\begin{eqnarray*}
p(\bu | \bw,\btheta^*,\bc,\Nb,A, \gamma) &\propto& p(\bu,\Nb |  \bw,\btheta^*,\bc) =  p(\bu,\Nb |   \btheta^*,\bc) = \prod_{i=1}^{N(T)} p(u_i,W_i | c_i,\theta_{c_i}^{\star}) \\
&=&   \prod_{i=1}^{N(T)} h_{\theta^{\star}_{c_i}}(W_i) \frac{1}{\xi_{c_i}}\ind_{[0,\xi_{c_i}]}(u_i) \propto \prod_{i=1}^{N(T)}   \frac{1}{\xi_{c_i}}\ind_{[0,\xi_{c_i}]}(u_i) \\
 \end{eqnarray*}
 where $\xi_{c_i} =\min(w_{c_i},\zeta)$.  So for every $i=1\dots N(T)$, $u_i \sim \mathcal{U}_{[0,\min(w_{c_i},\zeta)]}$.

\item[[$1^{\star}$.b]]  Secondly, we sample the weights of the empty classes $(w_{k})_{k\geq K_{N(T)}+1} | \Nb,\bw,\bc,\btheta^*,A, \gamma$. Without the slice sampler strategy, there is an infinite number to sample. But, thanks to the slice sampling, we only need to sample a finite number $K^{\star}$.  The weights of the empty classes follow the prior distribution (stick breaking).  For $k>K_{N(T)}$, 
\begin{eqnarray*}
v_k &\sim& \mathcal{B}(1,A)\\
w_k &=& v_k r_{k-1}\\
r_k &=& r_{k-1} (1-v_k) 
\end{eqnarray*}
As explained in \cite{fall2012}, we know that we have to represent all the components such that their weights $w_k > u_i$ for all the $u_i$. Assume that we have sampled $w_1,\dots, w_k$, then the weights of the following components can not exceed the rest $1- \sum_{j=1}^k w_j = r_{k}$. So if $r_k$ is such that $r_k < u_i$, for all $i=1\dots n$, i.e. if $r_k < u^{\star}$ with $u^{\star}= \min\{u_1,\dots,u_{N(T)}\}$ we are sure that there is no ``interesting component" after that,  ``interesting" meaning that they won't appear in joint the distribution of $(\Nb,\bu)$. We can stop and get
$$ K^{\star} = \min \{k, r_k < u^{\star}\}.$$
 In the end we have sampled $(w_{K_{N(T)}+1},\dots,w_{K^{\star}})$. 

\item[[$1^{\star}$.c]]   Sample the parameters of the empty classes, $(\theta^*_{K_{N(T)}+1},\dots,\theta^*_{K^{\star}})$
 $$\forall k=K_{N(T)}+1, \dots, K^{\star}, \quad \theta^*_{k} \sim_{i.i.d} G_{\gamma}  $$

\item[[$1^{\star}$.d]]  Sample the index $\bc = (c_1,\dots,c_{N(T)}) | \Nb,\bu,\btheta^*,\bw$, i.e. affect the observations to the classes  $\{1, \dots,K^{\star}\}$. Note that we will get new empty classes. The affectations are done using the following probabilities: 
\begin{eqnarray*}
p(\bc | \btheta^*,\bw,\bu,\Nb) &\propto&p(\bc, \btheta^*,\bw,\bu,\Nb)=  p(\Nb,\bu, \bc|  \btheta^*,\bw) p(\bw) p(\btheta^*) \\
&\propto& \prod_{i=1}^{ N(T)}  p(W_i,u_i, c_i|  \btheta^*,\bw)
\end{eqnarray*}
So, the $c_i$ are independent and,  for $i=1\dots N(T)$, for $k=1\dots, K^{\star}$
\begin{eqnarray*}
p(c_i = k |  \btheta^*,\bw,u_i,W_i) = w_{i,k} \propto  h_{\theta^{*}_{k}}(W_i) \frac{w_k}{\min(\zeta,w_{k}) }\ind_{ \{k | u_i< \min(\zeta,w_{k})\}}(k) 
\end{eqnarray*}

We obtain a new $K_{N(T)}$, which is the number of non-empty classes. We re-arrange the weights and the parameters by order of appearance in this  affectation. 

\item[[$1^{\star}$.e]]   Update $(w_1,\dots,w_{K_{N(T)}})$ and  $(\theta_1,\dots,\theta_{K_{N(T)}})$ for the non-empty classes. 

\begin{eqnarray}\label{loi theta}
 p(\theta_k | \bu,\Nb,\bw,\bc) &\propto& G_{\gamma}(\theta_k) \prod_{i=1, c_i=k}^n h_{\theta_k}(W_i), \quad \forall k=1\dots K_{N(T)}\\
w_1,\dots,w_{K_{N(T)}},r_{K_{N(T)}} &\sim& \mbox{Dir}(n_1,\dots,n_{K_{N(T)}},A) \nonumber
\end{eqnarray}
where $n_k = \#\{i | c_i=k\}$ and $r_{K_{N(T)}} = 1-\sum_{k=1}^{K_{N(T)}} w_k$.

Note that when $G_{\gamma}$ is the inverse of the translated inverse exponential distribution, $ p(\theta_k | \bu,\Nb,\bw,\bc)$ given in equation (\ref{loi theta}) is:
\begin{equation}\label{loi theta2}
p(\theta_k | \bu,\Nb,\bw,\bc) \propto  \frac{1}{\left(\frac{1}{\theta_k}-\frac{1}{T}\right)^{a-1}} e^{-\frac{\gamma}{\frac{1}{\theta_k}-\frac{1}{T}}}\frac{1}{\theta_k^{n_k}}\ind_{\left[\frac{1}{\max_{i | c_i=k} W_i} , + \infty \right[}(\theta_k)
\end{equation}
Its simulation can not be performed directly and we  resort to an accept-reject procedure to simulate exactly under this distribution. 
\end{enumerate}

\begin{remark}
The accept-reject procedure we propose is detailed and discussed here after. Note that $(\theta_k)$ changes of dimension at each iteration of the algorithm. As a consequence, a Metropolis-Hastings procedure can not be considered easily, since it could jeopardize the theoretical and practical convergence properties of the algorithm. 
\end{remark}

\paragraph{[$2^{\star}$.]  Sampling from  $A | \bw,\btheta^*,\bc, \bu, \gamma, \Nb$ }$\;$ \\


\noindent    Let  $K_{N(T)}$ be the current number of non-empty classes. \cite{west92} proves that under the prior distribution $A \sim \Gamma(a_\alpha, b_\alpha)$, 
we have 
\begin{equation}
A | x,K_{N(T)} \sim \pi_x \Gamma(a_{A}+K_{N(T)},b_{A}-\log(x))+(1-  \pi_x) \Gamma(a_{A}+K_{N(T)}-1,b_{A}-\log(x))
\end{equation}

where 
\begin{eqnarray*}
x| A,K _{N(T)} &\sim& \mathcal{B}(A+1,N(T))\\
\frac{\pi_x}{1-\pi_x} &=& \frac{a_{A}+K_{N(T)}-1}{n(b_A-\log(x))}
\end{eqnarray*}
Note that the generation of the new value of $A$ relies on the current value of $A$.

\paragraph{$[3^{\star}.]$  Sampling from  $\gamma | \bw,\btheta^*,\bc, \bu,A, \Nb$}  If a hierarchical level is set on $\gamma$,  we have to sample 
$$ \gamma | \btheta^*, \bw,\bu,\Nb,\bc$$
Using  the previous conditional distributions, we have:  
\begin{eqnarray*}
p( \gamma | \btheta^*, \bw,\bu,\Nb,\bc ) &\propto& p(\bu,\Nb,\btheta^*,\bw,\bc|\gamma)\pi(\gamma) =   p(\bu,\Nb| \btheta^*, \bc )p(\bc | \bw) p(\btheta^*|\gamma) \pi(\gamma) =\pi(\gamma) \prod_{k=1}^{K^{\star}} G_{\gamma}(\theta_k)\\
&\propto& \gamma^{a_\gamma-1} e^{-b_\gamma \gamma}   \prod_{k=1}^{K^{\star}}  \gamma^a e^{-\gamma /(\frac{1}{\theta^{\star}_k}-\frac{1}{T})} =   \gamma^{a_\gamma +  a K^{\star}-1} e^{-\gamma\left(b_\gamma + \sum_{k=1}^{K^{\star}} \frac{1}{(\frac{1}{\theta^{\star}_k} - \frac{1}{T})} \right)} 
\end{eqnarray*}
where $K^{\star}$ is the total number of classes used to represent $\overline{\lambda}$. Finally, we get: 

\begin{equation}\label{eq:postgamma}
\gamma | \btheta^*, \bw,\bu,\Nb,\bc \sim \Gamma\left(   a_\gamma +  a K^{\star}, b_\gamma + \sum_{k=1}^{K^{\star}} \frac{1}{(\frac{1}{\theta^{\star}_k} - \frac{1}{T})} \right)
\end{equation}

\bibliographystyle{apalike}
\bibliography{biblio.bib}
\end{document}